\documentclass[10pt]{amsart}
\usepackage[hidelinks]{hyperref}
\usepackage{doi,url}

\usepackage[utf8]{inputenc}
\usepackage{amsmath}
\usepackage{amssymb}
\usepackage{mathtools}
\usepackage[numbers]{natbib}
\DeclareRobustCommand{\noopsort}[1]{}

\usepackage{prettyref}
\usepackage{enumitem}

\numberwithin{equation}{section}

\newtheorem{theorem}{Theorem}[section]

\newtheorem{proposition}[theorem]{Proposition}
\newtheorem{corollary}[theorem]{Corollary}

\newtheorem{lemma}[theorem]{Lemma}
\newtheorem{fact}[theorem]{Fact}
\newtheorem{question}[theorem]{Question}

\newcounter{introthmcounter}
\newtheorem{thmA}{Theorem}[introthmcounter]

\theoremstyle{definition}

\newtheorem{definition}[theorem]{Definition}

\newtheorem{context}[theorem]{Context}
\newtheorem*{addendum*}{Addendum}

\theoremstyle{remark}
\newtheorem{remark}[theorem]{Remark}

\newcommand{\bP}{\mathbb P}
\newcommand{\bE}{\mathbb E}

\newcommand{\bR}{\mathbb R}
\newcommand{\bN}{\mathbb N}
\newcommand{\bK}{\mathbb K}

\DeclareMathOperator{\sign}{\mathrm{sgn}} 

\DeclareMathOperator{\std}{\mathrm{st}}

\newcommand{\bU}{\mathbb{U}}

\newcommand{\cN}{\mathcal{N}}

\def\Lhp{(\!(} 
\def\Rhp{)\!)} 

\def\-{\text{-}}

\newcommand{\No}{\mathbf{No}}

\newcommand{\tp}{\mathrm{tp}}

\newcommand{\dcl}{\mathrm{dcl}}
\newcommand{\convex}{\mathrm{convex}}

\newcommand{\tame}{\mathrm{tame}}
\newcommand{\fM}{\mathfrak{M}}

\newcommand{\cO}{\mathcal{O}}
\newcommand{\co}{\mathfrak{o}}

\newcommand{\var}[1]{\mathrm{#1}}

\newcommand{\RCF}{\mathrm{RCF}}

\DeclareMathOperator{\Br}{\mathrm{Br}}

\DeclareMathOperator{\val}{\mathbf{v}}
\DeclareMathOperator{\res}{\mathbf{r}}
\DeclareMathOperator{\rv}{\mathbf{rv}}
\DeclareMathOperator{\CH}{\mathrm{CH}}

\DeclareMathOperator{\Exponents}{\mathrm{Exponents}}
\DeclareMathOperator{\Span}{\mathrm{Span}}

\title[{$T$}-Convexity and Weakly Immediate types]{{$T$}-convexity, Weakly Immediate Types, and {$T$}-\texorpdfstring{$\lambda$}{λ}-Spherical Completions of o-minimal Structures}
\author{Pietro Freni}

\address{Institute of Mathematics, Czech Academy of Sciences, Žitná 25, 110 00 Praha 1, Czech Republic}



\thanks{This work has been supported by the Czech Science Foundation grant 26-22545I, by the Czech Academy of Sciences CAS (RVO 67985840), and by a Scholarship of the School of Mathematics of the University of Leeds.}

\begin{document}
	
\begin{abstract}
    It is well known that ordered exponential fields with a compatible non-trivial valuation cannot be spherically complete, but there are some that are ``complete enough''. This paper gives analogues of Kaplansky's theorem on maximally valued fields that hold for a suitable class of elementary extensions of some ordered exponential fields with a compatible valuation. More precisely it does so for models of any theory $T_{\text{convex}}$ given by the expansion of a fixed complete o-minimal theory of ordered fields $T$ by a predicate $\mathcal{O}$ for a non-trivial $T$-convex valuation ring (so, in particular, for o-minimal ordered exponential fields with a compatible valuation when $T=\mathrm{Th}(\mathbb{R}, \exp)$).
    
    For $\lambda$ an uncountable cardinal, say that a unary type $p(x)$ over a model of $T_{\text{convex}}$ is \emph{$\lambda$-bounded weakly immediate} if its cut is defined by an empty intersection of fewer than $\lambda$ many nested valuation balls.
	Call an elementary extension \emph{$\lambda$-bounded wim-constructible} if it is obtained as a transfinite composition of extensions each generated by one element whose type is $\lambda$-bounded weakly immediate.
		
	I show that $\lambda$-bounded wim-constructible extensions do not extend the residue field sort and that any two wim-constructible extensions can be amalgamated in an extension which is again $\lambda$-bounded wim-constructible over both.
		
	A consequence of this is that given an uncountable cardinal $\lambda$, every model of $T_{\text{convex}}$ has a unique-up-to-isomorphism $\lambda$-spherically complete $\lambda$-bounded wim-constructible extension providing an analogue of Kaplansky's theorem. I call this extension the $T$-$\lambda$-spherical completion.
    Another consequence is that $T_\convex$ is \emph{definably spherically complete}.
    
	When $T$ is power-bounded wim-constructible extensions are just the immediate extensions. I discuss the example of power-bounded theories expanded by $\exp$ (\emph{simply exponential} theories).
\end{abstract}
	
\maketitle
\tableofcontents
	
\section{Introduction}
	
\subsection{Motivation}

Spherically complete valued fields play an important role in the theory of valued fields ever since Kaplansky's classic result that every valued field $(\bK,\cO)$ of $0$-equicharacteristic has a spherically complete immediate extension and that such extension is unique up to isomorphism \cite{kaplansky1942maximal}.

Analogues of this result have been given for various expansions of a valued field.
For example in \cite{aschenbrenner2018maximal} all valued differential fields are shown to have an immediate spherically complete \emph{strict} extension and all $H$-fields to have an immediate spherically complete $H$-field extension. It is also shown there that if a valued differential field has a discrete valuation, then the immediate spherically complete strict extension is unique up to isomorphism \cite[Prop.~7.1]{aschenbrenner2018maximal}.

A perfect analogue of Kaplansky's result is known to hold for models of $T_\convex$ --- the common theory of models of a complete o-minimal theory $T$ expanded by a non-trivial $T$-convex valuation ring --- provided that $T$ is  \emph{power-bounded}: in that case all models $(\bE, \cO)$ of $T_\convex$ have a unique up to isomorphism spherically complete immediate \emph{elementary} extension (cf \cite[Cor.~1.12]{kaplan2023t} where this comes as a consequence of the so-called residue valuation property established in \cite{dries2000field} and \cite{tyne2003t}). 


A result such as the one above is known to be false if restated verbatim when $T$ is an o-minimal theory defining exponentiation such as $T_{\exp}$ or $T_{an,\exp}$, the theories, respectively, of the field of reals expanded by the natural exponentiation $\bR_{\exp}$ and of its expansion $\bR_{an, \exp}$ by all restricted analytic functions.
The reason is twofold: on the one hand when $T$ defines exponentiation, immediate elementary extensions of $(\bE, \cO)$ must in fact be dense extensions (see Remark~\ref{rmk:Imm_in_exp_theories}), on the other hand as a consequence of \cite[Thm.~1]{kuhlmann1997exponentiation} a non-trivial $T$-convex valuation on an o-minimal field defining exponentiation cannot be spherically complete.

Nevertheless, literature strongly suggests fields of transseries which enjoy ``some degree of spherical completeness'' as useful structural analogues of Hahn fields in the presence of an exponential. Examples are the field of surreal numbers $\No$, which can be canonically expanded to an elementary extension of $\bR_{an, \exp}$ (cf \cite{dries2001fields}\nocite{dries2001erratum}) and some so called $\kappa$-bounded Hahn fields $\bR\Lhp \fM\Rhp_\kappa$, which also can be expanded to models of $T_{an, \exp}$ provided they define an $\omega$-map (cf \cite{kuhlmann2005kappa} and \cite{berarducci2023exponential}).

This calls for an analogue of Kaplansky's result capable of encompassing the cases of $T_{\exp}$ and $T_{an, \exp}$.

The main goal of this paper is to provide such a result in the greater generality of all o-minimal fields expanded by a $T$-convex valuation by showing that for all uncountable cardinals $\lambda$, the class of $\lambda$-spherically complete elementary extensions has a unique prime element (weakly initial object) up to isomorphism which I call its \emph{$T$-$\lambda$-spherical completion}, and that such extension does not enlarge the residue field sort.

This will be done by analyzing what I will call \emph{$\lambda$-bounded wim-constructible extensions}, which are some extensions that embed (over their base model) in every  $\lambda$-spherically complete elementary extension.
More precisely an extension $(\bE_*, \cO_*)$ of some base model $(\bE, \cO)$ is \emph{$\lambda$-bounded wim-constructible} if there is an ordinal indexed $\dcl_T$ basis $(x_i: i< \alpha)$ of $\bE_*$ over $\bE$ with the property that the type of $x_i$ over the elementary substructure generated over $\bE$ by $\{x_j: j<i\}$ is given by an intersection of fewer than $\lambda$ many valuation balls.
Therefore the study of wim-constructible extensions will start from the study of unary types in $T_\convex$.

When $T$ is power-bounded, the rv-property (cf \cite{dries2000field}, \cite{dries2002minimal}, \cite{tyne2003t}, \cite{nowak2007proof}), ensures that for large enough $\lambda$, the $\lambda$-bounded wim-constructible extensions are precisely the immediate ones, which entails their closure both under left and right factors and allows for more streamlined inductive amalgamation arguments. These closure properties do not seem to be available in general whence the amalgamation needs a new ad-hoc strategy.
	
\subsection{Setting} Recall that given a complete o-minimal theory of fields $T$ and $\bE\models T$, a $T$-convex subring $\cO \subseteq \bE$ is an order-convex subset of $\bE$ closed under total continuous definable functions without parameters. The theory $T_\convex$ is defined as the expansion of $T$ with a predicate $\cO$ axiomatized as being $T$-convex and non-trivial (i.e.\ $\exists x \notin \cO$). 
As proven by van den Dries and Lewenberg in \cite[(3.10)-(3.15)]{dries1995t} such a theory is complete and weakly o-minimal.
		
Moreover, unary types over a model $(\bE, \cO)$ of $T_\convex$ are almost always completely determined by their reducts to the underlying o-minimal structure $\bE$. More specifically, they are, unless they force the realization of the cut between the valuation ring $\cO$ and the elements greater than $\cO$ (from now on $\cO$-\emph{special cut}), in which case their reduct to $\bE$ has exactly two extensions to a type over a model of $T_\convex$ (\cite[(3.6)-(3.7)]{dries1995t}). 
		 
This makes it convenient to rephrase the study of types over $(\bE, \cO)$ as the study of the relationship between types over its o-minimal reduct $\bE$ (that is cuts over $\bE$) and $\cO$.
		
Notice that since $T$ has definable Skolem functions, all types \emph{are} types over a model of $T$: namely given a set of parameters $A\subseteq \bE \models T$, any type over $A$ extends uniquely to a type over the definable closure of $A$, which is a model of $T$. We will therefore notationally distinguish between types for $T$ and types for $T_\convex$, by writing $\tp(x/\bE)$ for the type of $x$ over $T\models \bE$ and $\tp(x/(\bE, \cO))$ for the type of $x$ over the expanded structure $(\bE,\cO)$ whenever $\cO$ is a $T$-convex subring of $\bE$.
		
\subsection{Main results} I define a type over $\bE$ to be \emph{weakly $\cO$-immediate} if it is determined by a finitely consistent conjunction of $\cO$-valuation balls (recall that valuation balls, being order-convex, are type-definable over the o-minimal reduct $\bE$ by an ordered-convex partial type) and then show that if some type $\tp(x/\bE)$ over the o-minimal reduct $\bE$ is weakly $\cO$-immediate, then $\bE \langle x \rangle:=\dcl_T(\bE, x)$ does not realize the $\cO$-special cut (Proposition~\ref{prop:wim_main}), hence $\tp(x/\bE)$ extends to a unique $\tp(x/(\bE, \cO))$. This requires most of the work toward Theorem~\ref{introthm:A} below and is achieved through a delicate argument (Lemma~\ref{lem:lem0}) showing that given a regular enough monotone parameter-definable bijection $f$, if the image under $f$ of a translate of an $\cO$-module (as is the set of realizations of an immediate type) is contained in $\cO$, then either its standard part has an endpoint or it is an open interval.
		
With similar techniques and building on previous work of van den Dries and Lewenberg in \cite{dries1995t} and \cite{dries1997t}, I show also that if $\tp(x/\bE)$ is weakly $\cO$-immediate, then $\bE \langle x \rangle$ with the unique extension of $\cO$ does not properly extend the residue field sort (Proposition~\ref{prop:wim_cof}).
		
The resulting picture (Theorem~\ref{thm:types_rest}) is that unary elementary extensions $(\bE \langle x \rangle, \cO')\succ (\bE, \cO)$ can be partitioned into 3 classes:
\begin{enumerate}
	\item the ones generated by a weakly $\cO$-immediate $y$ (\emph{wim-generated});
	\item the ones properly extending the residue field sort (\emph{residual});
	\item the ones such that $\bE \langle x \rangle \setminus \bE= \bE + M$ where $M \subseteq \bE \langle x\rangle$ is such that $\val_{\cO'}M \cap \val_{\cO}\bE = \emptyset$ (\emph{purely valuational}). 
\end{enumerate}
		
The technical result is more precisely stated as a partition of the unary extensions of the o-minimal reduct $\bE$, once $\cO$ is fixed, into 5 \emph{mutually exclusive} classes: 
				
\begin{thmA}[\ref{thm:types}]\label{introthm:A}
	If $\bU \succ \bE$ and $x \in \bU \setminus \bE$, and $\cO$ is a non-trivial $T$-convex subring of $\bE$, then one and only one of the following holds
	\begin{enumerate}
		\item $\bE\langle x \rangle =\bE\langle y \rangle$ for some weakly $\cO$-immediate $y$;
		\item $\bE \langle x \rangle = \bE \langle b \rangle$ for some $b$ such that $\cO < b < \bE^{>\cO}$;
		\item $\bE \langle x \rangle = \bE \langle z \rangle$ for some $z$ such that for some (equiv.\ any) $\bK \preceq \bE$ maximal in $\{\bK \preceq \bE: \bK \subseteq \cO\}$, $\bK$ is cofinal in $\bK \langle z \rangle$;
		\item $\bE\langle x \rangle=\bE \langle d \rangle$ for some $d>\bE$;
		\item none of the above.
	\end{enumerate}
\end{thmA}

\smallskip

Finally, I show that any two $\lambda$-bounded wim-constructible extensions $(\bE_0, \cO_0)$ and $(\bE_1, \cO_1)$ of any $(\bE, \cO)$ amalgamate over $(\bE, \cO)$ in an extension $(\bE_2, \cO_2)$ that is $\lambda$-bounded wim-constructible over both $(\bE_0, \cO_0)$ and $(\bE_1, \cO_1)$ (Lemma~\ref{lem:wim_amalgamate}).
This is non-trivial because wim-constructibility may not be stable under taking left-factors. The amalgamation is thus given, roughly speaking, by an ad-hoc transfinite insertion procedure to turn a wim-construction for a ``saturated enough'' wim-constructible extension $(\bE_2, \cO_2)$ of $(\bE_1, \cO_1)$ into one which has an initial segment isomorphic to a wim-construction for $(\bE_0, \cO_0)$. Lemma~\ref{lem:wim_amalgamate} entails the second main result.

\begin{thmA}[\ref{thm:B}]\label{introthm:B}
	Let $(\bE, \cO)\models T_\convex$, and $\lambda$ be an uncountable cardinal. There is a unique-up-to-non-unique-isomorphism $\lambda$-spherically complete $\lambda$-bounded wim-constructible extension and it elementarily embeds in every $\lambda$-spherically complete extension of $(\bE, \cO)$.
\end{thmA}

The amalgamation provided by Lemma~\ref{lem:wim_amalgamate} is rather weak, so it is natural to ask if this can be strengthened (at least for some specific theories $T$). I therefore pose the question whether wim-constructible extensions are closed under taking right (or left) factors (Question~\ref{quest:wimc-factors}~(1) and (2) respectively). Positive answers to these questions would allow for better understanding of $T$-$\lambda$-spherical completions (Remark~\ref{rmk:quest-relevance}). 
		
When $T$ is power-bounded, by the rv-property (see \cite[Sec.~9 and 10]{dries2000field} and \cite[Ch.~12 and 13]{tyne2003t}), wim-constructible extensions are precisely the immediate ones and thus Questions~\ref{quest:wimc-factors}~(1) and (2) have a positive answer.
		
When $T$ is the expansion of a polynomially bounded theory by $\exp$, the following analogue of the rv-property holds: if $(\bE, \cO) \models T_\convex$, $f: \bE \to \bE$ is definable in $\bE$, and $x\in \bU \setminus \bE$ is weakly $\cO$-immediate, then there is an $\bE$-definable composition of exponentials, translations and changes of signs $g$ such that $g^{-1} f(x)$ is weakly $\cO$-immediate (Proposition~\ref{prop:gne_ppty}). This entails a very partial result (Corollary~\ref{cor:wimc_1wim}) toward a positive answer to Question~\ref{quest:wimc-factors}~(1) for such $T$.

\smallskip

The paper also includes an alternative proof of the rv-property for power-bounded o-minimal theories (Section~\ref{sec:pbdd}) and a straightforward generalization of results in \cite{dries2002minimal} from the polynomially bounded case to the power-bounded case (Subsection~\ref{ssec:known_results}) which is needed to point out a property of complete o-minimal theories that are obtained as expansions of a power-bounded theory by an exponential (Corollary~\ref{cor:equivalent_simply_exp}). As the author was recently made aware of, a similar generalization of the same results in \cite{dries2002minimal} can be found in \cite{foster2010power}, however it is stated in a slightly more restrictive context than the one needed in this paper.

\smallskip

The following corollary of the analysis of unary types in $T_\convex$ provided by Theorem~\ref{introthm:A} is also worth mentioning in this introduction as it connects with a question in \cite{bradley-williams2023spherically}.

\begin{thmA}[\ref{thm:definable-spherical-completeness}]\label{introthm:C}
    For every complete o-minimal $T$, the theory $T_\convex$ is \emph{definably spherically complete} in the sense that every definable family of nested valuation balls has a non-empty intersection.
\end{thmA}

This was previously known only when $T$ is power-bounded and in this new generality it entails a negative answer to \cite[Question~1.1.]{bradley-williams2023spherically} of whether definably spherically complete expansions of the theory of valued fields always admit spherically complete models: in fact, when $T$ defines exponentiation, as already mentioned, it is known that there are no spherically complete models of $T_\convex$.
\subsection{Structure of the paper} Section~\ref{sec:set-up} is dedicated to the set-up: the first two subsections recall some standard order- and valuation-theoretic facts, such as the relation between nested intersections of valuation balls and p.c.-sequences (Subsection~\ref{ssec:breadth}) and give a brief review of $T$-convexity, with particular regard of the basic facts needed throughout the note (Subsection~\ref{ssec:Tconvex_review}).	 Subsection~\ref{ssec:TechLem} contains some technical Lemmas used throughout the paper.
		
\smallskip
Section~\ref{sec:mainres} contains the proofs of Theorems~\ref{introthm:A}, \ref{introthm:B}, and \ref{introthm:C} (Subsections~\ref{ssec:unary_types}, \ref{ssec:wimext}, and \ref{ssec:definable_spherical_completeness} respectively) along with some concluding remarks together with Question~\ref{quest:wimc-factors} (Subsection~\ref{ssec:quest}).
		
\smallskip
Section~\ref{sec:pbdd} briefly reviews the rv-property for power-bounded theories that was already established in \cite[Ch.~12 and 13]{tyne2003t}, giving a shorter and somewhat different proof.
		
\smallskip
Section~\ref{sec:exp} analyzes weakly immediate types in \emph{simply exponential theories}, that is, exponential theories $T$ in some language $L$ for which there is a language $L_0$, such that $T|_{L_0}$ is power-bounded and $T$ is an expansion by definitions of $T|_{L_0\cup \{\exp\}}$.
	
\subsection{Acknowledgments}
The author is very grateful to his doctoral supervisors Vincenzo Mantova for his support and constant feedback and Dugald Macpherson for his support and feedback. This paper is part of a PhD project at the University of Leeds and is supported by a Scholarship of the School of Mathematics. Its content, except for Subsection~\ref{ssec:definable_spherical_completeness}, is presented in Chapter~3 of the author's doctoral dissertation \cite{wreo36168}.

Thanks to Marcus Tressl for pointing out connections of the material in this paper with the previous literature (see Remark~\ref{rmk:previous_lit}).
Thanks to Lou van den Dries for several comments and suggestions made in examining my thesis and to Pantelis Eleftheriou for some comments concerning exposition in the thesis. Such comments and suggestions thereby also improved the content of this paper.
Thanks also to David Bradley-Williams for making me aware of the question about definable spherical completeness and to David Bradley-Williams and Immi Halupczok for some encouraging comments about the application of the main results to that problem.
Thanks to an anonymous referee for making me aware of \cite{foster2010power}, Chapter~6 of which has overlaps with Subsection~\ref{ssec:known_results} (as explained therein).
Thanks to another anonymous referee for many suggestions and remarks, some of which pointed out some minor errors allowing the author to fix them and some others led to improvements in the presentation.

\section{Set-up}\label{sec:set-up}
	
This section gives the set-up, recalling some facts relevant for our analysis of $T_\convex$ and some technical lemmas.
	
\subsection{Breadth, p.c.-sequences and weakly immediate types}\label{ssec:breadth}

In the following a \emph{partial type} $p(\var{x})$ is just a set of formulae closed under consequences in the tuple of free variables $\var{x}$. The set of realizations of $p$ within a structure $\bU$ will be denoted by $p(\bU)$.
		
If $\bE$ is a field and $\cO\subseteq \bE$ is a valuation ring of $\bE$, we will denote by $\co$ its maximal ideal. If $\cO$ has a subscript or superscript, we will denote the corresponding maximal ideal as $\co$ with the same subscript/superscript. We also write the associated \emph{dominance relation} on $\bE$ as $x \preceq_\cO y \Leftrightarrow \big(x/y\in \cO \text{\;or\;} x=y=0\big)$ omitting the subscript $\cO$ if there is no ambiguity. Similarly for $x \asymp y \Leftrightarrow (x \preceq y\; \& \; y\preceq x)$, $x \prec y \Leftrightarrow (x\preceq y\; \& \; y \not \preceq x)$, and $x \sim y \Leftrightarrow x-y \prec x$.
		
We will denote the value-group of $(\bE,\cO)$ by $(\val(\bE, \cO), +):=(\bE^{\neq 0}/(\cO \setminus \co), \cdot)$ or $\val_\cO(\bE)$ and the valuation by $\val_{\cO} : \bE \to \val(\bE, \cO) \cup \{\infty\}$. We use the classical convention that the value-group is ordered by stipulating that $\val(\bE, \cO)^{\ge0}=\val_\cO(\cO^{\neq 0})$ and $\val_\cO(0)=\infty$. In particular $\val_{\cO}(x)>\val_{\cO}(y)$ if and only if $x \prec_{\cO} y$.

The residue field $\cO/\co$ will be denoted by $\res(\bE, \cO)$ or $\res_{\cO}(\bE)$ and the quotient map by $\res_{\cO}: \cO \to \cO/\co$. Finally the residue-value sort $(\bE^{\neq0}/(1+\co), \cdot)$ will be denoted by $\rv_\cO(\bE)$ and the quotient map by $\rv_{\cO} : \bE^{\neq 0} \to \rv_\cO(\bE)$. Again the subscript $\cO$ will be omitted if there is no ambiguity. A \emph{valuation ball} in $(\bE, \cO)$ will be a subset of the form $a+\cO b$ or $a+\co b$ with $a,b \in \bE$ and $b \neq 0$.
		
If $\bU$ is a first-order structure and $\bE \subseteq \bU$ is an elementary substructure we will write $\bE \preceq \bU$, or $\bE \prec \bU$ to signify that moreover $\bE \neq \bU$.
		
If $T$ is an o-minimal theory expanding the theory of densely ordered Abelian groups, $\bU \succ \bE \models T$, and $x \in \bU\setminus \bE$, then $\bE \langle x \rangle_T$ will denote the definable closure of the set $\bE\cup \{x\}$ (which is the minimum elementary substructure of $\bU$ containing $\bE$ and $x$, because $T$ has definable Skolem functions when expanded by a non-zero constant).

Throughout this paper, if $\bullet$ is a binary relation symbol and $S,C$ are sets, then $S^{\bullet C}$ will be an abbreviation for $\{y\in S: \forall c \in C,\; y\bullet c\}$ and if $C=\{c\}$ is a singleton we will write $S^{\bullet c}$ for $S^{\bullet C}$.

\begin{definition}
	Let $\bE$ be an expansion of a dense linearly ordered group. We say that a partial $1$-ary type $p(\var{x})$ over $\bE$ is \emph{convex} if 
	\[p(\var{x}_0), p(\var{x}_2), \var{x}_0< \var{x}_1 < \var{x}_2 \vdash p(\var{x_1}),\]
	or equivalently if it defines a convex subset in every elementary extension of $\bE$.
	If $p(\var{x})$ is a convex partial type over $\bE$, we will write $\cN(p)$ for the set of formulae
	\[\cN(p)(\var{x})=\{a<\var{x}<b: a,b \in \bE \cup \{\pm \infty\},\; p(\var{x}) \vdash a<\var{x}< b\}.\]
	In analogy with Kaplansky's terminology in \cite{kaplansky1942maximal}, we will call the \emph{breadth of $p$} the partial type $\Br(p)(\var{y})$ implied by the formulae
	\[\{|\var{y}|<b-a: a,b \in \bE \cup \{\pm \infty\},\; p(\var{x}) \vdash a<\var{x}< b\}.\]
	If $p=\tp(x/\bE)$ for $x$ in some elementary extension of $\bE$ and $p$ is convex, we write $\Br(x/\bE)$ for $\Br(p)$.
\end{definition}
	
\begin{remark}
	If $\bE$ is o-minimal, and $p$ is any unary, complete, and non-isolated type, then $p(\var{x})$ is convex and $\cN(p)(\var{x})\vdash p(\var{x})$.
\end{remark}
	
\begin{remark}
	For every non-realized convex unary partial type $p$ over $\bE$, $q:=\Br(p)$ is a convex subgroup, meaning that $q$ is convex and
	\[q(\var{x}_0), q(\var{x}_1) \vdash q(\var{x}_1-\var{x}_0).\]
	Moreover, for every elementary extension $\bU$ of $\bE$ realizing $p$,
	\[q(\bU)=\{y \in \bU: y+p(\bU)\subseteq p(\bU) \text{ or } -y+p(\bU)\subseteq p(\bU) \}.\]
    
    In fact, since by construction $q$ is convex and $q(\var{x}_0) \vdash q(-\var{x}_0)$, to see that $q(\var{x}_0), q(\var{x}_1) \vdash q(\var{x_0}+\var{x_1})$, it is enough to prove that $q(\var{x}_0)$ implies $q(2\var{x}_0)$. For this, in turn, it is enough to show that for all $b$ and $a$ in $\bE$ such that $p(\var{x})\vdash a<\var{x}<b$, there are $a', b' \in \bE$, such that $p(\var{x})\vdash a'<\var{x}<b'$ and $2(b'-a') \le b-a$. When $p(\var{x})$ is convex and non-realized, it must imply one of the two formulas $\var{x}>(a+b)/2$ and $\var{x}<(a+b)/2$, therefore we can choose such pair $(a', b')$ to be $(a, (a+b)/2)$ or $((a+b)/2, b)$ according to whether $p(\var{x})$ implies $\var{x}<(a+b)/2$ or $\var{x}>(a+b)/2$.
\end{remark}

\begin{remark}\label{rmk:basic_vt}
    Let $\bE \subseteq \bU$ be a field extension and $\cO\subseteq \bU$ be a valuation ring. Suppose $c \in \bE$ and $x \in\bU$.
    Then $\val_\cO(x-c)$ is the maximum of $\val_{\cO}(x-\bE)$ if and only if $\rv_\cO(x-c) \notin \rv_\cO(\bE)$. In fact for every $d \in \bU$, $\val_{\cO}(x-c)<\val_{\cO}(x-d)$ if and only if $\rv_{\cO}(d-c)=\rv_{\cO}(x-c)$.
\end{remark}
    
Recall the following two standard valuation-theoretic facts. 
	
\begin{lemma}\label{lem:wim_basic}
	If $\bU\supseteq \bE$ is a field extension, $\cO\subseteq \bU$ is a valuation subring, and $x \in\bU \setminus \bE$, then the following are equivalent:
	\begin{enumerate}
		\item $\val_{\cO}(x-\bE)$ has no maximum;
		\item $\val_\cO (x- \bE)\subseteq \val_\cO \bE$ and $\res_\cO(\cO\cap \bE)=\res_\cO(\cO \cap(x \bE^{\neq 0}+\bE))$.
	\end{enumerate}
	\begin{proof}
		Let $\val:=\val_\cO$, $\res:=\res_{\cO}$.
		$\lnot (2) \Rightarrow \lnot (1)$; if for some $c \in \bE$, $\val (x-c)\notin \val \bE$, then $\val(x-c)$ is maximum in $\val(x-\bE)$ by Remark~\ref{rmk:basic_vt}. Similarly, if $\res(\cO\cap \bE)\neq\res(\cO \cap(x \bE^{\neq 0}+\bE))$, then there are $c, d \in \bE$ such that $xc+d \in \cO$ and $(xc + d +\co) \cap \bE = \emptyset$. This, again by Remark~\ref{rmk:basic_vt} implies that $\val(x-d/c)$ is maximum in $\val(x-\bE)$.
		$(2) \Rightarrow (1)$; given $c \in \bE$, since $\val(x-c) \in \val \bE$ we can find $d\in \bE$ such that $(x-c)/d=x/d-c/d \in \cO$, but then since $\res(x/d-c/d) \in \res(\bE)$, $(x/d-c/d + \co) \cap \bE \neq \emptyset$ and we can find $b \in\bE$ such that $x/d-c/d-b \in \co$ whence $x-c-bd \in \co \cdot d$ and $\val(x-c-bd)> \val(x-c)$. 
	\end{proof}
\end{lemma}

\begin{lemma}\label{lem:wim_equiv}
	Let $\bU\supseteq \bE$ be fields and $\cO \subseteq \cO'$ be valuation subrings of $\bU$ such that $\bE \cap \cO=\bE\cap \cO'$. For $x \in \bU \setminus \bE$, the following are equivalent:
	\begin{enumerate}
		\item $\val_{\cO}(x-\bE)$ has no maximum;
		\item $\val_{\cO'}(x-\bE)$ has no maximum.
	\end{enumerate}
	\begin{proof}
		Let $\alpha: \val(\bU, \cO) \to \val(\bU, \cO')$ denote the ordered group homomorphism given by $\alpha (\val_{\cO}(y))=\val_{\cO'}(y)$.
			
		$\lnot(1) \Rightarrow \lnot(2)$ If $\val_{\cO}(x-\bE)$ has a maximum $\gamma=\val_\cO(x-c)$, then $\val_{\cO'}(x-c)=\alpha(\gamma)$ is maximum in $\val_{\cO'}(x-\bE)$ because $\alpha$ is an ordered group homomorphism.
			
		$(1) \Rightarrow (2)$ By Lemma~\ref{lem:wim_basic} it suffices to show $\val_{\cO'}(x-\bE)\subseteq \val_{\cO'}\bE$ and $\res_{\cO'}(\cO' \cap (\bE x + \bE))\subseteq \res_{\cO'}(\bE)$.
			
		To see the first assertion, let $c \in \bE$, so by Lemma~\ref{lem:wim_basic} $\val_{\cO}(x-c) \in \val_{\cO}(\bE)$; since $\alpha$ is a group homomorphism, this implies $\val_{\cO'}(x-c) \in \val_{\cO'}(\bE)$.
			
		For the second notice that if $c \in \bE$, $d \in \bE^{\neq 0}$, and $d(x-c) \in \cO'$, then $d(x-c) \in \cO$. In fact, $\val_{\cO}(d(x-c)) \in \val_{\cO}(\bE)$ by (1) and Lemma~\ref{lem:wim_basic} and since the kernel of $\alpha|_{\val_{\cO}(\bE)}$ is trivial by the hypothesis $\cO\cap \bE=\cO' \cap \bE$, we must have $\val_{\cO}(d(x-c))\ge 0$.
			
		Now, since by (1) $\val_{\cO}(d(x-\bE))$ has no maximum, this implies that there is $b \in \bE$ such that $d(x-c)-b \in \co$. We must have $d(x-c)-b\in \co'$, because otherwise $\val_{\cO}(d(x-c)-b) \in \val_{\cO}(\bE) \cap \ker(\alpha)$, which would contradict the injectivity of $\alpha|_{\val_{\cO}(\bE)}$, since by construction $\val_{\cO}(d(x-c)-b)>0$.

        Thus we showed that $\res_{\cO'}(\cO' \cap (\bE x + \bE))\subseteq \res_{\cO'}(\bE)$, which concludes the proof.
	\end{proof}
\end{lemma}
	
\begin{lemma}\label{lem:wim_basic_o}
	Let $\bU\supseteq \bE$ be ordered fields, $\cO\subseteq \bU$ be a convex valuation subring, and $x \in\bU \setminus \bE$. If $\val_{\cO}(x-\bE)$ has no maximum, then $\val_{\cO}((x-\bE)^{>0})$ and $\val_{\cO}((x-\bE)^{<0})$ are cofinal in $\val_{\cO}(x-\bE)$.
	\begin{proof}
		Let $\val:=\val_\cO$. Notice that if $\val(x-\bE)$ has no maximum, then by Lemma~\ref{lem:wim_basic}, $\val(x-\bE) \subseteq \val \bE$. Moreover $\val(x-c_1) \ge \val(x-c_0)$ if and only if $c_1 \in c_0+\cO (x-c_0)$. Thus it suffices to show that for every valuation ball $B:=c+\cO d$ with $c,d \in \bE$ and $x \in c+\cO d$, $B\cap \bE$ contains both elements above and below $x$.  
        Suppose not: then $\bE^{>B}=\bE^{>x}$ or $\bE^{<B}= \bE^{<x}$. In both cases we would get the contradiction that $\val(\bE-x)$ has a maximum. In fact, if $B':=c'+\cO d'$ with $c',d' \in \bE$ is such that $B' \subseteq B$ and $\bE^{>B}=\bE^{>B'}$, then $B'=B$. Thus any valuation ball $B'$ with center in $\bE$ containing $x$ would have radius $\val(d')\le \val(d)$.
	\end{proof}
\end{lemma}
	
From now on $\bE$ will be an o-minimal expansion of a real closed field.
	
\begin{definition}\label{defn:wim&pc}
	Let $\cO$ be a convex valuation ring on $\bE$ and $p(\var{x})$ be a complete type over $\bE$. We call $p(\var{x})$ \emph{weakly $\cO$-immediate} if for some (equiv.\ any) realization $x$ of $p$, $\val_{\cO'}(x-\bE)$ has no maximum for some (equiv.\ any, by Lemma~\ref{lem:wim_equiv}) extension $\cO'$ of $\cO$ to $\bE \langle x \rangle$.
		
	An \emph{$\cO$-pseudo-Cauchy sequence} in $\bE$ (\emph{$\cO$-p.c.\ sequence} for short) is an ordinal-indexed sequence $(x_i)_{i \in \lambda}\in \bE^\lambda$ such that for every $i<j<k<\lambda$, $x_j-x_k\prec x_i-x_j$. A \emph{pseudolimit} for such a $(x_i)_{i< \lambda}$ is any element $x$ such that $x-x_j\prec x-x_i$ for all $i<j$ or equivalently such that $x \in \bigcap_{i< \lambda}\big( x_i + \cO (x_i-x_{i+1})\big)$.
\end{definition}

\begin{remark}\label{rmk:previous_lit}
    An analysis of cuts of o-minimal structures with respect to convex valuations had been conducted in \cite{tressl2006pseudo}. In the terminology and notation of \cite[Def.~2.2 and 2.3]{tressl2006pseudo}, a type $p$ is weakly $\cO$-immediate if and only if $\sign(p)=0$, $V(p)\supseteq \cO$ and $p$ is not an $\cO$-limit.
\end{remark}

\begin{remark}
	For every convex valuation ring $\cO$ in $\bE$, valuation balls are type-definable in the o-minimal reduct $\bE$ \emph{by a convex partial type}. Indeed for $a, b \in \bE$ we have
	\[\var{x}a-b \in \co \leftrightarrow |\var{x}a-b|<\bE^{>\co}	\]
	\[\var{x}a-b \in \cO \leftrightarrow |\var{x}a-b|<\bE^{>\cO}	\]
	In particular, any intersection of valuation balls is still type-definable by a convex partial type in $\bE$.
\end{remark}

\begin{lemma}\label{lem:wim_pc}
	The following are equivalent for a complete unrealized unary type $p$ over $\bE$ and a convex valuation ring $\cO\subseteq \bE$:
	\begin{enumerate}
		\item $p$ is weakly $\cO$-immediate over $\bE$;
		\item there are families $(x_i)_{i \in I}\in \bE^I$, $(y_i)_{i \in I}\in \bE^{I}$ such that for every $\bU \succ \bE$ and every convex subring $\cO'\subseteq \bU$ with $\cO' \cap \bE=\cO$,
		\[p(\bU)=\bigcap \{ x_i+\cO' y_i : i \in I\}\]
		\item there is an \emph{increasing} $\cO$-p.c.-sequence $(x_i)_{i<\lambda}$ in $\bE^{<x}$ with $x\models p$, such that for every $\bU \succ \bE$ and every convex subring $\cO'\subseteq \bU$ with $\cO' \cap \bE=\cO$
		\[p(\bU)=\bigcap \{ x_i+\cO' (x_{i+1}-x_i) : i \in \lambda\}.\]
	\end{enumerate}
    \begin{proof}
		$(1)\Rightarrow (3)$ if $p$ is weakly $\cO$-immediate over $\bE$ and $\cO'$ is the maximal extension of $\cO$ to $\bU$, that is, the largest convex subring of $\bU$ such that $\cO=\cO'\cap \bE$, then $\val_{\cO'}(x-\bE)$ has no maximum by Lemma~\ref{lem:wim_equiv} and $\val_{\cO'}((x-\bE)^{>0})$ is cofinal in $\val_{\cO'}(x-\bE)$ by Lemma~\ref{lem:wim_basic_o}.
						
		So we can pick $(x_i)_{i< \lambda} \in (\bE^{<x})^\lambda$ such that $\val_{\cO'}(x-x_i)$ is strictly increasing and cofinal in $\val_{\cO'}(x-\bE)$. This will be an $\cO$-p.c.\ sequence in $\bE$ because for $i<j<k<\lambda$, 
		\[\val(x_i-x_j)=\val(x_i-x)<\val(x_j-x)=\val(x_j-x_k).\]
		Notice that for each $i< \lambda$, $x \in x_i + \cO' \cdot (x_i-x) = x_i + \cO' \cdot (x_i-x_{i+1})$. So $p(\bU) \subseteq \bigcap_{i<\lambda} x_i + \cO' \cdot (x_i-x_{i+1})$.
		Notice also that $\bigcap_{i<\lambda} x_i + \cO \cdot (x_i-x_{i+1}) = \emptyset$ because if $y \in \bigcap_{i<\lambda} x_i + \cO \cdot (x_i-x_{i+1})$, then $\val(y-x)>\val(y-x_i)$ for all $i<\lambda$ and $\val(x-x_i)$ would not be cofinal in $\val(x-\bE)$.
		Since $\cO'=\{t\in \bU: |t|<\bE^{>\cO}\}$ is defined by a convex partial type over $\bE$, $\bigcap_{i<\lambda} x_i + \cO' \cdot (x_i-x_{i+1})$ is also defined by a convex partial type $q$ over $\bE$, but since $q(\bE)=\bigcap_{i<\lambda} x_i + \cO \cdot (x_i-x_{i+1}) = \emptyset$, $q$ is complete and $p=q$.
			
		On the other hand if $\cO''$ is another convex extension of $\cO$ to $\bU$, we have $\cO'' \subseteq \cO'$ and thus
		\[ x_{i+1} + \cO'' ( x_{i+2} - x_{i+1} ) \subseteq x_{i+1} + \cO' ( x_{i+2} - x_{i+1} ) \subseteq x_{i} + \cO'' ( x_{i+1} - x_{i} ).\]
		where the second inclusion is because $\cO'(x_{i+2}-x_{i+1}) \subset \cO'' (x_{i+1}-x_i)$ and clearly $x_{i+1} \in x_i+\cO''(x_{i+1}-x_{i})$. 
			
		$(3) \Rightarrow (2)$ is obvious.
			
		$(2) \Rightarrow (1)$ for some $\bU$, $p(\bU)$ is non-empty, so the family of valuation balls $\{x_i+\cO y_i: i \in I\}$ must be finitely consistent and hence, by the ultrametric inequality, totally ordered by inclusion. Moreover since $p$ is not realized in $\bE$, we have that $\bigcap\{x_i+\cO y_i: i \in I\}=\emptyset$ so $\{\val_\cO (y_i): i \in I\}$ has no maximum because if $\val_\cO (y_j)$ was such a maximum, it would follow that $x_j+\cO y_j \subseteq x_i+ \cO y_i$ for all $i \in I$ by the finitary consistency of the family and the ultrametric inequality.
			
		Let $x \in p(\bU)$, and suppose toward contradiction that for some $c \in \bE$, $\val_{\cO'}(x-c)$ is maximum in $\val_{\cO'}(x-\bE)$. This would imply that $\val_{\cO'}(x-c) > \val_{\cO'} (y_i)$ for every $i \in I$, because $\val_{\cO'}(x-x_i) \ge \val_{\cO'} (y_i)$ and $\{\val_\cO (y_i): i \in I\}$ has no maximum.
		Thus, $c \in \bE \cap \bigcap\{ x_i+ \cO' y_i: i \in I\}$, contradicting $\bigcap\{x_i+\cO y_i: i \in I\}=\emptyset$.
	\end{proof}
\end{lemma}
	
\begin{lemma}\label{lem:br_Omod}
	Let $\bE\prec \bU\models T$, $\cO$ a convex subring of $\bU$, and $p$ a weakly $(\cO \cap \bE)$-immediate type over $\bE$. Then
	\begin{enumerate}
		\item $\Br(p)(\bU)$ is an $\cO$-submodule of $\bU$;
       \item for any realization $x\in p(\bU)$, $\Br(p)(\bU)=p(\bU) - x$.
	\end{enumerate}
	\begin{proof}
        Let $x$ be a realization of $p$.
        
		(1) Since $p$ is weakly $(\cO\cap \bE)$-immediate, $\val(x-\bE)$ has no maximum. Moreover by Lemma~\ref{lem:wim_basic_o}, $\val_{\cO}((x-\bE)^{>0})$ and $\val_{\cO}((x-\bE)^{<0})$ are cofinal in $\val(x-\bE)$. It follows that $\{\val(b-a): a,b \in \bE,\, a<x<b\}$ is also cofinal in $\val(x-\bE)$ and has no maximum. Thus
        \[\Br(p)(\bU)=\{y \in \bU: \val(y)> \val (x-\bE)\}\]
        is an $\cO$-submodule of $\bU$.
            
        (2) Notice that in particular if $z \in \Br(p)(\bU)$, then $\val(z)> \val(x-\bE)$. Therefore for each $c \in \bE$, we have
        \[ c > x \Longleftrightarrow c - x > 0 \Longleftrightarrow c-x > z\cO \Longleftrightarrow c > x + z\cO.\]
        Similarly $c<x$ if and only if $c<x+z\cO$. It follows that $x+z\cO \subseteq p(\bU)$. Thus we showed $x+\Br(p)(\bU) \subseteq p(\bU)$.

        To see the other inclusion let $y\in p(\bU)$. Then clearly $|y-x|< b-a$ for every $a,b \in \bE$, with $a<p(\bU)< b$, hence $y-x \in \Br(p)(\bU)$.
	\end{proof}
\end{lemma}
	
\subsection{Review of \texorpdfstring{$T$}{T}-convexity and tame extensions}\label{ssec:Tconvex_review}
 
Let $T$ be a complete o-minimal theory expanding RCF in a language $L$. If $\bE \models T$, recall that a \emph{$T$-convex} subring of $\bE$, as defined in \cite{dries1995t}, is a convex subring of $\bE$ which is closed under continuous $T$-definable functions $f: \bE \to \bE$ (by $T$-definable we mean $\emptyset$-definable in $T$). It is said to be \emph{non-trivial} if $\bE \neq \cO$. It is not hard to see that if $\bK \preceq \bE$, then the \emph{convex hull} $\cO := \CH_\bE(\bK)$ of $\bK$ in $\bE$ is a $T$-convex subring.
	
Following van den Dries and Lewenberg in \cite{dries1995t}, denote by $L_\convex$ the language obtained by extending $L$ by a unary predicate $\cO$ and by $T_\convex$ the theory given by $T$ together with an axiom scheme stating that $\cO$ is a non-trivial $T$-convex valuation ring.

Recall that an ordered structure $M$ is \emph{weakly o-minimal} if and only if every definable subset of $M$ is a finite union of order-convex subsets. A theory $T$ expanding the theory of linear orders is said to be \emph{weakly o-minimal} if all of its models are weakly o-minimal (cf \cite{macpherson2000weakly}, \cite{Dickmann1987Elimination}).
	
\begin{theorem}[van den Dries and Lewenberg, (3.10) and (3.15) in \cite{dries1995t}]\label{thm:DriesLewMain}
	The theory $T_\convex$ is complete and weakly o-minimal. Moreover if $T$ eliminates quantifiers and has a universal axiomatization (resp. is model complete), then $T_\convex$ eliminates quantifiers (resp. is model complete).
\end{theorem}
	
A key ingredient is the fact that if $p(\var{x})$ is a type over an o-minimal structure $\bE$, then it can extend in at most two ways to a type over the expanded structure $(\bE, \cO)$. This will be extremely important throughout this chapter.
	
\begin{lemma}[van den Dries and Lewenberg, (3.6) and (3.7) in \cite{dries1995t}]\label{lem:O_extensions}
	If $p(\var{x})$ is a unary type over $\bE\models T$, $(\bE, \cO)\models T_\convex$, $x \models p$ and $\cO'$ is a $T$-convex valuation subring of $\bE \langle x \rangle$ with $\cO' \cap \bE = \cO$, then $\cO' \in \{\cO_x^-, \cO_x\}$ where
	\[\cO_x^-:=\CH_{\bE\langle x \rangle}(\cO), \qquad \cO_x:=\{y \in \bE\langle x\rangle:  |y|< \bE^{>\cO}\}.\]
\end{lemma}
	
\begin{remark}
	Of course if $\bE\langle x\rangle$ does not contain any $b$ with $\cO<b<\bE^{>\cO}$, then $\cO_x^-=\cO_x$. If instead $\bE\langle x\rangle$ contains such a $b$, then, by the exchange property, $\bE \langle x \rangle = \bE \langle b \rangle$ and $\cO_x=\cO_b$. For simplicity, in the following we will denote $\cO_x^-=\CH_{\bE\langle x \rangle}(\cO)$ by $\cO$ if there is no ambiguity.
\end{remark}

\begin{remark}\label{rmk:tame&def}
    Recall from \cite[(1.12)]{dries1995t} that an elementary extension $\bK \preceq \bE$ of models of $T$ is said to be \emph{tame} (notation $\bK \preceq_\tame \bE$) if $\bK$ is definably Dedekind-complete in $\bE$ in the sense that for every $\bE$-definable subset $X$ of $\bE$, if $X \cap \bK$ is bounded, then it has a supremum in $\bK$.
    Also recall that a type $p(\var{x})$ over a model $M\models T$ is said to be \emph{definable} if for every formula $\varphi(\var{x}, \var{y})$, the set $\{m \in M^{\var{y}}: p(\var{x}) \vdash \varphi(\var{x}, m)\}$ is definable with parameters in $M$.
    It was proven by Marker and Steinhorn in \cite{marker1994definable} that $\bK \preceq_\tame \bE$ if and only if for every tuple $\overline{c}$ of elements in $\bE$, $\tp(\overline{c}/\bK)$ is a definable type (see also \cite{pillay1994definability}).
\end{remark}

\begin{remark}\label{rmk:unarydef}
    A unary type $p$ over $\bK\models T$ is definable if and only if, given a realization $x$ of $p$, $\bK^{<x}$ is an interval or is empty. Hence all tame extensions of $\bK$ with $\dcl$-dimension one are isomorphic to $\bK \langle d \rangle$ for some $d>\bK$.
\end{remark}
	
Tame extensions of o-minimal structures are closely related to $T$-convex valuation rings: 
	
\begin{fact}[van den Dries and Lewenberg, (2.12) in \cite{dries1995t}]	
    If $\bE \models T$ and $\cO$ is a $T$-convex valuation ring, then $\bK\preceq\bE$ is maximal in $\{\bK' \prec \bE: \bK'\subseteq \cO\}$ if and only if $\bK \preceq_\tame \bE$ and $\bK+\co= \cO$.
\end{fact}
	
For every tame extension $\bK \preceq_\tame \bE$ there is an associated standard part map $\std_\bK: \CH_\bE(\bK) \to \bK$ uniquely defined by the property that for every $x \in \CH_\bE(\bK)$, $|\std_\bK(x)-x|<\bK^{>0}$.
	
\begin{theorem}[van den Dries, Sec.~1 in \cite{dries1997t}]\label{thm:std_def}\nocite{dries1998correction}
	If $(\bE, \cO) \models T_\convex$, and $\bK \preceq_\tame \bE$ with $\CH(\bK)=\cO$, then $\std_\bK: \cO \to \bK$ induces an isomorphism between $\bK$ and the induced structure on the imaginary sort $\cO/\co$.
\end{theorem}
	
In what follows we will also need this other result of van den Dries and Lewenberg.
	
\begin{lemma}[van den Dries and Lewenberg]\label{lem:resrank}
	Let $\cO$ be a $T$-convex subring of $\bU$, $x \in \bU \setminus \bE$, $\bK$ maximal in $\{\bK'\preceq \bE: \bK\subseteq \cO\}$, and $\bK_x$ maximal in $\{\bK'\preceq \bE \langle x \rangle: \bK \preceq \bK'\subseteq \cO\}$. Then $\bK=\bK_x$ or $\bK_x$ is a 1-dimensional extension of $\bK$, that is $\bK_x=\bK \langle y\rangle$ for some $y \in \bK_x \setminus \bK$.
	\begin{proof}
		This is an immediate Corollary of \cite[Lemma~5.3]{dries1995t}.
	\end{proof}
\end{lemma}
	
\subsection{Some technical Lemmas}\label{ssec:TechLem}
	Throughout this subsection $T$ will be any complete o-minimal theory expanding $\RCF$, $(\bE, \cO)\models T_\convex$, $\bU\succ \bE$ an extension of $\bE$, and $b \in \bU \setminus \bE$ an element realizing the special cut $\cO < b < \bE^{>\cO}$. Also, $\bK\preceq \bE$ will be maximal such that $\bK \subseteq \cO$ and for every unary function $f$ we will write $f'$ or $\partial f$ for its derivative, and $f^\dagger:=f'/f$ for its logarithmic derivative (where they are defined).
	Recall the following folklore fact about o-minimal structures.
	
\begin{lemma}\label{lem:Tdefarecof}
	Let $x \in \bU$ be such that $x > \bK$ and $f: \bU \to \bU$ be a $\bK$-definable function. Then there is a $T$-definable function $g$ such that $g(x) > f(x)$.
	\begin{proof}
		Without loss of generality we can assume that $f$ is continuous and increasing. If $f$ is itself $T$-definable, there is nothing to prove. We show by induction on $n>0$ that if we can write $f(-)=F(\overline{c},-)$, with $\overline{c}:=(c_0, \ldots, c_{n-1})$ some $n$-tuple of parameters from $\bK$ and $F$ a $T$-definable function, then there is an $(n-1)$-tuple $\overline{d}$ in $\bK$ and a $\overline{d}$-definable function $g$ with $f(x)< g(x)$.
		Without loss of generality we can assume that $n$ is minimal for $f$, so $\overline{c}$ has dimension $n$, and whenever a $T$-definable cell of $\bK$ contains $\overline{c}$, it is open.
		Now consider a $T$-definable cell decomposition of $\bK^{n+1}$ such that on every cell $F$ is continuous and strictly monotone or constant in every variable. 
		By the minimal choice of $n$, there will be a $T$-definable open cell $A$ in $\bK^n$ containing $\overline{c}$ and a continuous $T$-definable function $\gamma: A \to \bK$ such that on the open $(n+1)$-cell $(\gamma, \infty)\coloneqq\{(\overline{y},z) \in A \times \bK: \gamma(\overline{y})<z\}$, $F$ is continuous, strictly increasing in the last variable and strictly monotone in the first $n$ variables. If the type of $c_{n-1}$ over $\dcl_T(c_0, \ldots, c_{n-2})$ is not definable (see Remark~\ref{rmk:unarydef}), then there are $a,b \in \dcl_T(c_0, \ldots, c_{n-2})$ such that $c_{n-1}\in [a,b] \subseteq A_{c_0, \ldots, c_{n-2}}:=\{s: (c_0, \ldots, c_{n-2}, s) \in A\}$. So if we set
		\[g(t):= 1+ \max\{ F(c_0, \ldots, c_{n-2}, s, t): s \in [a,b]\},\]
		then for $t>d := \max\{ \gamma(c_0, \ldots, c_{n-2}, s): s \in [a,b]\} \in \dcl_T(c_0, \ldots, c_{n-2})$ we have $g(t)> f(t)$, hence $g(x)> f(x)$.
		If instead $c_{n-1}$ has a definable type over $\dcl_T(c_0, \ldots, c_{n-2})$, then we can without loss of generality assume that $c_{n-1}> \dcl_T(c_0, \ldots, c_{n-2})$ and in that case it suffices to choose
		$g(t):=1+\max\{ F(c_0, \ldots, c_{n-2}, s, t): 0 \le s \le t\}$.
	\end{proof}
\end{lemma}

Part of the following Proposition is contained in \cite[Prop.~2.10]{dries1997t}, and the idea of the proof is similar.
	
\begin{proposition}\label{prop:coffact}
	Let $(\bE, \cO) \models T_\convex$. There are no $\bE$-definable maps restricting to a continuous increasing bijection between some final segment of $\cO$ and some final segment of $\bE$. Similarly for $\cO$ and $\co$ and for $\co$ and $\bE$ in place of $\cO$ and $\bE$.
	\begin{proof}
		Let $f$ be an $\bE$-definable map, so $f(t)=g(t, \overline{c})$ for some $T$-definable $g$ and some $n$-tuple $\overline{c}$ from $\bE$. The statement ``for every $\overline{c} \in \bE^n$, $g(-,\overline{c})$ does not restrict to a definable continuous increasing bijection between a final segment of $\cO$ and a final segment of $\bE$'' is equivalent to a first order sentence in the language of $T_\convex$. Since the theory of $(\bE,\cO)$ is $T_\convex$, to prove the thesis it suffices to show the statement holds in some model of $T_\convex$.
			
		Let $\bK\models T$ have cofinality at least $|T|^+$ and consider $\bE=\bK \langle d \rangle$ with $d> \bK$ and $\cO=\CH_\bE(\bK)$. 
			
		By Lemma~\ref{lem:Tdefarecof}, $\bE$ has then cofinality at most $|T|$ whereas $\cO$ has cofinality at least $|T|^+$ hence there can be no continuous increasing bijection between a final segment of $\bE$ and a final segment of $\cO$.
		
		Also notice that in the same model $\co$ has the same cofinality as $\bE$, because given increasing bijections $f,g: \bK \to \bK$, if $f(d)<g(d)$ then $1/f^{-1}(d)<1/g^{-1}(d)$; this proves the second assertion.
			
		To show the last assertion, it suffices to exhibit a model in which $\bE$ and $\co$ don't have the same cofinality. This can be constructed for example as $\bE:=\bK \langle \{d_i: i<|T|^+\}\rangle$ where $\bK\langle d_i : i > j \rangle < d_j$ and $\cO=\CH_\bE(\bK)$: in such a model the cofinality of $\co$ is $|T|^+$, whereas the cofinality of $\bE$ is the cofinality of $\bK\langle d_0\rangle$ which is bounded by $|T|$. 
	\end{proof}
\end{proposition}
	
\begin{lemma}\label{lem:transl}
	Let $f: \bU \to \bU$ be an $\bE$-definable function such that $f(b) \in \cO_b$. Then there is a $\bK$-definable function $g$ such that
	\begin{equation}\tag{$*$}\label{eq:stdp}
        (\bE, \cO) \models \exists a_- \in \cO,\; \exists a_+>\cO,\;\forall t \in (a_-, a_+),\; f(t)-g(t) \in \co.
    \end{equation}
	Moreover for any such function $g$,
	\begin{enumerate}
		\item $f(b)-g(b)\in \co_b$; in particular if $g_1$ is another $\bK$-definable function satisfying (\ref{eq:stdp}), then $g_1(b)=g(b)$;
		\item $f(b) \in \cO_b^-$ if and only if $g$ is eventually bounded by some $c \in \bK$;
		\item $f(b)\in \co_b^-$ if and only if $g$ is infinitesimal at infinity.
	\end{enumerate}
	\begin{proof}
		Recall that $\tp(c/\bE)=\tp(b/\bE)$ for all $c \in \cO_b^{>\cO}$, so since $f$ is $\bE$-definable and $f(b) \in \cO_b$, $f(\cO_b^{>\cO}) \subseteq \cO_b$. Since $f$ is monotone on an $\bE$-definable interval containing $\cO_b^{>\cO}$ there is $a_- \in \bE$, $a_-<b$, such that for all $t \in \cO$ if $t> a_-$, then $|f(t)|< \bE^{>\cO_b}=\bE^{>\cO}$, that is, $f(t) \in \cO$. Without loss of generality we can also assume that $f$ is monotone continuous on $\cO^{>a_-}$.
			
		By Theorem~\ref{thm:std_def}, the function $g:= \std_{\bK} \circ f|_{\bK^{>a_-}}$ is $\bK$-definable and monotone continuous. Since $f|_{\cO^{>a_-}}$ is monotone continuous, $f(t)=g(t)+\varepsilon(t)$ with $\varepsilon$ an $\bE$-definable function such that $\varepsilon(t) \in \co$ for $t \in \cO^{>a_-}$. By Proposition~\ref{prop:coffact}, $\varepsilon$ cannot define an order-preserving bijection between a final segment of $\cO$ and a final segment of $\co$. So there must be $a_+>b$ such that $\varepsilon(t) \in\co$ for all $t\in (a_-, a_+)$.
			
		As for (1), to see that $\varepsilon(b) \in \co_b$ notice that otherwise $1/\varepsilon(b) \in \cO_b \setminus \cO_b^-$, but this would imply that for all big enough $t \in \cO$, $1/\varepsilon(t)\in \cO\setminus \co$, contradicting that $\varepsilon((a_-, a_+)) \subseteq \co$.
			
		(2) and (3) are clear.
	\end{proof} 
\end{lemma}
	
\begin{lemma}\label{lem:specialder}
	Let $f: \bU \to \bU$ be $\bE$-definable. If $f(b) \in \co_b$, then $f'(b) \in \co_b$. If $f(b) \in \co_b^-$ then $f'(b) \in \co_b^-$.
	\begin{proof}
		Assume $f(b) \in \cO_b$. By Lemma~\ref{lem:transl}, there is an $\bE$-definable interval $J\ni b$, on which $f=g+\varepsilon$ where $g$ is $\bK$-definable and $\varepsilon(J)\subseteq \co$. We can also assume that $\varepsilon$ is differentiable and monotone on $J$. If $I \subseteq J$ is an $\bE$-definable interval such that $\varepsilon'(I)\cap\co=\emptyset$, then $\sup I - \inf I\in \co$ by the mean value theorem.
        It follows that if $I\subseteq J$ is a convex subset with $\varepsilon'(I)\cap\co=\emptyset$, then there is $c \in J\cap \bE$, such that $I \subseteq c+\co$.
            
        Since by weak o-minimality $\{t \in J: \varepsilon'(t)\notin \co\}$ is a finite union of convex subsets and, by the previous observation, each of such convex subsets must be included in some translate of $\co$, we must have, after possibly restricting $J$, that $\varepsilon'(J)\subseteq \co$. 
			
		Since $\varepsilon(J) \subseteq \co$ and $\varepsilon'(J) \subseteq \co$, by Lemma~\ref{lem:transl}~(1) applied to $\varepsilon$ and to $\varepsilon'$,  $\varepsilon(b) \in \co_b$ and $\varepsilon'(b) \in \co_b$. 
			
		Therefore to prove the statement it suffices to observe that $g(b) \in \co_b^- \Rightarrow g'(b) \in \co_b^-$ and $g(b) \in \co_b \Rightarrow g'(b) =0$.
		Remembering $g$ is $\bK$-definable, the first implication holds by Lemma~\ref{lem:transl}~(3) and because the derivative of a function infinitesimal at infinity is again infinitesimal. The second holds because $\bK \langle b \rangle \cap \co_b =\{0\}$.
	\end{proof}
\end{lemma}
	
We record the following easy consequences of the mean value theorem for reference.
	
\begin{lemma}\label{lem:basicmvt}
	Let $f: \bE^{>0} \to \bE^{>0}$ be a monotone $\bE$-definable differentiable function, and $J\subseteq \bE^{>0}$ be a convex subset such that $\forall t \in J, \; t f^\dagger(t) \in \cO$. Then
	\begin{enumerate}
		\item $\forall x, y \in J,\; x \sim y \rightarrow f(x) \sim f(y)$;
		\item $\forall x, y \in J,\; x\asymp y \rightarrow f(x) \asymp f(y)$.
	\end{enumerate}
	Similarly for every differentiable $\bE$-definable $f: \bE \to \bE$ and every convex $J\subseteq \bE^{>0}$ such that $\forall t \in J, \; f'(t) \in \cO$
	\begin{enumerate}[align = left, leftmargin=*]
		\item[1bis.] $\forall x, y \in J,\; x- y \in \co \rightarrow f(x) - f(y) \in \co$;
		\item[2bis.] $\forall x, y \in J,\; x- y \in \cO \rightarrow f(x) - f(y) \in \cO$.
	\end{enumerate}
	\begin{proof}
		(1). For every positive $x \neq y$, by the mean value theorem, there is $t$ between $x$ and $y$ such that $f(x)-f(y)=f'(t)(x-y)$ so
		\[\frac{f(x)-f(y)}{f(t)}=f^\dagger(t)t \frac{x-y}{t}\]
		Now note that $x\sim y$ if and only if $(x-y)/t \in \co$, and since $f$ is monotone, $f(x) \sim f(y)$ if and only if $(f(x)-f(y))/f(t)\in \co$, thus we get (1).
        
        (2). We will use that $\cO$ is $T$-convex. Suppose toward contradiction that there are $x \in J$, $s \in \cO^{>\co}$ such that $xs \in J$ and $f(xs)\succ f(x)$.
        We will assume without loss of generality that $s>1+\co$, so $|f|$ is increasing. Consider the increasing function $g: [1,s] \to \bE$ given by $g(r)=|f(xr)/f(x)|$. 
        Note that by (1) $g(r+\co) \subseteq g(r)(1+\co)$ for all $1+\co <r<s+\co$ and $g(s+\co^{\le 0}) \subseteq g(s)(1+\co^{\le 0})$.
        Let $F := \{r \in [1,s]: g(r) \in \cO\}$ and note that $\res_{\cO}(F)$ is bounded but cannot have a supremum: in fact it cannot have a maximum, because if $\res_\cO(s_1)$ was such a maximum, then $g(s_1+\co)$ would be cofinal in $\cO$ contradicting Proposition~\ref{prop:coffact}, and on the other hand if $g(s_1)>\cO$, since $g(s_1+\co^{\le 0})\subseteq g(s_1)(1+\co^{\le 0})$, and $g$ is continuous increasing on $[1,s]$ there would be $s_2<s_1+\co$ such that $g(s_2+\co) \subseteq g(s_2)(1+\co)=(g(s_1)/2)(1+\co)>\cO$.

        But then $\cO$ would not be $T$-convex as if it was, then by Theorem~\ref{thm:std_def}, $\res_\cO(\bE,\cO)$ would be o-minimal and $\res_{\cO}(F)$ would have a supremum. Contradiction.
			
		(1bis)-(2bis) Straightforward from the mean value theorem.
	\end{proof}
\end{lemma}

\begin{lemma}\label{lem:wimpres}
	If for some $c\in \bE$, $(t-c)f'^\dagger(t) \in \cO$ for all $t$ in some interval $(a,b)$, then for every weakly $\cO$-immediate $x\in \bU \setminus \bE$ such that $a<x<b$, $f(x)$ is weakly $\cO$-immediate.
	\begin{proof}
		Let $(x_i)_{i<\lambda}$ be a p.c.-sequence in $(\bE,\cO)$ whose pseudolimit is $x$. 
        It suffices to show that $(f(x_{\mu+i}))_{i<\lambda}$ is a p.c.-sequence for some $\mu< \lambda$. Note that for $i<j<k< \lambda$, by the mean value theorem, there are $\xi_{ij}, \xi_{ik} \in \bE$ respectively between $x_i$ and $x_j$ and between $x_j$ and $x_k$, such that  
		\[\frac{f(x_j)- f(x_k)}{f(x_j)-f(x_i)}=\frac{f'(\xi_{jk})}{f'(\xi_{ij})}\frac{x_j-x_k}{x_j-x_i}.\]
		Observe that for big enough $\mu< \lambda$, for all $k>j>i>\mu$, $\xi_{ij}$ and $\xi_{jk}$ are forced to be in $(a,b)$, and to be such that $\xi_{ij} - c \asymp \xi_{jk} - c$. Now set $g(t):=f'(t+c)$, and observe that by hypothesis $tg^\dagger(t)\in \cO$ for all $t \in (a-c, b-c)$, so by Lemma~\ref{lem:basicmvt}, $g(\xi_{ij} -c) \asymp g(\xi_{jk}-c)$, and the thesis follows.
	\end{proof}
\end{lemma}

\section{Main results}\label{sec:mainres}
	
This section contains the proofs of Theorem~\ref{introthm:A} and Theorem~\ref{introthm:B}. We maintain the same setting and notation of Subsection~\ref{ssec:TechLem}. So $\bE\prec \bU \models T$ for some fixed complete o-minimal theory $T$ expanding $\RCF$, $\cO\subset \bE$ is a non-trivial $T$-convex subring, $b\in \bU \setminus \bE$ is a realization of the $\cO$-special cut of $\bE$, and $\bK\preceq \bE$ is maximal such that $\bK \subseteq \cO$.
	
\subsection{Unary types in \texorpdfstring{$T_\convex$}{T convex}}\label{ssec:unary_types}

The goal of this subsection is to establish Theorem~\ref{thm:types}, which will, in turn, justify Definition~\ref{def:types}. The proof of Theorem~\ref{thm:types}, requires proving that some kind of types cannot be realized in 1-$\dcl_T$-dimensional extensions obtained by adjoining a realization of a given type, or equivalently, that there are no definable functions mapping the cut of one type to the cut of the other. For one of these arguments (Proposition~\ref{prop:wim_main}) the following Lemma of technical nature will be used.

\begin{lemma}\label{lem:lem0}
	Let $C$ be a non-empty (not necessarily definable) convex subset of $\cO$ such that $\std_\bK(C)$ has no maximum and no minimum. Let $f$ be an $\bE$-definable continuous strictly monotone function defined on some interval whose image contains $C$. Suppose furthermore that on $f^{-1}(C)$: $f$ is 2-fold differentiable with $f'$ continuous strictly monotone and never-zero and $f'^\dagger$ is either constant non-zero or continuous and strictly monotone on $f^{-1}(C)$. If $f(x) \in C$ and $f^{-1}(C)-x$ is an $\cO$-submodule of $\bE$, then $\std_\bK(C)$ is an open interval.
	\begin{proof}
		Set $z:=f(x)$ and let us refer to the hypothesis that $f'^\dagger$ is strictly monotone and continuous on $f^{-1}(C)$ as hypothesis (A) and to the hypothesis that $f'^\dagger$ is constant and non-zero as hypothesis (B). Also let $\std:=\std_\bK$ to lighten the notation.
		
		Let us consider the maps $\mu(r,t):=f(x+r(f^{-1}(t)-x))$ and $\sigma(t_1,t_2) := f(f^{-1}(t_1)+f^{-1}(t_2)-x)$, given by the push-forwards by $f$ respectively of the scalar multiplication and sum of the $\cO$-submodule $f^{-1}(C)-x$. The idea is to first show that this \emph{almost} induces a $\bK$-definable ordered $\bK$-vector space structure on $\std(C)$. More specifically we will show that:
		\begin{enumerate}
		      \item $\std \sigma(\std(t_1),\std(t_2)) = \std(\sigma(t_1,t_2))$ for all $t_1, t_2 \in C$.
			\item for every $r\in \cO$ and every $t\in C$, $\std \mu(\std(r),t)=\std(\mu(r,t))$;
			\item there is a finite $F \subseteq C$ such that for all $t \in C \setminus (F+\co)$ and for all $r\in \cO \setminus \co$, $\std(\mu(r, \std(t)))= \std \mu(r,t)$.
		\end{enumerate}
		Hence $\sigma$ induces a group operation $\bar{\sigma}$ on $\std(C)$ and $\mu$ induces a map $\bar{\mu}: \bK^{\neq 0} \times (\std(C) \setminus \std(F)) \to \std(C)$.
			
		To see $(1)$ observe that by Lemma~\ref{lem:basicmvt}~(1bis) it suffices to show that for all $t_1, t_2 \in C$, $\partial_i \sigma (t_1,t_2) \in \cO$ for $i\in \{1,2\}$.
		First, set $R(s_1, s_2):=f'(s_1+s_2-x)/f'(s_1)$ and notice that
		\[\partial_1R(s_1, s_2) =\frac{\partial}{\partial s_1} \frac{f'(s_1+s_2-x)}{f'(s_1)} = \frac{f''(s_1+s_2-x)}{f'(s_1)} - \frac{f'(s_1+s_2-x)f''(s_1)}{f'(s_1)^2}\]
		vanishes if and only if $f'^\dagger(s_1+s_2-x)=f'^\dagger(s_1)$. So if we are in hypothesis (A), and $s_1, s_2 \in f^{-1}(C)$, $\partial_1 R(s_1, s_2)=0$ if and only if $s_2=x$. On the other hand under hypothesis (B), we have $\partial_1 R(s_1, s_2)=0$ for all $s_1, s_2 \in f^{-1}(C)$.
			
		Observe that
		\[\partial_1 \sigma (t_1, t_2) = \frac{f'(f^{-1}(t_1) + f^{-1}(t_2) -x)}{f'(f^{-1}(t_1))}=R(f^{-1}(t_1), f^{-1}(t_2)),\]
		so, under the hypothesis (A), $\{t_1\in C: \partial_1^2 \sigma(t_1, t_2)= 0\}\neq \emptyset$ if and only if $t_2=z$, and under the hypothesis (B) $\partial_1^2\sigma(t_1, t_2)=0$ for all $t_1, t_2 \in C$.
			
		Under hypothesis (A) for each $t_2$ either $\partial_1\sigma (-,t_2)$ is strictly monotone on $C$, or $f^{-1}(t_2)=x$ and thus $\partial_1 \sigma (t_1, t_2)=1$ for all $t_1\in C$.
			
		Under hypothesis (B), instead, for each $t_2$, $\partial_1 \sigma(-, t_2)$ is constant on $C$.
			
		In either of the cases (A) or (B), for each $t_2$, $\partial_1 \sigma (-, t_2)$ is strictly monotone or constant on $C$.
			
		Since $\std(C)$ has no maximum nor minimum, this implies that $\{t\in C: \partial_1 \sigma(t, t_2)\notin \cO\}$ is either empty or contains an interval $(t_-, t_+)$ with $t_-+\co<t_+$. The latter option would imply, by the mean value theorem, that $\sigma((t_-, t_+) \times \{t_2\}) \not \subseteq \cO$, against the hypothesis that $\sigma(C \times C) \subseteq C$. So in either of the cases (A), (B), $\{t\in C: \partial_1 \sigma(t, t_2)\notin \cO\}$ is empty and $\partial_1 \sigma (C \times \{t_2\}) \subseteq \cO$.
			 
		Since this holds for all $t_2$, $\partial_1 \sigma (C\times C) \subseteq \cO$.
		By symmetry also $\partial_2 \sigma (C \times C) \subseteq \cO$ and (1) follows from Lemma~\ref{lem:basicmvt}~(1bis).
					
		Notice that this implies also that the standard part $\bar{\sigma}:=\std \sigma |_{\bK}$ of $\sigma$, gives an ordered Abelian group structure on $\bK \cap C$ with $0$-element $\std(z)$, in particular for every $t\in C\setminus (z+\co)$, $\sigma(t,t)-t \notin \co$.
			
		To prove (2) and (3) we will only use the hypothesis that $f$ and $f'$ are strictly monotone continuous on $f^{-1}(C)$ so we won't need to distinguish between (A) and (B). Observe that 
		\[\begin{aligned}
			\partial_1\mu (r, t) = f'(r(f^{-1}(t)-x)+x) \cdot (f^{-1}(t)-x), \\
			h(r,t):=\frac{1}{r} \partial_2\mu(r,t) = \frac{f'(r(f^{-1}(t)-x)+x)}{f'(f^{-1}(t))}.
		\end{aligned}\]
		Therefore by the hypothesis on $f$, $\partial_1\mu$ is monotone continuous in the first variable. It follows that for every $t\in C$ and $r\in \cO$, $\partial_1 \mu(r,t)\in \cO$, for otherwise $\partial_1 \mu(\cO^{\ge r},t)\cap \cO=\emptyset$ or $\partial_1\mu (\cO^{\le r},t) \cap \cO = \emptyset$, but both would imply by the mean value theorem that $\mu(\cO, t) \not \subseteq \cO$. This proves (2), again by Lemma~\ref{lem:basicmvt}~(1bis).
			
		Notice that for each $r$, $\std\{t \in C: \partial_2\mu(r,t)\notin \cO\}$ must be finite, because otherwise by weak o-minimality and by the mean value theorem one would get $\mu(r, C)\not \subseteq C$. So for each $r \notin \co$, $\std\{t \in C: h(r,t)\notin \cO\}$ is finite.
			
		On the other hand, $h$ is monotone in the first variable $r$ if restricted to $\cO^{>0}$ or to $\cO^{<0}$, so if for some $r,t$ with $r\notin \co$, $h(r,t)\notin \cO$, then $h(r_1,t)\notin \cO$ for all $r_1$ on one side of $r$ in $\cO^{>\co}$ or $\cO^{<\co}$.
			
		This implies that there is a finite $F\subseteq C$ such that for all $t \notin F+\co$, and all $r \in \cO \setminus \co$, $h(r,t)\in \cO$. Suppose not: by weak o-minimality the set $\{t \in C : \exists r\in \cO\setminus \co,\; h(r,t)\notin \cO\}$ would contain some closed interval $[t_-, t_+]$ with $t_-<t_+ + \co$ and $t_-, t_+ \in C$; by Theorem~\ref{thm:std_def}, $\std\{(r,t) \in (\cO \setminus \co) \times [t_-, t_+]: h(r,t) \notin \cO\}$ is $\bK$-definable so by definable choice there would be a $\bK$-definable function $\gamma : [t_-, t_+]\cap \bK \to \bK \setminus \{0\}$ with $h(\gamma(t)+\co, t+\co) \not\subseteq \cO$ for all $t \in [t_-, t_+]\cap \bK$. Without loss of generality, we can assume $\gamma ([t_-, t_+]\cap \bK) \subseteq \bK^{>0}$ and that $\gamma$ is continuous, hence it has a maximum $r_+$ and a minimum $r_-$. 

        By monotonicity of $h$ in the first variable, at least one among $\std\{t \in C: h(r_+ + 1 ,t)\notin \cO\}$ and $\std\{t \in C : h(r_-/2,t) \notin \cO\}$ contains $[t_-, t_+] \cap \bK$, so one of the two would be infinite, contradiction.
        
		Thus for all $t \notin F+\co$ and all $r \notin \co$, $\partial_2 \mu(r,t) \in \cO$. By Lemma~\ref{lem:basicmvt}~(1bis) we again deduce (3). 
			
		Now we prove the statement of the Lemma assuming for the sake of simplicity that $f$ is increasing on $f^{-1}(C)$. The case in which $f$ is decreasing is analogous.
			
		Since we assumed $\std(C)$ has no maximum, there is $t_0\in C$ with $C^{\ge t_0} \cap \big((F\cup \{z\})+\co \big)=\emptyset$. In particular for every $t\ge t_0>z+\co$, $\mu(2,t)=\sigma(t, t)>t+\co$.
			
		It follows that $\partial_1\mu(r, t_0)> \co$ for every large enough $r \in \cO$: if not, by weak o-minimality, there would be $r_0 \in \cO^{\ge 1}$ such that $\partial_1 \mu(\cO^{\ge r_0},t) \subseteq \co$, but this is absurd because then the mean value theorem would give $\mu(\cO^{\ge r_0}, t_0) \subseteq \mu(r_0, t_0)+\co$ whereas $\mu(r_0, t_0)\ge t_0$ would give 
		\[\mu(2r_0, t_0)=\mu(2, \mu(r_0, t_0))=\sigma(\mu(r_0, t_0), \mu(r_0, t_0))> \mu(r_0, t_0)+\co,\]
		contradiction.
			
		Since $\partial_1\mu(r, t_0)> \co$, for every large enough $r \in \cO$, some final segment of $\std(\mu(\cO^{>\co},t_0))$ is an open interval in $\bK$.

        Now if $C$ is bounded above in $\cO$, we can consider $t_1\coloneqq \sup_\bK \std(\mu(\cO^{>\co},t_0))\in \bK$. We must then have $t_1 \notin C$, because otherwise $t_2:=\mu(1/2, t_1)>\mu(\cO^{>\co},t_0)$, and by construction $t_2+\co< \sigma(t_2,t_2)=t_1$ contradicting the fact that $t_1=\sup \std(\mu(\cO^{>\co},t_0))$. Since $\std(C)$ is order-convex, it follows that then $t_1 > \std(C)$ and thus $t_1=\sup_{\bK} \std(C)$.

		A similar argument shows that $\std(C)$ is either unbounded below in $\bK$ or has an infimum in $\bK$. Thus $\std(C)$ must be an open interval.
	\end{proof}
\end{lemma}

\begin{addendum*}
    The following Proposition can also be derived as a consequence of \cite[Thm.~5.2]{tressl2005model}, which at the time of writing was unknown to the author. Here I give a proof based on Lemma~\ref{lem:lem0}.
\end{addendum*}

\begin{proposition}\label{prop:wim_main}
	Let $\bU\succ \bE$ and $(\bE, \cO)\models T_\convex$. If $x \in \bU \setminus \bE$ is weakly $\cO$-immediate, then $\bE \langle x \rangle$ does not contain any $b$ with $\cO < b < \bE^{>\cO}$.
	\begin{proof}
		Let $b$ be $\cO < b < \bE^{>\cO}$ and assume toward a contradiction that $\bE\langle x \rangle = \bE \langle b \rangle$. So there is an $\bE$-definable continuous monotone bijection $f: \bU \to \bU$ such that $f(x)=b$. Up to replacing $x$ with $-x$ which is still weakly $\cO$-immediate, we can also assume without loss of generality that $f$ is increasing. This means that $f$ must induce a continuous increasing bijection between $\tp(x/\bE)(\bE \langle x \rangle)$ and the set $C\coloneqq \tp(b/\bE)(\bE \langle x \rangle)=\cO_x^{>\cO}$. Since $x$ is weakly $\cO$-immediate, by Lemma~\ref{lem:br_Omod}, $M\coloneqq \Br(x/\bE)(\bE \langle x \rangle)$ is an $\cO_x$-submodule and $\tp(x/\bE)(\bE \langle x \rangle)=x+M$.

        Now let $\bK_x$ be maximal in $\{\bK'\preceq \bE \langle x\rangle: \bK \preceq \bK' \subseteq \cO_x\}$.
            
		If the hypotheses of Lemma~\ref{lem:lem0} apply to $f$, then $\std_{\bK_x}(C)$ should be an interval in $\bK_x$ which it is not, contradiction.
			
		If the hypotheses of Lemma~\ref{lem:lem0} do not apply, then $f'^\dagger=0$ on $x+M$, but then for some $c_0 \in \bE$, $c_1 \in \bE^{\neq 0}$, $f(t)=c_0+c_1t$ for all $t$ in some $\bE$-definable interval containing $x$. This would imply $x=f^{-1}(b)= (b - c_0)/c_1$, which contradicts the fact that $x$ is weakly $\cO$-immediate.
	\end{proof}
\end{proposition}
	
\begin{proposition}\label{prop:wim_cof}
	Let $\bU \succ \bE$, $(\bU, \cO) \models T_\convex$, $x \in \bU \setminus \bE$ weakly $\cO\cap \bE$-immediate over $\bE$, $\bK \prec \bE$ maximal in $\{\bK' \preceq \bE: \bK' \subseteq \cO\}$ and $\bK_x \preceq \bE\langle x \rangle$ maximal in $\{\bK' \preceq \bE\langle x \rangle: \bK \preceq \bK'\subseteq \cO\}$. Then $\bK=\bK_x$.
	\begin{proof}
		By Lemma~\ref{lem:resrank} we have $\bK_x=\bK \langle f(x) \rangle$ for some $\bE$-definable $f$ and, by Proposition~\ref{prop:wim_main}, $\bE\langle x \rangle$ does not realize the special cut, so $\bK$ is cofinal in $\bK_x$ and $\tp(f(x)/\bK) \vdash \tp(f(x)/\bE)$.
			
		By Lemma~\ref{lem:wim_pc}, there are $(x_i^+)_{i< \lambda}$ and $(x_i^-)_{i<\lambda}$ $(\cO\cap \bE)$-p.c.-sequences for $x$ in $\bE$ with $x_i^+>x$ decreasing and $x_i^-<x$ increasing. Assume without loss of generality that $f$ is differentiable with $f'>0$ and set $z_i^+\coloneqq f(x_i^+)$, $z_i^-\coloneqq f(x_i^-)$, and $z\coloneqq f(x)$.
			
		Now suppose toward contradiction that $\bK \neq \bK_x$, so $f(x) \notin \bK$. Then up to extracting a subsequence we can assume that $z_i^+ - z_{i+1}^+ \asymp 1 \asymp z_i^- - z_{i+1}^-$. Notice that 
		\[x_{i+1}^+ - x_i^+  = f^{-1}(z_{i+1}^+) -  f^{-1}(z_{i}^+) = (f^{-1})'(\zeta_i^+) (z_{i+1}^+-z_{i}^+) \asymp (f^{-1})'(\zeta_i^+)\]
		for some $\zeta_i^+$ with $z_{i}^+<\zeta_i^+<z_{i+1}^+$. The sequence $(\zeta_i^+)_{i<\lambda}$ is decreasing and coinitial in $\bE^{>z}$, and $\val \big((f^{-1})' (\zeta_i^+)\big)$ is an increasing sequence because $(x_i^+)_{i< \lambda}$ is a p.c.-sequence. Similarly there is a sequence $(\zeta_i^-)_{i<\lambda}$ increasing and cofinal in $\bE^{<z}$ with $\val \big( (f^{-1})' (\zeta_i^-) \big)$ again increasing. But this would imply that $(f^{-1})'$ is decreasing on a final segment $\bE^{<z}$ and increasing on an initial segment of $\bE^{>z}$, which contradicts the hypothesis $z \notin \bE$. 
	\end{proof}
\end{proposition}
	
\begin{theorem}\label{thm:types}
	If $x \in \bU \setminus \bE$, and $\cO$ is a non-trivial $T$-convex subring of $\bE$, then the following are mutually exclusive
	\begin{enumerate}
		\item $\bE\langle x \rangle =\bE\langle y \rangle$ for some weakly $\cO$-immediate $y$;
		\item $\bE \langle x \rangle = \bE \langle b \rangle$ for some $b$ such that $\cO < b <\bE^{>\cO}$;
		\item $\bE \langle x \rangle = \bE \langle z \rangle$ for some $z$ such that for some (equiv.\ any) $\bK \preceq \bE$ maximal in $\{\bK \preceq \bE: \bK \subseteq \cO\}$, $\bK \langle z \rangle$ is a cofinal extension of $\bK$;
		\item $\bE\langle x \rangle=\bE \langle d \rangle$ for some $d>\bE$.
	\end{enumerate}
	\begin{proof}
		The fact that (1) implies that none of (2) or (3) holds is given respectively by Propositions~\ref{prop:wim_main} and \ref{prop:wim_cof}. To see that (2) and (3) are mutually exclusive, let $\cO_x:=\{t \in \bE\langle x \rangle: |t|<\bE^{>\cO}\}$, $\bK$ be as in $(3)$, and $\bK_x$ be maximal in $\{\bK'\preceq \bE \langle x \rangle: \bK \prec \bK'\subseteq \cO_x\}$. By Lemma~\ref{lem:resrank}, $\bK_x=\bK \langle z \rangle$ for some $z \in \bE \langle x \rangle\setminus \bE$: only one among (3) and (2) can occur because either $\bK$ is cofinal in $\bK_x$, or there is $b> \bK$ such that $\bK_x=\bK\langle b \rangle$, in which case $\cO=\CH_{\bE}(\bK)< b \in \bK_x <\bE^{>\cO}$.
			
		Finally (4) clearly implies that none of (1), (2), or (3) holds because if $(4)$ holds, then for every $h \in \CH_{\bE \langle d \rangle}(\bE)$, there is $c \in \bE$ such that $|c-h|<\bE^{>0}$, which is contradicted by the choice $h=y$ for (1), $h=b$ for (2), and $h=z$ for (3).
	\end{proof}
\end{theorem}
	
\begin{definition}\label{def:types}
	If $x$ satisfies $(1)$, then we call $x$ \emph{weakly immediately generated} with respect to $\cO$ (or \emph{weakly $\cO$-immediately generated}).
		
	If $x$ is in case $(2)$ we call it \emph{$\cO$-special}.
		
	If $x$ is in case $(3)$ we call it \emph{$\cO$-cofinally residual}.
		
	If $x$ is as in $(4)$ then we say $x$ is \emph{tame} (because in this case $\bE\langle x \rangle \succ_\tame \bE$).
		
	If $x \in \bU \setminus \bE$ is not in the cases above, then we call $x$ \emph{strictly purely $\cO$-valuational}.
		
	We say that $x$ is \emph{dense} over $\bE$ if $\bE$ is order dense in $\bE \langle x \rangle$.
\end{definition}
	
\begin{remark}
	If $\bE$ admits a non-trivial $T$-convex valuation ring $\cO$, then $x$ is dense if and only if $x$ is weakly $\cO$-immediate and $\Br(x/\bE)(\bE)=0$.
\end{remark}
	
\begin{remark}\label{rmk:strictpv}
	An element $x \in \bU \setminus \bE$ is strictly purely $\cO$-valuational if and only if $\val(\bE\langle x \rangle, \cO_x)$ is a cofinal extension of $\val(\bE, \cO)$ that does not realize the type $0^+$ and moreover for every $y \in \bE \langle x \rangle\setminus \bE$, there is $c \in \bE$ such that $\val_{\cO}(y-c) \notin \val_{\cO} \bE$.
\end{remark}
	
Theorem~\ref{thm:types} can be restated in terms of unary types over $(\bE, \cO)$. To do that it is convenient to introduce some terminology for the types over $(\bE, \cO) \models T_\convex$ as well.
	
\begin{definition}\label{def:types1}
	Let $(\bE, \cO)\models T_\convex$, $(\bU, \cO')\succ (\bE, \cO)$, and $x \in \bU \setminus \bE$.
	We say that $x$ is
	\begin{enumerate}
		\item \emph{weakly immediately generated over $(\bE, \cO)$} if there is $y \in \bE \langle x \rangle\setminus \bE$, such that for every $c_0 \in \bE$ there is $c_1 \in \bE$ with $y-c_0 \prec y-c_1$.
		\item \emph{residual} over $(\bE, \cO)$ if $\res_{\cO'}(\bE\langle x \rangle)\neq \res_{\cO'}(\bE)$;
		\item \emph{purely valuational} if there is $M\subseteq \bE \langle x \rangle$ such that $\val(M) \cap \val (\bE)= \emptyset$ and $\bE \langle x \rangle \setminus \bE = \bE + M$.
	\end{enumerate}
	Notice that whether $x$ satisfies $(1)$, $(2)$, or $(3)$ only depends on $\tp(x/(\bE,\cO))$.
    Also note that in (3) $M$ is just given by $\{y \in \bE \langle x \rangle: \val(y) \notin \val(\bE)\}$.
\end{definition}
	
\begin{remark}
	The terminology given in Definition~\ref{def:types1} is consistent with the one given in Definition~\ref{def:types}, more specifically:
	\begin{enumerate}
		\item $x$ is weakly $\cO$-immediately generated if and only if $x$ is weakly immediately generated over $(\bE, \cO)$;
		\item if $x$ is $\cO$-special, then it is residual (but not cofinally) when $\cO'\cap \bE \langle x \rangle =\{t\in \bE \langle x \rangle : |t|<\bE^{>\cO}\}$, and it is purely valuational (but not strictly) when $\cO' \cap \bE \langle x \rangle = \CH_{\bE \langle x \rangle} (\cO)$;
		\item $x$ is cofinally residual if and only if $x$ is residual over $(\bE, \cO)$ and the underlying (proper) extension of residue fields is cofinal;
		\item if $x$ is tame, then it is purely valuational (but not strictly);
		\item $x$ is strictly purely $\cO$-valuational if and only if $x$ is purely valuational and $\val (\bE, \cO)^{>0}$ is both coinitial and cofinal in $\val(\bE\langle x \rangle, \cO')^{>0}$ (Remark~\ref{rmk:strictpv}).
	\end{enumerate}
\end{remark}
	
With this Theorem~\ref{thm:types} can be restated as
	
\begin{theorem}\label{thm:types_rest}
	Let $(\bU, \cO)\models T_\convex$, $\bE \prec \bU$ with $\bE\not \subseteq \cO$, and $x \in \bU \setminus \bE$. Then $x$ is either residual, or weakly immediately generated, or purely valuational.
\end{theorem}

\subsection{Wim-constructible extensions}\label{ssec:wimext}
In this subsection we define wim-constructible extensions and prove an amalgamation result (Lemma~\ref{lem:wim_amalgamate}). This, after some extra smallness assumptions on the extensions are imposed, allows for consideration of spherical completions of models.
		
\begin{remark}
	Recall that any extension $\bE'\succ \bE$ containing no $\cO$-special elements admits a unique $T$-convex extension $\cO'$ of $\cO$ to $\bE'$.
\end{remark}
	
\begin{remark}
	Let $(\bE', \cO')\succeq (\bE, \cO)\models T_\convex$. Given a valuation ball $x+\cO'y$ of $(\bE', \cO')$, $(x+\cO'y) \cap \bE$ is an intersection of valuation balls in $(\bE, \cO)$: in fact either $(x+\cO'y) \cap \bE=\emptyset$ or, without loss of generality, we can assume that $x \in \bE$, and the assertion follows from the fact that $\cO' y \cap \bE = \bigcap\{\cO z: z \in \bE, z \succeq y\}$. 
		
	Notice also that if $\bE'\succ \bE$ contains no $\cO$-special elements, then whenever $(x+\cO'y)\cap \bE \neq \emptyset$, the convex partial type defining $x+\cO'y$ is the unique extension of the partial type of $(x+\cO'y) \cap \bE$ to $\bE'$ both with respect to $T$ and to $T_\convex$.
\end{remark}
	
\begin{definition}\label{defn:wim_cof}
	If $(\bE, \cO) \models T_\convex$ and $p$ is a weakly immediate type over $(\bE, \cO)$, then we define the \emph{cofinality} of $p$ as the cofinality of $\val_{\cO'}(x-\bE)$ for $x \in p(\bU)$, $\bU \succ \bE$ and $\cO'$ any extension of $\cO$ to $\bU$, or equivalently the least cardinal $\lambda$ such that there is a p.c.\ sequence $(x_i)_{i \in \lambda}$ for $x$.
    We call an $x$ and a $p$ as above \emph{$\lambda$-bounded} if they have cofinality $<\lambda$.
\end{definition}
		
\begin{remark}\label{rmk:cof_equiv}
	If $y \in \bE \langle x \rangle \setminus \bE$ and $x,y$ are weakly $\cO$-immediate, then the types of $x$ and $y$ have the same cofinality. To see this, notice that the cofinality of $x$ is in fact the cofinality of $\bE^{<x}$,  and of $-(\bE^{>x})$. Then recall that if $y \in \bE \langle x \rangle \setminus \bE$, there is a monotone $\bE$-definable bijection $f: \bE \to \bE$ such that $y=f(x)$.
\end{remark}

\begin{definition}\label{defn:wimc_extension}
	Let $(\bE, \cO)\preceq (\bE_*, \cO_*)$ be an elementary extension of models of $T_\convex$.
	A sequence $(x_i: i\in I)$ in $\bE_*$, with $(I,<)$ a well ordered set and such that for every $j$, $x_j$ is weakly immediate over $\bE_j:=\bE \langle x_i: i < j\rangle$, is called a \emph{wim-construction} (over $(\bE, \cO)$). A wim-construction $(x_i: i\in I)$ in $\bE_*$ is said to be \emph{$\lambda$-bounded} if every $\tp(x_j/(\bE_j, \cO_* \cap \bE_j))$ has cofinality $<\lambda$. 
	We say that $(\bE_*, \cO_*)$ is:
	\begin{enumerate}
		\item \emph{wim-constructible} if there is a wim-construction $(x_i:i<\mu)$ in $\bE_*$ such that $\bE_*=\bE \langle x_i :i<\mu \rangle$;
		\item \emph{$\lambda$-bounded wim-constructible} if there is a $\lambda$-bounded wim-construction $(x_i:i<\mu)$ in $\bE_*$ such that $\bE_*=\bE \langle x_i :i<\mu \rangle$;
		\item \emph{strictly wim} if it has no non-weakly immediate factor, that is if every $\bE_1$ with $\bE \preceq \bE_1 \prec \bE_*$ is such that every $x \in \bE_* \setminus \bE_1$ is weakly immediate over $\bE_1$.
		\item $\kappa$-\emph{wim} for $\kappa$ a (possibly finite) cardinal if for every ordinal $\mu < \kappa$, the definable closure of every $(\mu+1)$-tuple in $\bE_*$ is wim-constructible over $(\bE, \cO)$.
	\end{enumerate}
	A (possibly infinitary) type over $(\bE, \cO)$ is said to be (\emph{$\lambda$-bounded}) \emph{wim-constructible} if it is the type of a ($\lambda$-bounded) wim-construction over $(\bE, \cO)$.
\end{definition}
	
\begin{remark}
	Wim-constructible types over $(\bE, \cO)$ are uniquely determined by their reduct to $\bE$, because, by Proposition~\ref{prop:wim_main}, they don't realize the special cut.
\end{remark}

\begin{remark}\label{rmk:wimc-prserve}
    By Theorem~\ref{thm:types}, wim-constructible extensions do not enlarge the residue field sort and preserve the cofinalities of the field and of the valuation ideal.
\end{remark}
	
\begin{remark}
	Every strictly wim extension is wim-constructible.
	If, for a given $T$, strictly wim extensions are transitive (i.e.\ closed under composition) then the converse holds as well. This is the case for example when $T$ is power-bounded (see Section~\ref{sec:pbdd}, in particular Corollary~\ref{cor:pbddwimc_im}).
\end{remark}
	
\begin{lemma}\label{lem:wimstat}
	Let $(\bU, \cO') \models T_\convex$, $\bE \preceq \bE_1 \prec \bU$ with $\bE \not \subseteq \cO'$ and $x \in \bU \setminus \bE$. Set $\cO\coloneqq\cO'\cap \bE$ and $\cO_1\coloneqq \cO'\cap \bE_1$. Suppose $\tp(x/\bE)(\bE_1)=\emptyset$. Then $x$ is weakly immediate over $(\bE, \cO)$ if and only if it is weakly immediate over $(\bE_1, \cO_1)$. If this is the case, then $\tp(x/(\bE,\cO))$ and $\tp(x/(\bE_1, \cO_1))$ have the same cofinality.
	\begin{proof}
		Let $p=\tp(x/\bE)$.
		If $x$ is weakly $\cO$-immediate over $\bE$ with cofinality $\kappa$, there is a p.c.-sequence $(x_i)_{i< \kappa}$ in $\bE$ such that 
		\[p(\bU)=\{x: x-x_i \succ x-x_{j}: i<j< \kappa\}\]
		If $p(\bU)$ has empty intersection with $\bE_1$, it defines a weakly immediate type over $(\bE_1, \cO_1)$ with cofinality $\le \kappa$.
			
		Conversely, suppose $x$ is weakly immediate over $(\bE_1, \cO_1)$ and let $p_1:=\tp(x/\bE_1)$ and $\kappa_1$ be its cofinality. Then $p_1(\bU)$ is an intersection of valuation balls $(B_i)_{i< \kappa_1}$ with center and radius in $\bE_1$ and w.l.o.g.\ $B_j \subset B_i$ for all $j < i$.
		For every $i$, $B_i \cap \bE$ is either empty or an intersection of valuation balls in $(\bE, \cO)$.
        If $B_i \cap \bE = \emptyset$ for some $i$, then for all $y \in B_i$, $\tp(y/\bE)=p$ and $p(\bE_1)\neq \emptyset$ contradicting the hypothesis that $p(\bE_1)=\emptyset$.
        Therefore, $B_i\cap \bE\neq \emptyset$ for all $i< \kappa_1$.
			
		Since $\kappa_1$ is regular, up to extracting a subsequence we can assume that $B_j \cap \bE \subset B_i \cap \bE$ for all $i<j$. But then for each $i$, there is a valuation ball $B'_i$ with radius and center in $\bE$ such that $B_{i+1} \subseteq B'_i \subseteq B_{i}$ and the cofinality of $p$ is $\le \kappa_1$.	
	\end{proof}
\end{lemma}
	
Recall that two types $p$, $q$ over $\bE$ are said to be \emph{weakly orthogonal} (cf \cite[Def.~3.12]{mennuni2020product}) if, given \emph{disjoint} tuples of variables $\overline{\var{x}}$, $\overline{\var{y}}$ of the appropriate length, $p(\overline{\var{x}})\cup q(\overline{\var{y}})$ defines a complete type over $\bE$. In particular, if $a \in \bE_*\setminus \bE$ and $\overline{b} \in \bE_*^n$, then $\tp(a/\bE)$ and $\tp(\overline{b}/\bE)$ are weakly orthogonal if and only if $\tp(a/\bE)$ is not realized in $\bE \langle \overline{b}\rangle$.
	
In the following lemma given two ordered sets $(A,<)$, $(B, <)$ we denote by $A\sqcup B$ the disjoint union of $A$ and $B$ with the smallest \emph{partial} order extending the orders on $A$ and $B$. 
	
\begin{lemma}\label{lem:toname00}
	Let $(\bE, \cO) \models T_\convex$, $(I,<)$, $(J,<)$ be well-orders and $(x_i: i\in I)$ and $(y_j : j\in J)$ be $\kappa$-bounded wim-constructions over $(\bE, \cO)$. If their types are weakly orthogonal over $\bE$, then for every bijective order preserving $h: I \sqcup J \to \beta$ to an ordinal $\beta$, the sequence $(z_k: k< \beta)$ defined by $z_{k} = x_{h^{-1}(k)}$ for $k \in h(I)$ and $z_{k} = y_{h^{-1}(k)}$ for $k \in h(J)$ is a $\kappa$-bounded wim-construction over $\bE$.
	\begin{proof}
		We need to show that $z_k$ is wim over $\bE \langle z_l: l<k\rangle$ with cofinality $< \kappa$. Suppose without loss of generality that $k \in h(I)$, thus $z_k$ is wim over $\bE \langle z_{l}: l\in h(I)^{<k}\rangle$ with cofinality $< \kappa$.
			
		From the orthogonality hypothesis we deduce that $\tp\big( z_k / \bE \langle z_l: l \in h(I)^{<k} \rangle \big)$ is weakly orthogonal to $\tp\big((z_{h(j)}: j \in J) /\bE \langle z_l: l \in h(I)^{<k} \rangle\big)$, so it is not realized in $\bE \langle z_l:l<k\rangle$. Thus we can conclude by Lemma~\ref{lem:wimstat}.
	\end{proof}
\end{lemma}
	
\begin{lemma}\label{lem:wimcont}
	Let $(\bU, \cO)\models T_\convex$, $(\bE_i)_{i< \lambda}$ be an 
    increasing sequence of elementary substructures of $\bU$, and $\bE_\lambda := \bigcup_{i< \lambda} \bE_i$. If $y\in \bU$ is weakly immediate over $\bE_{i}$ for all $i<\lambda$ then $y$ is weakly immediate over $\bE_\lambda$.
	\begin{proof}
		If $\lambda$ is not a limit ordinal the statement is trivial. So assume $\lambda$ is a limit ordinal.

        We will first prove the statement under the additional hypothesis that $(\bE_i)_{i<\lambda}$ is strictly increasing and continuous, i.e.\ such that $\bE_j=\bigcup_{i<j} \bE_{i+1}$ for all $j<\lambda$.
			
		If for some $j<\lambda$, $\tp(y/\bE_j)$ is not realized in $\bE_\lambda$, then we can invoke Lemma~\ref{lem:wimstat} and conclude.
			
		Otherwise we may suppose that for every $j<\lambda$, there is some $h(j) \ge j$ such that $\tp(y/\bE_{j})(\bE_{h(j)+1})\neq \emptyset$ and $\tp(y/\bE_{j})(\bE_{h(j)}) = \emptyset$. Build a sequence $(y_j)_{j<\lambda}$ with $y_j \in \tp(y/\bE_{j})(\bE_{h(j)+1})$, so in particular $y_j \notin \bE_{h(j)}\supseteq \bE_j$. 
			
		Notice that $\cO(y-y_j) \subseteq \Br(y/\bE_{h(j)})(\bE_\lambda)$, in particular $\big(y_j+\cO(y-y_j)\big) \cap \bE_{h(j)} = \emptyset$. So for every $\beta$, given $\gamma$ with $h(\gamma)>\beta$, 
		\[\big(y_\gamma + (y-y_\gamma) \cO\big) \cap \bE_\beta \subseteq \big(y_\gamma + (y-y_\gamma) \cO\big) \cap \bE_{h(\gamma)}=\emptyset.\]
		It follows that the sequence of valuation balls $(y_\gamma + (y-y_\gamma)\cO)\cap \bE_\lambda$ has empty intersection in $\bE_\lambda$ and $y$ is weakly immediate over $\bE_\lambda$.

        To get rid of the additional hypothesis, note that given any increasing $(\bE_i)_{i<\lambda}$ we can build an increasing sequence $(\bE_i')_{i\le\mu}$ of elementary substructures of $\bE_\lambda$ inductively by setting $\bE'_0=\bE_0$, $\bE_{i+1}'=\min \{\bE_j: \bE_j\supsetneq \bE_i'\}$ and $\bE_{i}'=\bigcup_{j<i}\bE'_j$ for limit $i$, stopping at $i=\mu$ when $\bE_i=\bE_\lambda$. Then, $y$ is weakly immediate over each $\bE_{i+1}'$ because $\bE_{i+1}'=\bE_j$ for some $j<\lambda$, and it will also be weakly immediate over $\bE_i'$ for every limit $i$: indeed, if this wasn't true then there would be a least limit ordinal $i$ such that $y$ is not weakly immediate over $\bE_i'=\bigcup_{j<i} \bE_j'$, we would then have that $(\bE_j')_{j<i}$ is an increasing sequence such that $y$ is weakly immediate over each $\bE_j'$ for $j<i$ and satisfies the additional hypothesis of being strictly increasing and continuous, thus contradicting what was already proven.
	\end{proof}
\end{lemma}

We are now ready to prove the main technical Lemma of this section. This is an amalgamation result for wim-constructible extensions and will be the main ingredient of the proof of Theorem~\ref{introthm:B} (Corollary~\ref{thm:B} below). Sometimes in the remainder of this section we will use the notation $\varphi_*p$ for the image (or push-forward) of a type $p$ over some set of parameters $A$ by an elementary map $\varphi:A \to B$, i.e.\
\[\varphi_* p(\overline{\var{x}}) = \{\alpha(\overline{\var{x}}, \varphi \overline{a}): \alpha(\overline{\var{x}}, \overline{\var{y}}) \text{ formula without parameters},\; \overline{a} \in A^{\overline{\var{y}}}, p(\var{x})\vdash \alpha(\overline{\var{x}}, \overline{a})\}\]
which is a type over $B$ (see \cite[p.~50]{tent2012course}).
	
\begin{lemma}\label{lem:wim_amalgamate}
	Let $(\bE, \cO)\models T_\convex$ and let $\bE_0, \bE_1$ be $\kappa$-bounded wim-constructible extensions of $\bE$. Then there is a $\kappa$-bounded wim-constructible extension $\bE_2 \succeq \bE_1$ and an elementary embedding $j:\bE_0 \to \bE_2$ over $\bE$ such that $\bE_2$ is $\kappa$-bounded wim-constructible over $j\bE_0$.
	\begin{proof}
		Let $\bE_2$ be a maximal $\kappa$-bounded wim-constructible extension within some $\kappa+|\bE_0|^+$ saturated extension of $\bE_1$. Let $(y_i: i < \beta)$ be a $\kappa$-bounded wim-construction of $\bE_1$ over $\bE$ and $(y_i: i< \beta+\alpha)$ an extension of it to a $\kappa$-bounded wim-construction of $\bE_2$ over $\bE$.
			
		Take a $\kappa$-bounded wim-construction $(x_i: i< \lambda)$ for $\bE_0$. We inductively build a sequence $(z_j: j < \lambda)$ in $\bE_2$ and a decreasing sequence $(S_j)_{j\le \lambda}$ of subsets of $\alpha+\beta$ such that for all $j\le \lambda$, $\bigcap_{i < j}S_{i+1}=S_j$ and
		\begin{enumerate}[label=(\arabic*), ref=(\arabic*)]
			\item \label{wim_amalgamate_induction:1} $(z_i : i < j)$ is a $\kappa$-bounded wim-construction over $\bE$ inside $\bE_2$ such that $\varphi(x_i):=z_i$ for $i<j$ defines an elementary map over $\bE$ and thus an embedding $\varphi: \bE_0^j:=\bE \langle x_i:i<j\rangle \to \bE_2$;
			\item \label{wim_amalgamate_induction:2} $(y_i: i \in S_j)$ is a $\kappa$-bounded wim-construction for $\bE_2$ over $\varphi(\bE_0^j)$;
			\item \label{wim_amalgamate_induction:3} for every $k\le j$ and $l \in S_j$, $(\varphi\bE_0^k)\langle y_i : i \in S_k^{<l}\rangle \subseteq (\varphi \bE_0^j) \langle y_i : i \in S_j^{<l}\rangle$;
			\item \label{wim_amalgamate_induction:4} for every $k \le j$, and every $l \in S_j$, $\tp(y_{l}/(\varphi\bE_0^k)\langle y_i : i \in S_k^{<l}\rangle)$ is not realized in $(\varphi\bE_0^j)\langle y_i : i \in S_j^{<l}\rangle$.
		\end{enumerate}
			
		If $j=0$ then just set $S_0=\alpha+\beta$.
			
		Suppose $(z_i : i < j)$ and $S_j$ have been built. Notice $\varphi_*\tp(x_j/\bE_0^j)$ is realized in $\bE_2$, and that by \ref{wim_amalgamate_induction:2} of the inductive hypothesis $(y_i:i \in S_j)$ is a $\dcl_T$-basis of $\bE_2$ over $\varphi(\bE_0^j)$. So there is $\iota(j)\in S_j$, minimal such that $(\varphi\bE_0^j)\langle y_i: i \le \iota(j) \rangle$ realizes $\varphi_*\tp(x_j/\bE_0^j)$. Set $z_j$ to be a realization of $\varphi_*\tp(x_j/\bE_0^j)$ in $(\varphi\bE_0^j)\langle y_i: i \le \iota(j) \rangle$ and set $S_{j+1}:=S_j\setminus \{\iota(j)\}$.
        Clearly such a choice satisfies \ref{wim_amalgamate_induction:1}.
        Moreover, by construction $\tp(z_j/\varphi(\bE_0^j))$ is weakly orthogonal to $\tp((y_i: i  \in S_{j}^{<\iota(j)})/\varphi(\bE_0^j))$, because $\iota(j)$ was chosen so that $\tp(z_j/\varphi(\bE_0^j))=\varphi_*(\tp(x_j/\bE_0^j))$ is not realized in $\varphi(\bE_0^j)\langle y_i: i \in S_j^{<\iota(j)}\rangle$.
        It thus follows from Lemma~\ref{lem:toname00} that $(y_i: i \in S_{j}^{<\iota(j)})$ is a $\kappa$-bounded wim-construction over $\varphi(\bE_0^{j+1})$. Since, by the exchange property, $(\varphi\bE_{0}^{j+1})\langle y_i: i \in S_{j}^{<\iota(j)}\rangle=(\varphi\bE_{0}^{j})\langle y_i: i \in S_{j}^{\le\iota(j)}\rangle$, it follows that $(y_i: i \in S_{j+1})$ is a $\kappa$-bounded wim-construction for $\bE_2$ over $\varphi\bE_{0}^{j+1}$. Thus \ref{wim_amalgamate_induction:2} holds too.
				
		As for \ref{wim_amalgamate_induction:3}, if $l \le \iota(j)$, then $S_{j+1}^{<l}=S_j^{<l}$ so trivially $(\varphi \bE_0^{j}) \langle y_i: i \in S_j^{<l}\rangle \subseteq (\varphi \bE_0^{j+1}) \langle y_i: i \in S_{j+1}^{<l}\rangle$.
		On the other hand if $l > \iota(j)$, then $(\varphi \bE_0^{j}) \langle y_i: i \in S_j^{<l}\rangle = (\varphi \bE_0^{j+1}) \langle y_i: i \in S_{j+1}^{<l}\rangle$ by the choice of $z_j$ and by the exchange property. So by inductive hypothesis if $k \le j$, then
		\[(\varphi \bE_0^{k}) \langle y_i: i \in S_k^{<l}\rangle \subseteq (\varphi \bE_0^{j}) \langle y_i: i \in S_j^{<l}\rangle \subseteq (\varphi \bE_0^{j+1}) \langle y_i: i \in S_{j+1}^{<l}\rangle.\]
			 
		As for \ref{wim_amalgamate_induction:4}, if $l \le \iota(j)$, by construction $\tp(z_j / \varphi \bE_0^j)$ is not realized in $(\varphi \bE_0^j)\langle y_i: i \in S_j^{<l}\rangle$ and $S_{j+1}^{<l}=S_{j}^{<l}$. Moreover if $l \in S_{j+1}$, then in fact $l< \iota(j) \notin S_{j+1}$. 
        
		It follows that $\tp(y_l/ (\varphi \bE_0^j)\langle y_i: i \in S_j^{<l}\rangle)$ is not realized in $(\varphi \bE_0^{j+1})\langle y_i: i \in S_j^{<l}\rangle$: for if it was and $y_l'$ was a realization then by the exchange property $z_j=f(y_l')$ for some $(\varphi \bE_0^j)\langle y_i: i \in S_j^{<l}\rangle$-definable function $f$, but then $f(y_l) \in (\varphi \bE_0^j)\langle y_i: i \in S_j^{\le l}\rangle$ would realize $\tp(z_j / (\varphi \bE_0^j)\langle y_i: i \in S_j^{< l}\rangle)\vdash \tp(z_j / \varphi \bE_0^j)$, and $l=\iota(j)\notin S_{j+1}$, contradiction.

        On the other hand if $l> \iota(j)$ then by construction $(\varphi \bE_0^{j}) \langle y_i: i \in S_j^{<l}\rangle = (\varphi \bE_0^{j+1}) \langle y_i: i \in S_{j+1}^{<l}\rangle$ so for every $l\in S_{j+1}$, $\tp(y_l/(\varphi\bE_0^{j})\langle y_i: i \in S_{j}^{<l}\rangle)$ is not realized in $(\varphi\bE_0^{j+1})\langle y_i: i \in S_{j+1}^{<l}\rangle)$.
        Finally, by inductive hypothesis if $k \le j$, then $\tp(y_l/(\varphi \bE_0^k) \langle y_i: i \in S_k^{<l}\rangle)$ is not realized in $(\varphi \bE_0^j) \langle y_i: i \in S_j^{<l}\rangle$ and thus extends uniquely to $\tp(y_l / (\varphi \bE_0^j) \langle y_i: i \in S_j^{<l}\rangle)$ which in turn is not realized in  $(\varphi \bE_0^{j+1}) \langle y_i: i \in S_{j+1}^{<l}\rangle$. Thus (4) is preserved at successor steps.
			
		At limit steps instead we set $S_{j}=\bigcap_{k<j}S_{k+1}$ and take the union of the sequences $(z_i)_{i<k+1}$ for $k<j$.
		To see \ref{wim_amalgamate_induction:3} we will show that for $l \in S_j$, 
		\begin{equation}\label{eq:2} \tag{$*$}
		(\varphi\bE_0^j) \langle y_i: i \in S_{j}^{<l}\rangle = \bigcup_{k<j} (\varphi\bE_0^k) \langle y_i: i \in S_{k}^{<l}\rangle.
		\end{equation}
		In fact, for every $k<j$, $z_k \in (\varphi\bE_0^{k+1}) \langle y_i: i \in S_{k}^{<l} \rangle$ and for every $i \in S_j$, $i \in S_{k}$ for every $k<j$ which gives the containment $\subseteq$. For the other containment, it suffices to show that for every $m \in (S_0\setminus S_j)^{<l}$, $y_m \in (\varphi\bE_0^j) \langle y_i: i \in S_{j}^{<l}\rangle$. 
			
		Toward contradiction let $m \in (S_0 \setminus S_j)^{<l}$ be minimum such that $y_m \notin (\varphi\bE_0^j) \langle y_i: i \in S_{j}^{<l}\rangle$ and let $k<j$ be such that $y_m \in S_{k} \setminus S_{k+1}$. By construction $y_m \in (\varphi \bE_0^{k+1}) \langle y_i: i \in S_{k+1}^{<m}\rangle$, moreover by the minimality of $m$, $S_{k+1}^{<m}=S_j^{<m}$ and thus
		\[y_m \in (\varphi \bE_0^{k+1}) \langle y_i: i \in S_{j}^{<m}\rangle \subseteq (\varphi \bE_0^j) \langle y_i : i \in S_{j}^{<l}\rangle,\]
		because $m<l$, contradiction.
			
		To see \ref{wim_amalgamate_induction:4} suppose $\tp(y_l/(\varphi\bE_0^k) \langle y_i: i \in S_{k}^{<l}\rangle)$ was realized in $(\varphi\bE_0^j) \langle y_i: i \in S_{j}^{<l}\rangle$. Then by (\ref{eq:2}) it would be realized in some $(\varphi\bE_0^{m}) \langle y_i: i \in S_{m}^{<l}\rangle$ for some $j>m>k$, contradicting the inductive hypothesis.

		To see \ref{wim_amalgamate_induction:2} notice that $\tp(y_l/(\varphi \bE_0^j)\langle y_i: i\in S_j^{<l}\rangle)$ is the unique extension of $\tp(y_l/\bE\langle y_i: i\in S_0^{<l}\rangle)$, which is $\kappa$-bounded weakly immediate so we can apply Lemma~\ref{lem:wimstat}.	
	\end{proof}
\end{lemma}
	
\begin{remark}
	If one wants to disregard the $\kappa$-bounded condition, the proof can be somewhat simplified by only showing (1)-(2)-(3) hold for the construction and using Lemma~\ref{lem:wimcont} and (3) to get (2) in the limit step.
\end{remark}

\begin{remark}\label{rmk:lambda_spherical_completeness}
    Recall that a valued field is said to be \emph{$\lambda$-spherically complete} if every finitely consistent intersection of fewer than $\lambda$ many valuation balls is non-empty. Thus a model $(\bE, \cO)\models T_\convex$ is $\lambda$-spherically complete if and only if for every $(\bU, \cO')\succ (\bE, \cO)$, and every $x \in \bU\setminus \bE$, wim over $(\bE, \cO)$, one has that $\val_{\cO'}(x-\bE)$ has cofinality greater or equal to $\lambda$.
    In other words $(\bE,\cO)$ is $\lambda$-spherically complete if and only if it has no proper $\lambda$-bounded wim-constructible extensions.
        
	Thus, a $\lambda$-bounded wim-constructible extension of $(\bE, \cO)$ is \emph{maximal} (i.e.\ has no proper $\lambda$-wim-constructible extension) if and only if it is $\lambda$-spherically complete.
\end{remark}
    
\begin{theorem}\label{thm:uniquebddwimm}
	For every $(\bE, \cO) \models T_\convex$ and every uncountable cardinal $\lambda$, maximal $\lambda$-bounded wim-constructible extensions of $(\bE, \cO)$ exist, are unique-up-to-non-unique isomorphism, are universal (i.e.\ weakly terminal) for $\lambda$-bounded wim-constructible extensions, and embed into every $\lambda$-spherically complete elementary extension of $(\bE, \cO)$.
	\begin{proof}
        First notice that if $(\bU, \cO')\succ (\bE, \cO)$ is a $\lambda$-saturated elementary extension of $(\bE,\cO)$ and $(x_i:i<\mu)$ is a maximal $\lambda$-bounded wim-construction within $(\bU, \cO')$, then it is maximal in every elementary extension of $(\bU, \cO')$: in fact if there were unary $\lambda$-bounded wim types over $(\bE\langle x_i: i<\mu\rangle, \cO' \cap \bE\langle x_i: i<\mu\rangle)$ they would be realized in $(\bU, \cO')$ by $\lambda$-saturation and $\big(\bE\langle x_i: i<\mu\rangle, \cO' \cap \bE\langle x_i: i<\mu\rangle\big)$ would not be maximal within $(\bU, \cO')$.

        It follows, for such a $\lambda$-bounded wim-construction $(x_i:i<\mu)$, that if $\bE_2:=\bE \langle x_i:i<\mu\rangle$ and $\cO_2:=\bE_2 \cap \cO'$, then $(\bE_2, \cO_2)$ is a maximal $\lambda$-bounded wim-constructible extension of $(\bE,\cO)$.
            
        Now let $(\bE_1, \cO_1)$ be another $\lambda$-bounded wim-constructible extension of $(\bE,\cO)$.
        By Lemma~\ref{lem:wim_amalgamate} there is an embedding of $j:\bE_1\to\bE_2$ over $\bE$ such that $\bE_2$ is $\lambda$-bounded wim-constructible over $j\bE_1$. If, furthermore, $(\bE_1, \cO_1)$ is a maximal $\lambda$-bounded wim-constructible extension, then $j$ is also surjective.
            
		Finally if $(\bE_3, \cO_3)$ is a $\lambda$-spherically complete extension of $(\bE,\cO)$, then $(\bE_2, \cO_2)$ embeds in $(\bE_3, \cO_3)$. To see this it suffices to show that $\tp((x_i: i< \mu)/(\bE, \cO))$ is realized in $(\bE_3, \cO_3)$.

        For this, in turn, it suffices to show that if $\alpha<\mu$, $(z_i:i< \alpha)$ is a realization of $\tp((x_i: i < \alpha)/(\bE, \cO))$ in $(\bE_3, \cO_3)$, and $\varphi : \bE \langle x_i: i<\alpha\rangle\to \bE \langle z_i: i<\alpha\rangle$ is the corresponding elementary embedding over $\bE$ given by $\varphi(x_i):=z_i$, then $\varphi_*\tp(x_\alpha/\bE \langle x_i:i<\alpha\rangle)$ is realized in $\bE_3$.
            
        Now, since $\varphi_*\tp(x_\alpha/\bE \langle x_i:i<\alpha\rangle)$ is $(\cO_3 \cap \bE \langle z_i:i<\alpha\rangle)$-wim and $\lambda$-bounded, and $(\bE_3, \cO_3)$ is $\lambda$-spherically complete, we have that $\varphi_*\tp(x_\alpha/\bE \langle x_i:i<\alpha\rangle)$ is realized in $(\bE_3, \cO_3)$ and the proof is complete.
	\end{proof}
\end{theorem}
	
In view of Remark~\ref{rmk:lambda_spherical_completeness}, Theorem~\ref{thm:uniquebddwimm} can be restated as
	
\begin{corollary}\label{thm:B}
	Let $(\bE, \cO)\models T_\convex$, and $\lambda$ be an uncountable cardinal. There is a unique-up-to-non-unique-isomorphism $\lambda$-spherically complete $\lambda$-bounded wim-constructible extension $(\bE_\lambda, \cO_\lambda)$ and it elementarily embeds in every $\lambda$-spherically complete extension of $(\bE, \cO)$.
\end{corollary}
	
\begin{definition}
	If $(\bE, \cO) \models T_\convex$, call a maximal $\lambda$-bounded weakly immediate extension the \emph{$T$-$\lambda$-spherical completion of $(\bE, \cO)$}.
\end{definition}
	
We conclude this subsection by observing that maximal $\lambda$-bounded wim-constructible extensions are homogeneous over wim-constructibly embedded $\lambda$-bounded wim-constructions.
	
\begin{corollary}
	Let $(\bE, \cO) \prec (\bE_\lambda, \cO_\lambda)$ be the $\lambda$-bounded spherical completion of $(\bE, \cO)$, and let $(x_i: i< \alpha)$, $(y_i: i< \alpha)$ be $\lambda$-bounded wim-constructions over $(\bE, \cO)$ with the same type over $(\bE, \cO)$, such that $(\bE_\lambda, \cO_\lambda)$ is wim-constructible over $\bE \langle x_i :i< \alpha \rangle$ and $\bE \langle y_i :i< \alpha \rangle$. Then there is an automorphism $\varphi$ of $(\bE_\lambda, \cO_\lambda)$ over $(\bE, \cO)$ such that $\varphi(x_i)=y_i$ for each $i<\alpha$.
	\begin{proof}
		Let $\psi$ be the isomorphism over $\bE$ between $\bE \langle x_i :i< \alpha \rangle$ and $\bE \langle y_i :i< \alpha \rangle$, given by $\psi(x_i)=y_i$. By Theorem~\ref{thm:uniquebddwimm}, $\psi$ extends to an automorphism $\varphi$ of $(\bE_\lambda, \cO_\lambda)$.
	\end{proof}
\end{corollary}

\subsection{On definable spherical completeness}\label{ssec:definable_spherical_completeness}In this brief subsection\footnote{The material of this subsection is originally from a separate note of the author on the arxiv \cite{freni2024t2}. The author has no intention of publishing such note on its own.} we show how Theorem~\ref{thm:types} entails definable spherical completeness of $T_\convex$. For this we will first need a Lemma.

\begin{lemma}\label{lem:cofinalities-in-vg}
    Let $(\bE, \cO)\models T_\convex$.
    Every $(\bE, \cO)$-definable subset $X \subseteq \val(\bE, \cO)$ either has a maximum or has the same cofinality as one of the three sets $\cO$, $\co$ and $\bE$.
    \begin{proof}
        Let $b$ be a realization of the special cut of $(\bE, \cO)$ so that $\cO$ is externally definable in the language of $T$ using $b$. By \cite[(3.10)]{dries1995t}, every $(\bE, \cO)$-definable subset $X$ of $\bE^{>0}$ can be written as a boolean combination of intervals and preimages of $\cO$ by monotone $\bE$-definable functions. In particular $X= \bE \cap Y$ for some $\bE\langle b\rangle$-definable set $Y\subseteq \bE \langle b\rangle$. Thus by o-minimality $X$ is a finite union of points and convex sets of the form $\bE \cap (a_0,a_1)$ where $a_0, a_1 \in \bE \langle b\rangle$. Now if $a_1 \in \bE$, then $(a_0, a_1) \cap \bE$ has the cofinality of $\bE$. Otherwise it has either the cofinality of $\bE^{<b}$ or of $\bE^{>b}$ and these are easily seen to be the cofinalities of $\cO$ and $\co$ respectively. Thus every $(\bE, \cO)$-definable subset of $\bE^{>0}$ either has a maximum or one of the three cofinalities in the statement. Finally observe that $\val_\cO|: \bE^{<0} \to \val(\bE, \cO)$ is weakly increasing, so if $X \subseteq \bE^{<0}$, then the cofinality of $\val_\cO(X)$ is either the cofinality of $X$ or $1$.
    \end{proof}
\end{lemma}

\begin{theorem}\label{thm:definable-spherical-completeness}
    $T_\convex$ is definably spherically complete.
    \begin{proof}
        Let $(\bE, \cO)\models T_\convex$ and $\lambda$ be a cardinal greater than the cofinalities of $\bE$, $\cO$ and $\co$. Let $(\bE_*, \cO_*)\succ (\bE, \cO)$ be the $T$-$\lambda$-spherical completion of $(\bE,\cO)$. It suffices to observe that $(\bE_*, \cO_*)$ is then definably spherically complete: by Remark~\ref{rmk:wimc-prserve} we have that $\cO$ and $\cO_*$ have the same cofinality as well as $\bE$ and $\bE_*$ and $\co$ and $\co_*$. Thus by Lemma~\ref{lem:cofinalities-in-vg}, in $(\bE_*, \cO_*)$ every nested $(\bE_*, \cO_*)$-definable family of valuation balls, when ordered under reverse inclusion, has cofinality strictly smaller than $\lambda$, and thus has non-empty intersection in $(\bE_*, \cO_*)$.
    \end{proof}
\end{theorem}

\begin{remark}\label{rmk:dependency}
    Notice the above proof didn't use the uniqueness of the $T$-$\lambda$-spherical completion and conceptually relies solely on Theorem~\ref{thm:types}. In fact, in the proof we could have chosen $(\bE_*, \cO_*)$ to be a maximal $\lambda$-bounded wim-constructible extension within some $\lambda$-saturated elementary extension of $(\bE, \cO)$.
\end{remark}

\subsection{Some Remarks and Questions}\label{ssec:quest}
It is natural to ask whether stronger results than Lemma~\ref{lem:wim_amalgamate} hold for specific theories $T$.
	
\begin{question}\label{quest:wimc-factors}
	Let $\bE \prec \bE_1$ be a wim-constructible extension.
	\begin{enumerate}
		\item[(1)] Are all $\bE_0$ with $\bE \prec \bE_0 \prec \bE_1$ wim-constructible over $\bE$?
		\item[(1w)] If $\bE \prec \bE_0 \prec \bE_1$ and $\bE_1$ is wim-constructible over $\bE_0$, is $\bE_0$ necessarily wim-constructible over $\bE$?
		\item[(2)] Is $\bE_1$ necessarily wim-constructible over all $\bE_0$ with $\bE \prec \bE_0 \prec \bE_1$?
		\item[(2w)] If $\bE \prec \bE_0 \prec \bE_1$ and $\bE_0$ is wim-constructible over $\bE$, is $\bE_1$ necessarily wim-constructible over $\bE_0$?
	\end{enumerate}
\end{question}
	
\begin{remark}
	For $i \in\{1,2\}$ a positive answer to $(i)$ implies a positive answer to $(iw)$.
\end{remark}
	
\begin{remark}
	Question~\ref{quest:wimc-factors}~(1) can be equivalently restated as: \textit{
	\begin{enumerate}
		\item[(1e)] are wim-constructible extensions $\lambda$-wim for all $\lambda$?
	\end{enumerate}}
\end{remark}
	
\begin{remark}
	The answer to both $(1)$ and $(2)$ in Question~\ref{quest:wimc-factors} is affirmative in the case $T$ is power-bounded (see Remark~\ref{rmk:pbdd_wimcfact}). We will give a very partial result toward (1) for some exponential theories in Section~\ref{sec:exp}.
\end{remark}
	
\begin{remark}\label{rmk:quest-relevance}
	Let $(\bE_\lambda, \cO_\lambda) \succeq (\bE, \cO)$ be the $\lambda$-spherical completion of $(\bE, \cO)$.
	Suppose Question~\ref{quest:wimc-factors}~(1) has an affirmative answer. Then every $\lambda$-spherically complete $(\bE_1, \cO_1)$ such that $(\bE_\lambda, \cO_\lambda) \succeq (\bE_1, \cO_1) \succeq (\bE, \cO)$ is isomorphic to $(\bE_\lambda, \cO_\lambda)$ over $(\bE, \cO)$.
	If furthermore Question~\ref{quest:wimc-factors}~(2) has an affirmative answer, then $\bE_1 = \bE_\lambda$.
		
	In particular an affirmative answer to Question~\ref{quest:wimc-factors}~(1), would imply that a $\lambda$-spherical completion for $(\bE, \cO)$ could be obtained as a union $\bigcup_{i< \lambda}\bE_i$ within a sufficiently saturated extension $(\bE_*, \cO_*)\succ (\bE, \cO)$ of a sequence of $\bE_i$ such that $\bE_i=\bigcup_{j<i} \bE_{j+1}$ and $\bE_{i+1}$ is the definable closure (for $T$) of a $\lambda$-spherical completion of $\bE_i$ \emph{qua real closed valued field} in $(\bE_*, \cO_*)$.
\end{remark}
	
\begin{remark}
	Adapting some results from \cite{mourgues1993every} and \cite{Ressayre1993Integer} (see also \cite{DAquino2012Real}), it can be shown that for every $T$ and every $(\bE, \cO) \models T_\convex$ there is big enough $\lambda$, such that the $\lambda$-bounded $\lambda$-spherical completion $(\bE_\lambda, \cO_\lambda)$ of $(\bE, \cO)$ is isomorphic as a real closed valued field (RCVF) to a $\lambda$-bounded Hahn field in the sense of \cite{kuhlmann2005kappa} with coefficients from the residue field of $(\bE, \cO)$. Moreover, the isomorphism can be chosen so that the image of $\bE \subseteq \bE_\lambda$ is truncation-closed and --- in the case $T$ defines an exponential --- also so that the logarithm of the monomial group is the set of purely infinite elements. This is done in more generality in \cite{freni2024structure} and won't be used here.
\end{remark}
	
The goal of the subsequent sections is to start the study of maximal $\lambda$-bounded wim-constructible models of $T_\convex$. For power-bounded theories this is essentially known, because if $T$ is power-bounded, all weakly immediately generated elements over some $(\bE, \cO)\models T_\convex$ are in fact weakly immediate. Therefore, if $x$ is weakly $\cO$-immediate over $\bE$, then $(\bE\langle x \rangle, \cO_x)$ is an \emph{immediate} extension of $(\bE, \cO)$, and thus a maximal $\lambda$-bounded wim-constructible extension of $(\bE, \cO)$ is an expansion of the $\lambda$-spherical completion of the reduct of $(\bE, \cO)$ to RCVF. This was first observed in \cite{kaplan2023t} and follows from results in \cite{dries1995t}, \cite{dries1997t}, and \cite{tyne2003t}, and we report more precisely on this in Section~\ref{sec:pbdd}.
	
For exponential o-minimal theories, the situation is not as well known and we commence its study in Section~\ref{sec:exp} by focusing on theories arising as expansions of power-bounded theories by an exponential function.

\section{The power-bounded case}\label{sec:pbdd}
	
In this section we review the known results about weakly immediate types in the case in which $T$ is power-bounded. The results are virtually all included in \cite{dries1995t}, \cite{dries1997t}, and \cite{tyne2003t} but we give somewhat different and shorter proofs.
		
\subsection{Introductory remarks}
	
Let $T$ be a complete o-minimal theory and $\bK \models T$. Recall that a \emph{power} of $\bK$ is a $\bK$-definable endomorphism $\theta: \bK^{>0} \to \bK^{>0}$ of the group $(\bK^{>0}, \cdot)$. By o-minimality such an endomorphism must be monotone and differentiable on its whole domain with derivative $\theta'(t)=\theta'(1)\theta(t)/t$. The \emph{exponent} of $\theta$ is defined as the element $\theta'(1)\in \bK$. It is easy to see that the set $\Exponents(\bK)$ of exponents of $\bK$ forms a subfield of $\bK$. The standard notation for powers is to write a power $t \mapsto \theta(t)$ as $t^{\alpha}:=\theta(t)$ where $\alpha=\theta'(1) \in \Exponents(\bK)$. Note that with this, the derivative of $t^\alpha$ has the familiar expression $\alpha t^{\alpha-1}$.
		
There is a fundamental dichotomy established by Miller (see \cite{miller1993growth}) regarding the subfield $\Exponents(\bK)$:
	
\begin{theorem}[Miller]
	The following are equivalent for an o-minimal expansion $\bK$ of a real closed field.
	\begin{enumerate}[label=(\arabic*), ref=(\arabic*)]
		\item $\bK$ does not define an exponential, that is, an isomorphism $\exp : (\bK,+) \to (\bK^{>0}, \cdot)$;
		\item \label{Miller:exponents_invariant} for every $\bK_1 \succeq \bK$, $\Exponents(\bK)=\Exponents(\bK_1)$; 
		\item every $\bK$-definable function $f: \bK \to \bK$ is eventually bounded by some power, that is there is $\beta \in \Exponents(\bK)$, such that for big enough $t \in \bK$, $f(t) < t^\beta$;
		\item for every $\bK$-definable function $f: \bK \to \bK$ there is an exponent $\beta \in \Exponents(\bK)$, such that $\lim_{s\to \infty} f(s)/s^\beta \in \bK$. 
	\end{enumerate}	
\end{theorem}
		
\begin{definition}[Miller]
    If $\bK$ satisfies one of the equivalent conditions above, then it is said to be \emph{power-bounded}. A complete o-minimal theory $T$ is said to be \emph{power-bounded} if one (or equivalently all) of its models is power-bounded. For a power-bounded theory $T$, $\Exponents(T)$ is defined to be the field of exponents in its prime model (this is then, by \ref{Miller:exponents_invariant} above, the field of exponents of every $\bK \models T$).
\end{definition}
		
The main thing we will need besides Miller's Theorem above is the following Lemma, also due to Miller.
		
\begin{lemma}[Miller]\label{lem:pbdd_asymp}
	If $T$ is power-bounded, $\bK \models T$ and $f: \bK \to \bK$ is $\bK$-definable and $\gamma$ is an exponent of $T$ then the following are equivalent
	\begin{enumerate}
		\item $\displaystyle{\lim_{t \to \infty}t f^{\dagger}(t)=\gamma}$;
		\item there is $k \in \bK^{\neq 0}$ such that $\displaystyle{\lim_{t \to \infty} \frac{f(t)}{kt^\gamma}}=1$.
		\end{enumerate}
	\begin{proof}
		Recall that for every $\bK$-definable $f$ there is an exponent $\gamma$ such that $(2)$ holds.
		If $\lim_{t \to \infty} f(t)= \pm \infty$, by power-boundedness, there is a power $t^\gamma$ such that $\lim_{t \to \infty} f(t)t^{-\gamma} \in \bK^{\neq 0}$ and by L'Hôpital's rule
		\[\lim_{t\to \infty}\frac{f(t)}{t^{\gamma}}=\lim_{t \to \infty} \frac{f'(t)}{\gamma t^{\gamma-1}}.\]
		So the limit in $(1)$ is $\gamma$.
		If $\lim_{t \to \infty} f(t) = 0$ the argument is the same.
		If $\lim_{t \to \infty} f(t) = c \in \bK^{\neq 0}$, then $f(t)= c+ g(t)$ where $\lim_{t \to \infty} g(t)=0$. Then it suffices to observe that \[tf^\dagger(t)=tg^\dagger(t) g(t)/f(t),\]
		so the limit in $(1)$ is in fact $0$.
	\end{proof}
\end{lemma}
	
\subsection{The rv-property}

As observed in \cite{dries1997t} the hypothesis of being power-bounded has strong consequences on $T_\convex$, in particular with regard to the structure induced on the value-group sort (cf \cite[Thm.~4.4]{dries1997t}). The starting point of this analysis is \cite[Prop.~4.2]{dries1997t} which is here convenient to report with some of its consequences in the following

\begin{lemma}\label{lem:pbdd_char}
	Let $T$ be a complete o-minimal theory. Then the following are equivalent
	\begin{enumerate}
		\item $T$ is power-bounded;
		\item for every $(\bE, \cO)\models T_\convex$ and every unary definable function $f: \bE \to \bE^{> 0}$
		\[(\bE, \cO) \models \exists a \in \bE,\; \exists M \in \cO, \; \forall t\in \cO,\; t > M \rightarrow f(t)\asymp a\]
		\item for every $(\bE, \cO)\models T_\convex$ and every unary definable function $f: \bE \to \bE^{> 0}$
		\[(\bE, \cO) \models \exists M \in \cO, \; \forall t\in \cO,\; t > M \rightarrow f^\dagger(t)\in \cO\]
		\item for every $(\bE, \cO) \models T_\convex$ and every unary definable function $f: \bE \to \bE^{> 0}$ there is a natural number $k$ such that
		\[(\bE, \cO) \models \exists t_0, \ldots t_{k-1},\; \forall t, \bigwedge_{i\in k} |t -t_i|> \co \rightarrow  f^\dagger(t) \in \cO.\]
		\item for every $(\bE, \cO) \models T_\convex$ and every unary definable function $f: \bE^{>0} \to \bE^{> 0}$ there is a natural number $k$ such that
		\[(\bE, \cO) \models \exists t_0, \ldots t_{k-1},\; \forall t,\, \left(t>0 \;\&\; \bigwedge_{i\in k}t \not\sim t_i \right) \rightarrow t f^\dagger(t) \in \cO.\]
	\end{enumerate}
	\begin{proof}
		$(1) \Rightarrow (2)$. This is essentially \cite[Proposition~4.2]{dries1997t}. We repeat the argument here for the convenience of the reader. Since $T_\convex$ is complete it suffices to prove it holds for every $\bE$-definable function $f:\bE\to \bE$ in some model $(\bE, \cO)\models T_\convex$. Take $\bK\models T$ of cofinality greater than $|T|$, let $d>\bK$ and consider $(\bK\langle d \rangle, \cO)\models T_\convex$ with $\cO=\CH(\bK)$. Since 
		\[\cO \ni t\mapsto \val_{\cO}f(t)\in \val_{\cO} \bE = \val_{\cO} d^{\Exponents(T)},\]
		is eventually a monotone function by o-minimality of $\bE$, and the cofinality $\cO$ is greater than the cardinality of $\val_{\cO}(\bE)$, such a function must be eventually constant.
			
		$(2) \Rightarrow (3)$. Suppose $(3)$ fails. Then for some $M \in \cO$, for every $t\in \cO^{>M}$, $|f^\dagger(t)|>\cO$. We can also assume, up to possibly choosing a larger $M\in \cO$ and working with one among $-f$, $1/f$ or $-1/f$ instead of $f$, that for every $t_1, t_2\in \cO$ with $t_2 > t_1 > M$, $f(t_2)>f(t_1)>0$, so in particular for $t>M$, $f^\dagger(t)>\cO$. Now for every $t_1, t_2 \in \cO^{>M}$ with $t_1<t_2+\co$, we can pick $t_0$ such that $t_1< t_0+\co<t_2$. Then for some $\xi \in (t_0, t_2)$
		\[\frac{f(t_2)}{f(t_1)} = \frac{f(t_0)+f(t_2)-f(t_0)}{f(t_1)} = \frac{f(t_0)}{f(t_1)} + f^\dagger(\xi) (t_2-t_0) \frac{f(\xi)}{f(t_1)} > (t_2-t_0)f^\dagger(\xi),\]
		by the mean value theorem and because $f$ is increasing and positive. Since, by hypothesis $f^\dagger(\xi)> \cO$ and $t_2-t_0>\co$, it follows that $f(t_2)/f(t_1)>\cO$. Since $t_1$ and $t_2$ were arbitrary with $t_2> t_1+\co$, this contradicts $(2)$.
			
		$(3) \Rightarrow (4)$. Assume $\lnot (4)$, then by weak o-minimality of $(\bE, \cO)$ there would be $t_0, t_1 \in \bE$, $t_1>t_0+\co$ such that for every $t \in (t_0, t_1)$, $f^\dagger(t)>\cO$. Consider $t_2:=(t_0+t_1)/2$ and define $g(t)=f(t_2-1/t)$. Then for $t>2/(t_1-t_0)$ and $t \in \cO\setminus \{0\}$,
		\[g^\dagger(t) = t^{-2} f^\dagger(t_2-1/t)>\cO\]
		contradicting $(3)$.
			
		$(4) \Rightarrow (5)$. Assume $(5)$ fails. Then, again by weak o-minimality, there are some $t_1 > t_0 > 0$, $t_1 \not \sim t_0$ such that for every $t \in (t_0, t_1)$, $|tf^\dagger(t)| > \cO$. Without loss of generality, we can assume $t_1\asymp t_0$, so considering $g(t):=f(t_0 t)$ we would have that for $t\in (1, t_1/t_0)$, $g^\dagger(t) = t_0 f^\dagger(t_0 t) \asymp t_0 t f^\dagger(t_0t) \succ \cO$ and this contradicts $(4)$ as $t_1/t_0>1+\co$.
			
		$(5) \Rightarrow (3)$ is clear.
			
		$(3) \Rightarrow (1)$ follows from Miller's dichotomy, because if $T$ is not power-bounded it defines an exponential $\exp: \bE \to \bE$, satisfying $\exp^\dagger(t)=1$ for all $t$, but then for $d> \cO$, the definability of $f(t):=\exp(dt)$ contradicts $(3)$ as $f^\dagger(t)=d$. 
	\end{proof}
\end{lemma}
    
\begin{corollary}\label{cor:pbdd_rescof}
	If $T$ is power-bounded, $(\bE, \cO) \models T_\convex$, $x$ is $\cO$-cofinally residual over $\bE$, and $\cO':=\cO_x=\cO_x^-$, then $\val(\bE \langle x\rangle, \cO') = \val(\bE, \cO)$.
	\begin{proof}
		By hypothesis there is $z \in \bE \langle x \rangle$ such that $\bK\langle z \rangle \preceq_\tame \bE\langle x \rangle$ and $\bK \langle z \rangle$ is a cofinal extension of $\bK$. Let $f$ be an $\bE$-definable function. By (5) in Lemma~\ref{lem:pbdd_char}, there are $t_1, t_2\in \bE$ such that $\co<t_1+\co<z<t_2+\co\subseteq \cO$, $f$ is increasing on $(t_1, t_2)$, and for every $t \in(t_1, t_2)$, $tf^\dagger(t) \in \cO$. By Lemma~\ref{lem:basicmvt}~(2), it then follows that for each $t \in (t_1, t_2)$, $f(t)\asymp f((t_2+t_1)/2)$, so $f(z)\asymp f((t_2+t_1)/2)$.
	\end{proof}
\end{corollary}

Another consequence of Lemma~\ref{lem:pbdd_char}, and of the relative quantifier elimination for $T_\convex$ is that, if $T$ is power-bounded, then the value-group sort in $T_\convex$ is o-minimal. We include a proof sketch for the convenience of the reader.

\begin{proposition}[Prop.~4.3 in \cite{dries1997t}]\label{prop:pbdd_vg_omin}
    If $T$ is power-bounded and $(\bE, \cO)\models T_\convex$ then the imaginary sort $\val(\bE,\cO)$ with its induced structure is o-minimal.
    \begin{proof}[Proof Sketch]
        Write a definable subset of $\val(\bE, \cO)$ as $\val(C)$ for some $C\subseteq \bE$ definable in the structure $(\bE, \cO)$. Since $(\bE, \cO)$ is weakly o-minimal, it suffices to show that $\val (C)$ has an infimum in $\val(\bE, \cO)\cup \{-\infty\}$. By \cite[(3.10)]{dries1995t} (here Theorem~\ref{thm:DriesLewMain}), $C$ is a boolean combination of intervals and sets of the form $C':=\{x \in \bE: f(x) \in \cO\}$ for some $\bE$-definable function $f$. Observe that, therefore, it suffices to show that each convex component of $\val(C')$ for such $C'$ has an infimum. For this use $(1)\Rightarrow (2)$ of Lemma~\ref{lem:pbdd_char} above.
    \end{proof}
\end{proposition}

\begin{lemma}\label{lem:vginducedaff}
    Suppose that $T$ is power-bounded and $\beta \in \Exponents(T)$.
	Let $f$ be $\bE$-definable, differentiable and never $0$. If $f^\dagger(t)t \sim \beta$ for all $t$ in some interval $J$, then on that interval $f(t)\sim c t^{\beta}$ for some $c$.
	\begin{proof}
		Let $h(t):=f(t)/t^\beta$ and $r<s \in J$. Then there is $r<\eta<s$ s.t.
		\[\frac{h(s)-h(r)}{h(\eta)}=h^{\dagger}(\eta)\eta \cdot \frac{s-r}{\eta} = \frac{s-r}{\eta} \left(\eta f^\dagger(\eta) - \beta \right)\]
		So, since $\eta f^\dagger(\eta) - \beta \prec \beta \asymp 1$, if $r \asymp s$, then $h(r)\sim h(s)$.
        Notice since $h$ is continuous and never $0$, we can without loss of generality assume that $h>0$.
        Moreover by o-minimality and by continuity of $h$, we can reduce to the case in which $h$ is continuous and monotone on $J$ and $0\notin J$. Finally, up to replacing $f$ with $1/f$, and precomposing with a sign-change we can also assume that $h$ is in fact continuous increasing and $0<J$.
        
        With these assumptions, the induced function $\tilde{h} : \val_\cO(J) \to \rv_\cO(\bE)$, given by $\tilde{h}(\val_\cO (t)) = \rv_\cO (h(t))$ is a continuous decreasing function. Let $\gamma_0 \in \val_{\cO}(J)$. Then $D:=\{\gamma \in \val_{\cO}(J): \tilde{h}(\gamma)>\tilde{h}(\gamma_0)\}$ is convex and open (for the order topology) in $\val_{\cO}(J)$. Moreover $D$ is definable in the imaginary sort $\val(\bE,\cO)$, which by \cite[Prop.~4.3]{dries1997t} (here Proposition~\ref{prop:pbdd_vg_omin}), is o-minimal. Thus the set $D$, if non-empty, is an open interval in the order $\val_{\cO}(J)$.
        
        Suppose toward contradiction that $D$ is non-empty. Note that since $\tilde{h}$ is continuous decreasing, $D$ is bounded above in $\val_\cO(J)$, e.g.\ by $\gamma_0\notin D$, thus we can find $r \in J$ such that $\val_\cO(r)$ is the supremum of $D$ and not an endpoint of $\val_\cO(J)$.
        Consider the map $g(t) \coloneqq h(rt)/h(r)$. Since $g$ is continuous increasing and by construction $g(\cO^{>\co})=1+\co$ there is an initial segment of $\bE^{>\cO}$ on which $g$ is continuous strictly increasing that gets mapped to an initial segment of $1+\bE^{>\co}$ contradicting Proposition~\ref{prop:coffact}. 
        
        Thus $D$ is empty and $\tilde{h}(\gamma) \le \tilde{h}(\gamma_0)$ for all $\gamma \in \val_\cO(J)$. An analogous argument shows that in fact $\tilde{h}(\gamma) \ge \tilde{h}(\gamma_0)$ for all $\gamma \in \val_\cO(J)$, whence $\tilde{h}$ must be constant and we have the thesis.
	\end{proof}
\end{lemma}

\begin{lemma}\label{lem:pbdd-spec}
    Let $T$ be power-bounded, $(\bU,\cO) \models T_\convex$, $\bE \prec \bU$, and $\cO<b<\bE^{>\cO}$. Then, for every $\bE$-definable $f: \bE \to \bE$ not constant around $b$, there is a unique $\beta_f\in \Exponents (T)\setminus \{0\}$ and some $c_0, c_1 \in \bE$ such that $f(b)-c_0 \sim_\cO c_1 b^{\beta_f}$.
	Moreover if $\{f_s :s \in D\}$ is an $\bE$-definable family of unary functions, $\{\beta_{f_s}: s \in D\}$ is finite.
	\begin{proof}
        Without loss of generality we can assume that $\bU$ is $|T|^+$-saturated.
        Let $\cO_b = \{x \in \bE\langle b \rangle: |x|<\bE^{>\cO}\}$ and observe that $\bE\langle b \rangle \cap \cO=\CH_{\bE\langle b\rangle}(\cO \cap \bE)=:\cO_b^-$.
        Recall that, by Theorem~\ref{thm:types}, if $\bK \prec \bE$ is maximal in $\{\bK'\prec \bE: \bK' \subseteq \cO\}$, then $\bK$ is also maximal in $\{\bK'\prec \bE \langle b \rangle: \bK' \subseteq \cO_b^-\}$.
            
		Since $f$ is not constant around $b$, $f(b) \notin \bE$ and by Theorem~\ref{thm:types_rest}, $f(b)$ is neither cofinally residual nor weakly immediate, hence for some $c_0$, $\val_{\cO}(f(b)-c_0) \notin \val_{\cO}(\bE)$. Since the theory is power-bounded, by (2) of Lemma~\ref{lem:pbdd_char}, there is $a\in \bE$ and $\sigma \in \{\pm 1\}$ such that $a(f(b)-c_0)^\sigma \in \cO_b^{>\cO}$.
		Now $g(t):=a (f(t)-c_0)^\sigma$ is such that $\bE^{>\cO}>g(b)> \cO$, so by Lemma~\ref{lem:transl}, there is a $\bK$-definable $h$ infinite at infinity and such that $h(b)-g(b)\in \co_b$. By power-boundedness, $h(b)\asymp_\cO b^{|\beta|}$ for some positive exponent $|\beta|$, hence $g(b)/b^{|\beta|} \in k+\co$ for some $k \in \bK$. Setting $c_1=(k/a)^{\sigma}$ and $\sign(\beta)=\sigma$, then yields $f(b) - c_0 \sim_\cO c_1 b^\beta$.
			
		Now we show that if $f(b)-c_0' \sim_\cO c_1'b^{\beta'}$ then $c_0 - c_0'\preceq_\cO 1$, $c_1 \sim_\cO c_1'$ and $\beta = \beta'$. In fact clearly we must have $\val_{\cO}(c_0 - c_0')> \val_{\cO}(c_1 b^\beta)$ and $c_1b^\beta \sim_\cO c_1'b^{\beta'}$, but then $b^{\beta'-\beta} \sim_\cO c_1/c_1'$ so $\rv_{\cO}(b^{\beta-\beta'}) \in \rv_{\cO}\bE$ and this is possible only if $\beta = \beta'$.
			
		For the last bit first note that up to enlarging $\bE$, we can assumed that $\bK\neq \Exponents(T)$, so in order to prove that $\{\beta_{f_s}: s \in D\}$ is finite it is enough to show that it is $\bK$-definable. This will follow from the claim that $b f_s'^\dagger(b)\sim_\cO \beta_{f_s}-1$ if $\beta_{f_s} \neq 1$ and $bf_s'^\dagger(b) \prec_\cO 1$ if and only if $\beta_{f_s} = 1$.

        To see why, recall $f(b) = c_0 + a^{-\sigma} \cdot (h(b)+\varepsilon(b))^\sigma$ where $h$ is $\bK$-definable and infinite at infinity, and $\varepsilon(b) \in \co_b$. So if $\sigma=1$, $f'^\dagger(b)=(h+\varepsilon)'^\dagger(b)$, if instead $\sigma=-1$, then $f'^\dagger(b)= (h+\varepsilon)'^\dagger - 2(h+\varepsilon)^\dagger$.
        
        Now note that Lemma~\ref{lem:specialder} also implies that if $g_1, g_2$ are unary $\bE$-definable functions with $g_1(b) \notin \co_b$ and $g_1(b)-g_2(b) \in \co_b$, then $g_1^\dagger(b)-g_2^\dagger(b)\in \co_b$, in fact since $g_1(b)/g_2(b) \in 1+\co_b$ we get $(g_1^\dagger-g_2^\dagger)(b)=(g_1/g_2)^\dagger(b) \in \co_b$.
        
        It thus follows that since $h(b), h'(b) \notin \co_b$, $(h+\varepsilon)^\dagger-h^\dagger \in \co_b$ and $((h+\varepsilon)'^\dagger-h'^\dagger)(b) \in \co_b$ and  so for both $\sigma \in \{\pm1\}$, we have $f'^\dagger(b)-(\tilde{h})'^\dagger(b)  \in \co_b$ where $\tilde{h}=h^\sigma$.
        
        Thus it suffices to show $b\tilde{h}'^\dagger(b)\sim_\cO \beta_f-1$ if $\beta_f\neq 1$ and $b\tilde{h}'^\dagger(b)\prec_\cO 1$ if $\beta_f=1$, but this follows immediately from the fact that $(\bK \langle b\rangle, \cO_b^- \cap \bK\langle b \rangle)$ with the derivation $r(b)\mapsto r'(b)$ for all $\bK$-definable 1-ary function $r$ is the Hardy field of the o-minimal structure $\bK$ (cf \cite[Sec.~3]{miller1993growth}).
        In fact then we get by L'Hôpital's rule for $\bK$-definable functions that, as long as $\beta_f \neq 1$, $\tilde{h}(b) \sim_\cO cb^{\beta_f}$ implies $\tilde{h}'(b) \sim_\cO \beta_f cb^{\beta_f-1}$ and $\tilde{h}''(b) \sim_\cO \beta_f (\beta_f-1)cb^{\beta_f-2}$. Similarly when $\beta_f=1$, $\tilde{h}(b)'\sim_\cO c$ and $\tilde{h}''(b) \prec_\cO 1$.

        We thus have shown $\{\beta_{f_s} : s \in D\}=\std(\{1+bf_s'^{\dagger}(b): s \in D\})$ where $\std=\std_\bK: \cO_b^- \to \bK$ is the standard part map.
        Now note that since $\cO\cap \bE<b<\bE^{>\cO}$, the set
        \[\begin{aligned}
        X=\{(s,\varepsilon, 1+t) \in D\times \bE^2: |t-bf_s'^\dagger(b)|<\varepsilon\}
        \end{aligned}\]
        is definable with parameters in $(\bE, \cO\cap \bE)$ 
        and for all $t \in \bE$ and $s \in D$ we have
        \[\big(\forall \varepsilon \in \bE^{>\co}, |t-bf_s'^\dagger(b)|<\varepsilon\big)\Leftrightarrow \std(bf_s'^\dagger(b))=\std(t).\]
        Thus the $\bK$-standard part of the $(\bE, \cO\cap \bE)$-definable set $\{t\in \bE: \exists s \in D, \forall \varepsilon \in \bE^{>\co} , (s, \varepsilon, t) \in X\}$ is exactly $\std(\{1+b f_s'^\dagger(b): s \in D\})=\{\beta_{f_s}: s \in D\}$ which is therefore definable in $\bK$.
	\end{proof}
\end{lemma}

\begin{lemma}\label{lem:pbdd_generic}
	Let $T$ be power-bounded and $(\bE, \cO) \models T_\convex$.  For every $\bE$-definable $f: \bE \to \bE$, and for every $a\neq 0$ there is a finite set $F_a\subseteq \bE$ such that
	\[\forall t \in \bE \setminus (F_a+a\cO),\; \exists c,\; c-t \notin a\cO,\; (t-c)f'^\dagger(t) \in \cO,\]
	or equivalently,
	\[\forall t \in \bE \setminus (F_a+a\cO),\; f'^\dagger(t)\prec 1/a.\]
	\begin{proof}
		Observe that it suffices to prove the statement for $a=1$, as the statement for general $a$ follows by applying the statement with $a=1$ to the function $f_a(t):=f(at)$.
				
		For every $x$, by (5) of Lemma~\ref{lem:pbdd_char} applied to $t \mapsto f'(x+t)$ there is a finite set $F$ such that for all $t \notin F(1+\co)$, $tf'^\dagger(x+t)\in \cO$. In particular there are $t_-<\cO<t_+$ such that for all $t \in (t_-, t_+)\setminus \cO$, $tf'^\dagger(x+t)\in \cO$, and setting $x_-=x+t_-$ and $x_+=x+t_+$, we get $x_-< x+\cO<x_+$ such that $(t-x) f'^\dagger(t)\in\cO$ for all $t \in (x_-, x_+) \setminus x+\cO$ and thus $f'^\dagger(t)\prec 1$ for all $t \in (x_-, x_+)\setminus x+\cO$.
				
		Consider the $(\bE, \cO)$-definable set $D=\{t \in \bE : f'^\dagger(t)\notin \cO\}$. By weak o-minimality of $(\bE, \cO)$, $D$ must be a finite union of order-convex sets, but by the observation above, any order convex subset $C\subseteq D$ must be contained in a single coset of $\cO$ hence the thesis.
	\end{proof}
\end{lemma}
	
\begin{corollary}\label{cor:pbdd_mod}
	Let $T$ be power-bounded, $(\bE, \cO) \models T_\convex$, and $f: \bE \to \bE$ be $\bE$-definable. Then there is a natural number $n\in \bN$ such that for every $\cO$-submodule $M$ of $\bE$, there is $F_M\subseteq \bE$ with $|F_M|<n$ such that for every $x \in \bE$, if $x \notin F_M+M$, then $f'^\dagger(x)M \subseteq \co$.
	\begin{proof}
        Note by Lemma~\ref{lem:pbdd_generic} and a routine compactness argument, given any $(\bE_1, \cO_1)\models T_\convex$ and an $\bE_1$-definable function $f: \bE_1 \to \bE_1$, there is a natural number $n$ such that for all $a \in \bE_1$, there is a finite set $F_a$ with fewer than $n$ elements such that $\forall t \in \bE_1 \setminus (F_a+a\cO), \, af'^\dagger(t)\prec 1$.

        Also note that every $\cO$-submodule $M$ of $\bE$ always has the form $\val_\cO^{-1}(U)$ where $U$ is an upward-closed subset of $\val_\cO(\bE) \cup \{\infty\}$.
        
		Now let $(\bE_1,\cO_1)$ be an elementary extension of $(\bE,\cO)$ containing an element $z$ such that $\val_{\cO_1}(\bE)^{<U}<\val_{\cO_1}(z) \le U$ so that $M=z\cO_1\cap \bE$ with $z \in \bE_1$. There is a finite $F \subseteq \bE_1$ with fewer than $n$ elements such that whenever $x \notin F + z\cO_1$, $zf'^\dagger(x) \cO_1 \subseteq \co_1$.
			
		Let $F_M\subseteq \bE$ be a set with fewer than $n$ elements such that $F_M+M=(F+z\cO_1) \cap \bE$. Such set can be obtained picking representatives of $(t+ z\cO_1)\cap \bE$ for every $t \in F$ such that $(t+z\cO_1) \cap \bE \neq \emptyset$.
			
		Now if $x \in \bE$ and $x \notin F_M+M$, then $x \notin F+z\cO_1$ so $Mf'^\dagger(x) \subseteq z \cO_1 f'^\dagger(x) \subseteq \co$.
	\end{proof}
\end{corollary}

We are now ready to give a short proof of the residue-valuation property. This was proven first by van den Dries and Speissegger in \cite[Proposition~9.2]{dries2000field} for polynomially bounded structures over $\bR$ and later in full for power-bounded structures in \cite[]{tyne2003t}. An alternative, yet different, proof of the residue-valuation property in the polynomially bounded case is in \cite{nowak2007proof}.

\begin{proposition}[Tyne]\label{prop:rv-prop}
	Let $T$ be power-bounded, $x \in \bU \succ \bE \models T$ and $\cO$ a non-trivial $T$-convex valuation subring of $\bE$.
	If $x$ is weakly $\cO$-immediate then every $y \in \bE \langle x \rangle$ is weakly $\cO$-immediate.
	\begin{proof}
		Assume $(x_i)_{i< \lambda}$ is an $\cO$-p.c.\ sequence in $\bE$ for $x$ and let $f: \bE \to \bE$ be an increasing $\bE$-definable function.
		Without loss of generality we can assume that $\Br(x/\bE)(\bE) \neq 0$, for otherwise $x$ is dense and $f(x)$ is still dense by continuity of $f$ at generic points.
		Let $p:=\tp(x/\bE)$, $M := \Br(p)(\bE)$, and $M_x := \Br(p)(\bE\langle x \rangle)$, and recall that $p(\bE\langle x \rangle)=x+M_x$.
		By Corollary~\ref{cor:pbdd_mod} there is a finite $F_M\subseteq \bE$ such that for every $t \notin F_M+M$, $f'^\dagger(t)M \subseteq \co$.
			
        Notice that by definition of breadth, for each $c\in \bE$, either $c+M<x$ or $c+M>x$. Also notice that $\{c+M: c+M<x\}\subseteq \bE/M$ does not have a maximum because if $c_m+M$ was such a maximum, then $\val_{\cO_x}(x-c_m)\notin \val_{\cO_x} \bE$, against the hypothesis that $x$ is weakly $\cO$-immediate. Similarly $\{c+M: c+M>x\}$ does not have a minimum.
        It follows that there are $a_0, a_1 \in \bE$ such that $a_0<x<a_1$, $(a_0, a_1) \cap (F_M + M)= \emptyset$. By o-minimality of $\bE$, we can furthermore choose such $a_0, a_1$ so that $f'^\dagger|_{[a_0, a_1]}$ is continuous strictly monotone or constant. Therefore $|f'^\dagger(x)| \le \max\{|f'^\dagger(t)|: t \in [a_0, a_1]\}=:d$.
			
		Now $M=\Br(p)(\bE) \subseteq \co/d$, that is $\Br(p)(\bE)(t) \vdash d|t| < r$ for all $r \in \bE^{> \co}$.
		Since by Theorem~\ref{thm:types}, $\bE \langle x\rangle$ does not contain any $y$ with $\co<y<\bE^{>\co}$, $M_x=\Br(p)(\bE\langle x \rangle) \subseteq \co_x/d$ as well and $f'^\dagger(x)M_x \subseteq \co_x$.
			
		It follows that there is $\tilde{c} \in \bE \langle x \rangle$ such that $|x-\tilde{c}|>M_x$ and $(x-\tilde{c})f'^\dagger(x) \in \cO_x$: since $f'^\dagger(x)M_x\subseteq \co$, $1/(f'^\dagger(x)) \notin M_x$, and $\tilde{c}:=x-1/(f'^\dagger(x))$ does the job. Finally, since $p(\bE\langle x \rangle)=x+M_x$, and $\tilde{c} \notin x+M_x$, there is $c \in \bE$ between $x$ and $\tilde{c}$. For such $c$, $|x-c|<|x-\tilde{c}|$, so $|x-c|f'^\dagger(x) \in \cO_x$ and again by Theorem~\ref{thm:types}, $|(x-c)f'^{\dagger}(x)|<s$ for some $s \in \cO$, whence $|(t-c) f'^\dagger(t)| \in \cO$ for all $t$ in an $\bE$-definable interval containing $x$ and thus $f(x)$ is weakly immediate by Lemma~\ref{lem:wimpres}.
	\end{proof}
\end{proposition}

As a Corollary we get the following which to the knowledge of the author was first stated in \cite[Cor.~1.12]{kaplan2023t}.
	
\begin{corollary}
	If $T$ is a power-bounded theory, and $(\bE, \cO) \models T_\convex$, then there is $(\bE_1, \cO_1)\succeq (\bE, \cO)$ such that $\rv_{\cO_1}(\bE_1) = \rv_{\cO}(\bE)$ and $(\bE_1, \cO_1)$ is spherically complete.
	\begin{proof}
		Use Theorem~\ref{thm:uniquebddwimm} with a cardinal $\lambda$ greater than the cofinality of any subset of $\val_{\cO}(\bE)$ to get a $\lambda$-spherically complete weakly immediate $(\bE_1, \cO_1)\succeq (\bE, \cO)$ and observe that by Proposition~\ref{prop:rv-prop} $\rv_{\cO_1}(\bE_1)=\rv_{\cO}(\bE)$.
	\end{proof}
\end{corollary}

It is also worth mentioning that Proposition~\ref{prop:rv-prop} implies that the trichotomy of unary types over some $(\bE, \cO)\models T$ of Definition~\ref{def:types1} and Theorem~\ref{thm:types_rest}, in the case when $T$ is power bounded, agrees with the standard trichotomy of unary types over a power-bounded o-minimal structure (as presented in the proof of \cite[Thm.~5.2]{daquino2015note}, cf \cite[Rmk.~3.4]{daquino2015note}).

In conjunction with the result \cite[Thm.~4.4]{dries1997t}, that for power-bounded $T$, the induced theory of $\val_{\cO}\bE$ is the (o-minimal) theory of an ordered vector space over $\Exponents (T)$, Proposition~\ref{prop:rv-prop} implies:
	
\begin{corollary}\label{cor:rvpptcor}
	If $T$ is power-bounded, $x\in \bU \succ \bE \models T$, and $\cO$ is a non-trivial $T$-convex valuation subring of $\bU$, then if for some $y \in \bE\langle x\rangle$, $\val_\cO(y) \notin \val_\cO \bE$, then there are $c_1,c_2 \in \bE$ and $\beta \in \Exponents(T)$, such that $y\sim c_2 |x-c_1|^\beta$.
	\begin{proof}
		If $x$ is $(\cO\cap \bE)$-special then the thesis follows from Lemma~\ref{lem:pbdd-spec}. So we may assume $x$ is not $(\cO\cap \bE)$-special.
			
		If for some $y \in \bE\langle x\rangle$, $\val_\cO(y) \notin \val_\cO \bE$, then $y$ is not weakly immediate over $(\bE, \cO)$, so by Proposition~\ref{prop:rv-prop} it is also not weakly immediately \emph{generated} and $x$ is not weakly immediate. It is also not cofinally residual by Corollary~\ref{cor:pbdd_rescof}. It follows that there is $c_1$ such that $\val_{\cO} (x-c_1) \notin \val_{\cO}(\bE)$. Now $y=f(x)$ for some $\bE$-definable $f$.
			
		Set $g(t) := f(c_1+t)$. Since $\val_{\cO}(x-c_1) \notin \val_{\cO}(\bE)$, by (5) in Lemma~\ref{lem:pbdd_char}, $tg'^\dagger(t) \in \cO$, for all $t$ in some $\bE$-definable neighborhood of $x-c_1$. In particular, by Lemma~\ref{lem:basicmvt}~(2), $g$ induces a function $G : J_0 \to J_1$, with $J_0, J_1$ $(\bE, \cO)$-definable intervals in $\val_\cO(\bE\langle x \rangle)$ containing $\val_{\cO}(x-c_1)$, such that $G(\val_{\cO}(t)) = \val_\cO(g(t))$.
        
        Moreover by \cite[Thm.~4.4]{dries1997t} $G$ is affine in a neighborhood of $\val_\cO(x-c_1)$, so $\val_{\cO}(y) = G(\val_\cO(x-c_1)) = \val_{\cO}\tilde{c}_2+\beta\val_\cO(x-c_1)$, for some $\tilde{c}_2 \in \bE^{\neq 0}$.
        Finally since $y$ is not residual and $y^{-1}\tilde{c}_2|x-c_1|^\beta \in \cO \setminus \co$, we can find $c_2\in\bE$ such that $y \sim c_2 |x-c_1|^\beta$.
	\end{proof}
\end{corollary}

As shown in \cite{dries2002minimal} in the case of polynomially bounded theories over $\bR$, this implies a (weak) preparation theorem, generalizing the one from \cite{parusinski1994lipschitz}. The argument here is the same as \cite[Lemma~2.2]{dries2002minimal}, with a little extra care needed due to the fact that the theory may not have an Archimedean prime model.
	
\begin{proposition}
	Let $T$ be a power-bounded o-minimal theory, $\bK \models T$, and $f: \bK^{n+1} \to \bK$ be a definable function. Then there is a covering of $\bK^{n+1}$ by finitely many definable cells $\{C_i: i \in M\}$ such that for every $i$, there are $\theta_i: \bK^n \to \bK$, $a_i: \bK^n \to \bK$, $\lambda_i \in \Exponents(T)$ and $\emptyset$-definable constants, $0<u_i^-< u_i^+$ such that for every $C_i$ either $f$ is identically $0$ on $C_i$ or 
	\[\bK \models \forall (\overline{c},t) \in C_i,\; \theta_i(\overline{c}) \neq t \; \& \; u_i^-<\left|\frac{f(\overline{c},t)}{a_i (\overline{c})|t-\theta_i(\overline{c})|^{\lambda_i}}\right|<u_i^+\]
	\begin{proof}
		By compactness it suffices to show that for every $\bE \models T$, and every $1$-ary $\bE$-definable function, there is a partition of $\bE$ into finitely many $\bE$-definable sets $C_i$ such that for every $i$ either $f(C_i)=\{0\}$ or there are $a_i, \theta_i \in \bE$ with $\theta_i \notin C_i$, $\lambda_i \in \Exponents(T)$ and $0$-definable $0<u_i^-< u_i^+$ such that
		\[\bE \models \forall t \in C_i,\; u_i^-<\left|\frac{f(t)}{a_i(t-\theta_i)^{\lambda_i}}\right|<u_i^+.\]
		This in turn is equivalent to showing that for every $\bU \succ \bE$ and $x \in \bU \setminus \bE$, such that $f(x) \neq 0$, there are $a, \theta \in \bE$, $\lambda \in \Exponents(T)$, and $u^-, u^+$ $0$-definable such that $u^-<|f(x)/(a(x-\theta)^\lambda)|<u^+$. Let $\cO$ be the minimal $T$-convex subring of $\bE \langle x \rangle$ (that is $\cO$ is the convex hull of the prime model of $T$). If $\val_{\cO}f(x) \in \val_{\cO}(\bE)$, then there is $a \in \bE$ such that $f(x)/a \asymp_{\cO} 1$. On the other hand if $\val_{\cO}f(x)\notin \val_{\cO} \bE$, then by Corollary~\ref{cor:rvpptcor}, there are $a, \theta \in \bE$ and $\lambda$ such that $\val_{\cO} a + \lambda \val_{\cO} (x-\theta) = \val_{\cO}(f(x))$.
        In both cases there are $a, \theta \in \bE$ and $\lambda\in \Exponents(T)$ such that $f(x)/(a(x-\theta)^\lambda) \in \cO\setminus \co$. Since $\cO$ is by definition the convex hull of the prime model, this means there are $\emptyset$-definable $u_-, u_+$, such that $u_-<|f(x)/(a(x-\theta)^\lambda)|<u_+$.
	\end{proof}
\end{proposition}

In the last subsection we will also need the following further Corollaries of the rv-property.

\begin{corollary}\label{cor:rvpptcor2}
    Let $T$ be power-bounded, $(\bE, \cO)\models T_\convex$ and $(y_i)_{i<n}$ be a tuple of elements of some elementary extension of $(\bE,\cO)$ such that $(\val(y_i):i<n)$ is $\Exponents(T)$-linearly independent over $\val(\bE, \cO)$. Then for every $x \in \bE \langle y_i: i<n\rangle\setminus \bE$, there are $c,d \in \bE$ and a tuple $(\beta_i)_{i<n}$ of exponents of $T$, at least one of which non-zero, such that $\val(x-c) = \val(d) + \sum_i \val(y_i)\beta_i$.
    \begin{proof}
        By induction on $n$. If $n=0$, there is nothing to prove, so suppose $n>0$ and $x \in \bE\langle y_i: i<n\rangle \setminus \bE\langle y_i:i<n-1\rangle$. Then by Corollary~\ref{cor:rvpptcor}, there are $x_0, x_1\in \bE\langle y_i: i<n-1\rangle$ and $\beta_{n-1}\in \Exponents(T)\setminus \{0\}$ such that $x-x_0 \asymp x_1 \cdot |y_{n-1}|^{\beta_{n-1}}$.
        
        By inductive hypothesis, there are $c_0,c_1, d_0, d_1 \in \bE$ and tuples of exponents $(\beta_{i,0})_{i<n-1}$, $(\beta_{i,1})_{i<n-1}$, such that for $j \in \{0,1\}$, $\val(x_j-c_j) = \val(d_j) + \sum_{i<n-1} \val(y_i)\beta_{i,j}$.
        Now $x - c_0 = (x-x_0) + (x_0-c_0)$ and
        \[
        \begin{aligned}
        \val(x-x_0)= \val(d_1) + \beta_{n-1}\val(y_{n-1}) + \sum_{i<n-1}\beta_{i,1}\val(y_{i}) \\
        \neq \val(d_0)+ \sum_{i<n-1}\val(y_i)\beta_{i,0}=\val(x_0-c),
        \end{aligned}\]
        where the inequality is due to the fact that $(\val(y_i))_{i<n}$ is linearly independent over $\val(\bE)$ and $\beta_{n-1}\neq 0$.
        It follows that $\val(x-c_0)$ is either $\val(x-x_0)$ or $\val(x_0-c)$ which are both of the required form, so the proof is completed.
    \end{proof}
\end{corollary}

\begin{corollary}\label{cor:pbdd_residual}
    Let $T$ be power-bounded, $(\bE, \cO\cap \bE)\prec (\bU, \cO)\models T_\convex$, and suppose $x \in \bU \setminus \bE$ is residual. Then for all $z \in \bE \langle x \rangle\setminus \bE$, there are $d,c \in \bE^{\neq0}$ such that $d(z-c)\in\cO\setminus (\bE+\co)$.
    \begin{proof}
        If $x$ is $(\bE \cap \cO)$-special, then this follows from Lemma~\ref{lem:pbdd-spec}.
        If $x$ is cofinally residual, then $\val_\cO(\bE\langle x \rangle)=\val(\bE)$, by Corollary~\ref{cor:pbdd_rescof} but $\rv_{\cO}(\bE\langle x \rangle)\neq \rv_{\cO}(\bE)$ by Proposition~\ref{prop:rv-prop}. It follows that if $z \in \bE \langle x \rangle$, then there is $c \in \bE^{\neq 0}$ such that $\rv(z-c) \notin \rv(\bE)$. Since however $\val(z-c) \in \val(\bE)$, there is $d\in \bE^{\neq 0}$ such that $d(z-c) \in \cO \setminus (\bE+\co)$.
    \end{proof}
\end{corollary}
    
\subsection{Immediate extensions}
An easy consequence of Proposition~\ref{prop:rv-prop} and Lemma~\ref{lem:wim_basic} is that an extension of models of $T_\convex$ is wim-constructible if and only if it is immediate.
	
\begin{corollary}\label{cor:pbddwimc_im}
	If $T$ is power-bounded then $(\bE, \cO) \preceq (\bE_1, \cO_1)\models T_\convex$ is wim-constructible if and only if $\val(\bE)=\val(\bE_1)$ and $\res(\bE)=\res(\bE_1)$.
	\begin{proof}
		If $\val(\bE, \cO)=\val(\bE_1, \cO_1)$ and $\res(\bE, \cO)=\res(\bE_1, \cO_1)$, then by Lemma~\ref{lem:wim_basic} all elements of $(\bE_1, \cO_1)$ must be weakly immediate over $(\bE, \cO)$. On the other hand since immediate extensions are closed under left factors, $(\bE\langle x \rangle, \cO_1 \cap \bE \langle x \rangle)$ is again immediate over $(\bE, \cO)$ for any $x \in \bE_1$. A transfinite inductive argument thus shows that $(\bE_1, \cO_1)$ must be constructible.
			
		Conversely Proposition~\ref{prop:rv-prop} implies that if $x$ is weakly immediate over $(\bE, \cO)$ then $(\bE\langle x \rangle, \cO_x)$ is immediate over $(\bE, \cO)$ and again a transfinite induction shows that constructible extensions are immediate.
	\end{proof}
\end{corollary}
	
\begin{remark}\label{rmk:pbdd_wimcfact}
	This implies that when $T$ is power-bounded both the Questions~\ref{quest:wimc-factors}~(1) and (2) have positive answers because immediate extensions are closed both under left and right factors.
	It also follows that for any $(\bE, \cO)$, for all big enough $\lambda$ the $\lambda$-spherical completions of $(\bE, \cO)$ are spherically complete, contrary to what happens when $T$ is exponential.
\end{remark}

\section{The simply exponential case}\label{sec:exp}

In this section we initiate the study of wim-constructible extensions in some exponential o-minimal structures.
	
\subsection{Basic facts on exponential o-minimal theories}\label{ssec:expintro}
This subsection is dedicated to some fundamental known facts on exponential o-minimal structures which can be found in \cite[Sec.~2]{miller1993growth} and in \cite{kuhlmann2000ordered}.
	
\begin{definition}
	A complete o-minimal theory $T$ expanding $\RCF$ is \emph{exponential} if for some model $\bK\models T$, there is a $\bK$-definable \emph{exponential function}, i.e.\ a \emph{non-constant} function $F: \bK \to \bK$ such that $F(x+y)=F(x)F(y)$ for every $x,y\in \bK$.
\end{definition}

\begin{remark}[Prop.~3.2 in \cite{miller1993growth}]
    Any exponential function $F: \bK \to \bK$ must be differentiable, strictly monotone, and such that $F(\bK)=\bK^{>0}$. In fact we must have $F(x)=F(x/2)^2\ge 0$ and $F(x)\neq 0$ (as otherwise $F(\bK)=\{0\}$). Moreover $F: \bK \to \bK^{>0}$ must be a bijection, because $F^{-1}(1)$ and $F(\bK)$ are definable subgroups of $(\bK,+)$ and $(\bK^{>0}, \cdot)$ respectively and thus by o-minimality and the fact that $F$ is non-constant, $F^{-1}(1)=\{0\}$ and $F(\bK)=\bK^{>0}$. Finally from the property $F(x+y)=F(x)F(y)$ it follows that $F$ is differentiable somewhere if and only if it is differentiable everywhere. Since by o-minimality $F$ must be differentiable on a final segment of $\bK$, it is differentiable everywhere, moreover $F'(x)=F(x)F'(0)$, so each such $F$ is either increasing (when $F'(0)>0$) or decreasing (when $F'(0)<0$).
\end{remark}

\begin{remark}\label{rmk:exp_zerodef}
    As observed in \cite[Sec.~2]{miller1993growth}, if $T$ is exponential, then there is a $\emptyset$-definable exponential function $\exp$ with $\exp'(0)=1$. In fact, if $\bK\models T$ and $F$ is a $\bK$-definable exponential function, then $F=f(-,c)$ for some $\emptyset$-definable function $f$ and some tuple $c$ from $\bK$. But then
    \[\varphi(z):=``\big(\forall x, y,\, f(x,z)f(y,z)=f(x+y,z)\big)\, \&\, \big(\exists t,\; f(t,z) \neq 1\big)"\]
    is consistent in $\bK$. Since $\dcl_T(\emptyset)\preceq \bK$, there is a tuple $c'$ in $\dcl_T(\emptyset)$ such that $T \vdash \varphi(c')$, whence $G:=f(-,c')$ is a $\emptyset$-definable exponential. The sought $\emptyset$-definable exponential function $\exp$ with $\exp'(0)=1$ is then given by $\exp(t):=G(t/G'(0))$.
\end{remark}
 
\begin{remark}\label{rmk:exp_rescaling}
	Given two exponential functions $F$ and $G$ defined in a model $\bK \models T$, there is a constant $k\in \bK$ such that $F(t)=G(kt)$: indeed $F^{-1} \circ G: \bK \to \bK$ would be a definable non-zero endomorphism of $(\bK,+)$ hence, by o-minimality, a homothety. Thus all $\bK$-definable exponential functions are of the form $F(t)=\exp(at)$ for $a=F'(0) \in \bK^{\neq0}$ and $\exp$ defined in Remark~\ref{rmk:exp_zerodef} above is the unique exponential with $\exp'(0)=1$. We will refer to $\exp$ as \emph{the exponential} of $T$. We will write $\log: \bK^{>0} \to \bK$ for the compositional inverse of $\exp$.
\end{remark}
	
\begin{remark}
    As observed in \cite[Note~2.6]{miller1993growth}, if $T$ is an exponential theory then its exponential $\exp$ satisfies the following \emph{growth axiom}
    \[\forall c\ge 1,\, \forall t>c^2, \exp(t)>\exp(c\log(t)).\]
    This in particular entails that if $(\bE,\cO)\models T_\convex$, then for all $t>\cO$, we have $\exp t > t \cO$.
\end{remark}

\begin{remark}
	If $T$ is an exponential theory and $(\bE, \cO)\models T_\convex$ then $\exp(\cO)=\cO^{>\co}$ (by $T$-convexity and because $\exp: \bE \to \bE^{>0}$ is $\emptyset$-definable continuous). It follows that for every $a \in \bE$, $\exp(\cO+a)= \val_{\cO}^{-1}(\val_\cO(\exp(a))) \cap \bE^{>0}$.
\end{remark}

\begin{remark}[Immediate elementary extensions]\label{rmk:Imm_in_exp_theories}
    If $T$ is exponential and $(\bE, \cO)\models T_\convex$ then the immediate elementary extensions of $(\bE, \cO)$ must be dense. In fact if $(\bE\langle x\rangle, \cO')\succ (\bE, \cO)$ and $x$ is weakly immediate but not dense over $(\bE, \cO)$, then for some $a \in \bE$, we have $\cO'\subseteq \Br(ax/\bE)(\bE\langle x\rangle)$. But this implies that $\val_{\cO'} \exp(ax)\notin \val_{\cO'}(\bE)$, contradicting immediateness.
\end{remark}

The first expansions of $\bR$ defining an exponential that have been proven to be o-minimal were $\bR_{\exp}$, the ordered field of the reals expanded with the natural exponential \cite{wilkie1996model}, and its expansion by all restricted analytic functions $\bR_{an, \exp}$ \cite{dries1994real}. In fact $\bR_{an, \exp}$ also eliminates quantifiers in the language of restricted analytic functions together with $\exp$ and $\log$ \cite{dries1994elementary}. 
	
In \cite[Thm~3.2]{dries2002minimal} it was shown more generally that if $T_0$ is any theory of a power-bounded expansion $\bR_0$ of the field of reals and $T_0$ defines the restricted natural exponential $\exp|_{[0,1]}$ and eliminates quantifiers, then the expansion $\bR_{0, \exp}$ of $\bR_{0}$ by the unrestricted exponential function is o-minimal and its theory $T$ eliminates quantifiers in the language $L=L_0 \cup \{\exp, \log\}$ where $L_0$ is the language of $T_0$. Moreover $\bR_{0, \exp}$ is universally axiomatized provided that $T_0$ was. We will be essentially concerned with a slight generalization of this situation.
	
\begin{definition}\label{defn:simply_exp}
	We will call an exponential o-minimal theory $T$ \emph{simply exponential} if it has a power-bounded reduct $T_0:=T|_{L_0}$ to a language $L_0$ such that:
    \begin{enumerate}
        \item for all $T$-definable closed bounded intervals $I$, $\exp|_I$ is $T_0$-definable;
        \item $T$ is an expansion by definitions of $T':=T|_{(L_0\cup \{\exp\})}$.
    \end{enumerate}
\end{definition}

\begin{remark}
    Whether (1) is actually needed in the definition above seems to be a non-trivial matter.
    For once, when $T$ does not admit Archimedean prime models, it seems unknown whether $T$ needs to have a power-bounded reduct $T_0$ defining all restrictions of $\exp$ to closed bounded intervals (see \cite[p.~3]{berarducci2026elementary} which indicates that proving power-boundedness for the reduct of all such $T$s to $\{+, \cdot, \exp|_{[-1,1]}\}$ is an open problem). I thank the anonymous referee for pointing this issue out to me.
    
    In the case $T$ has an Archimedean prime model, whether the existence of a power bounded reduct $T_0$ satisfying (2) guarantees that there exists one satisfying (1) and (2) seems also non-trivial. In fact it seems plausible that some exponential $T$ might have power bounded reducts $T_0$ satisfying (2) which when expanded by $\exp|_I$ allow to define the unrestricted exponentiation, so given a $T_0$ satisfying (2), finding a new $T_0'$ satisfying (1) and (2) may require more than just expanding the given $T_0$.
\end{remark}

\subsection{A brief review of known results}\label{ssec:known_results}
The main goal of this subsection is to point out that if $T'$ is a simply exponential theory in some language $L'$, then there is an expansion by definition $T^*$ of $T'$, with a power-bounded reduct $T_0=T^*|_{L_0}$ such that $T^*$ is an expansion by definition of $T''=T^*|_{L_0 \cup\{\exp, \log\}}$ \emph{and moreover $T''$ eliminates quantifiers and is universally axiomatized} (Corollary~\ref{cor:equivalent_simply_exp}).

To do this we will need to review a generalization of \cite[Thm.~B]{dries2000field} and \cite[Thm~3.2]{dries2002minimal} to power-bounded o-minimal theories: more precisely what will be described is a relative quantifier elimination  (Proposition~\ref{prop:qe_for_log_exp}) for the expansion of a power-bounded structure $\bK$ by an exponential $\exp$ which is \emph{compatible} in an appropriate sense (Definition~\ref{def:T(exp,log)}).
The whole argument is essentially the one used in \cite{dries1994elementary} for $\bR_{an, \exp,\log}$, see also \cite{lion1997theoreme} (in French) which gives a preparation theorem for such structure. In particular the ideas and arguments in this subsection are virtually all contained in the previous literature and are not original to the author (cf \cite{Ressayre1993Integer}, \cite{dries1994real},  \cite{dries1994elementary}, \cite{kuhlmann2000ordered}, and \cite{dries2002minimal}). Furthermore, as the author was recently made aware of, a similar generalization was presented in Foster's thesis \cite[Ch.~6]{foster2010power}, which however is stated in a slightly more restrictive setting than the one needed here (see Remark~\ref{rmk:foster_comparison}).

Since some intermediate results (Lemmas~\ref{lem:log_scale_properties} and \ref{lem:exp_step}) will be used in the last subsection and since the content presented so far in this paper sets the proof of Proposition~\ref{prop:qe_for_log_exp} just a few lemmas away, the author felt that presenting its proof in full is a sensible choice in terms of readability.

\begin{context}\label{defn:T_with_restricted_exp}
    In this subsection we fix a power-bounded o-minimal theory $T$ in some language $L$ such that $T$ is universally axiomatized and eliminates quantifiers.
    We denote by $\Lambda:=\Exponents(T)$ the field of exponents of $T$.
\end{context}

\begin{definition}\label{def:T(exp,log)}
    Given a model $\bE$ of $T$ we will call an isomorphism $\exp: (\bE, +, <) \to (\bE^{>0},\cdot, <)$ a \emph{compatible} exponential if
    \begin{enumerate}
        \item for every \emph{bounded} $\emptyset$-definable interval $I$, $\exp|_I$ is $\emptyset$-definable and moreover
        $(\exp|_{[-1,1]})'(0)=1$;
        \item $\exp(t)^\beta =\exp (t\beta)$ for all $t \in \bE$ and $\beta \in \Lambda$; 
        \item it satisfies the growth axiom, that is for all $c,t \in \bE$ with $t>c^2>1$, we have $\exp(t)>\exp(c\log(t))$.
    \end{enumerate}
    We will denote by $\log$ the compositional inverse of such an $\exp$.
    We will denote by $T(\exp, \log)$ the common theory of all expansions $(\bE, \exp, \log)$ of some model $\bE$ of $T$ by a compatible exponential $\exp$ and its inverse $\log$.
\end{definition}

\begin{remark}\label{rmk:unique_bdd_exp}
    Note that if $\exp_1$ and $\exp_2$ are compatible exponentials on some $\bE\models T$, then $\exp_1|_I = \exp_2|_I$ for all $\emptyset$-definable bounded interval $I$. 
    To see why, note first that if $\exp$ is a compatible exponential, then $\exp|_I$ must satisfy $(\exp|_I)'=\exp|_I$: this is because given $x \in I$, and $\varepsilon$ such that $(x-\varepsilon, x+\varepsilon)\subseteq I$, $\exp|_{(x-\varepsilon, x+\varepsilon)}(x+t)=\exp(x) \exp(t)$ for all $t \in (-\varepsilon, \varepsilon)$, thus since $(\exp|_{(-\varepsilon, \varepsilon)})'(0)=1$, $\exp|_I$ is differentiable at $x$ and $(\exp|_I)'(x)=(\exp|_I)(x)$.
    
    It follows that each $\exp_i|_I$ should satisfy $\exp_i'|_I= \exp_i|_I$ for $i\in \{1,2\}$, thus $h_I:=(\exp_1|_I)/(\exp_2|_I)$ should satisfy $h_I'=0$ and thus all $h_I$ are constants $h_I=:c_I\in\bE$. Since given any $\emptyset$-definable interval $I$, there is a bounded $\emptyset$-definable interval containing $I \cup [-1,1]$, it follows that $c_I=1$ for all $\emptyset$-definable intervals $I$.
\end{remark}

\begin{remark}\label{rmk:T(exp,log)_extension_axiom_scheme}
    If some model $\bE\models T$ admits a compatible exponential $\exp$, then Remark~\ref{rmk:unique_bdd_exp} entails that there is a (non-definable) family of $\emptyset$-definable unary functions $e_I$ defined on bounded intervals, such that condition (1) of Definition~\ref{def:T(exp,log)} is equivalent to requiring $\exp|_{I}=e_I$ for all bounded $\emptyset$-definable $I$.
    In particular if an $\bE\models T$ admitting a compatible exponential exists, then the prime model $\bP$ of $T$ also admits a compatible exponential $\exp$, and $(\bP, \exp, \log)$ is a prime model for $T(\exp, \log)$.
\end{remark}

\begin{remark}\label{rmk:T(exp,log)univ_axiomatized}
     Note that by Remark~\ref{rmk:T(exp,log)_extension_axiom_scheme}, if $T(\exp, \log)$ is consistent then it is universally axiomatized in the language $L\cup \{\exp, \log\}$ and has an initial model $(\bP, \exp, \log)$.
\end{remark}

\begin{remark}[Comparison with the setting of \cite{foster2010power}]\label{rmk:foster_comparison}
    As already mentioned \cite[Ch.~6]{foster2010power} deals with a similar setting: there it is assumed that $\Exponents(T)$ is cofinal in the prime model $\bP$ of $T$ and that there is a $T$-definable restricted exponential on $[0,1]$. This entails that $\bP$ admits a compatible exponential (see \cite[6.4.1 and 6.4.3]{foster2010power}), but is not equivalent to the latter condition: for example one could have $T$ be artificially expanded by constants greater than any exponent of $T$. It should be noted however that the proof of Proposition~\ref{prop:qe_for_log_exp} presented below is essentially the same as the proof of \cite[Thm.~6.5.2]{foster2010power} and, to some extent, of \cite[Thm.~4.1, Cor.~4.5]{dries1994elementary}.
\end{remark}

In the following Lemmas we fix a model $(\bE, \exp, \log)$ of $T(\exp, \log)$ and a (possibly trivial) $T$-convex valuation subring $\cO\subseteq \bE$ such that $\exp(\cO)\subseteq \cO$. Notice that the growth axiom implies that then $\exp(\cO)=\cO^{>\co}$ and that for all $x\in \bE^{>\cO}$, $\log|x|\prec x$.

\begin{lemma}\label{lem:log(wim)}
    Let $\bE_0\subseteq \bE$ be an $L(\log)$-substructure and $\bE_1\subseteq \bE$ be a $L$-substructure containing $\bE_0$. If $\bE_1$ is immediate over $\bE_0$, then $\bE_1$ is $L(\log)$-closed. Furthermore if $\bE_* \supseteq \bE_0$ is a model of $T(\exp, \log)$ extending $(\bE_0,\log)$ and $\iota :\bE_1 \to \bE_*$ is an $L$-embedding over $\bE_0$, then $\iota$ is an $L(\log)$-embedding.
    \begin{proof}
        Let $z \in (\bE_1)^{>0}$. Notice that since $\bE_1$ is immediate over $\bE_0$, $z=z_0(1+z_1)$ for some $z_0 \in \bE_0$ and some $z_1 \in \bE_1\cap \co$, whence $\log(z)=\log(z_0) + \log(1+z_1)$. But then $\bE_0 \langle \log (z)\rangle_L = \bE_0\langle z_1\rangle_{L}\subseteq \bE_1$ because $\log|_{[\exp(-1),\exp(1)]}$ is in $L$, thus $\bE_1$ is $\log$-closed. Furthermore
        \[\iota(\log z)=\log(z_0) + \iota\log(1+z_1)= \log(z_0)+ \log(1+\iota(z_1))= \log(\iota(z)).\]
        Thus $\iota$ is $\log$-preserving.
    \end{proof}
\end{lemma}

\begin{lemma}\label{lem:log(res)}
    Let $\bE_0\subseteq \bE$ be an $L(\log)$-substructure and $x \in \bE\setminus \bE_0$. Suppose furthermore that $\cO=\CH(\bP)$. If $\res_\cO(\bE_0\langle x \rangle)\neq \res_\cO(\bE_0)$, then $\bE_0\langle x \rangle$ is $L(\log)$-closed. Moreover if $\bE_* \supseteq \bE_0$ is a model of $T(\exp, \log)$ extending $(\bE_0,\log)$ and $\iota:\bE_0\langle x \rangle \to \bE_*$ is an $L$-embedding over $\bE_0$, then $\iota$ is an $L(\log)$-embedding.
    \begin{proof}
        Suppose $z \in\bE_0 \langle x \rangle$. Then by Corollary~\ref{cor:pbdd_residual} for some $c,d\in \bE_0^{\neq 0}$ such that $z_1:=(z-c)/d \in \cO\setminus (\bE_0+\co)$. Notice that since $\cO=\CH(\bP)$, $\log$ is $T$-definable around every point of $\cO^{>\co}$. It follows that
        \begin{itemize}
            \item if $c \preceq d$, then $z_1+c/d \in \cO\setminus \co$ so $\log|z_1+c/d|\in \bE_0\langle x\rangle$ and $\iota(\log|z_1+c/d|)=\log|\iota(z_1+c/d)|$;
            \item if $c \succeq d$, then $1 + dz_1/c \in \cO \setminus \co$ so $\log|1+z_1d/c|\in\bE_0 \langle x \rangle$ and $\iota(\log|1+z_1d/c|)=\log|\iota(1+z_1d/c)|$. 
        \end{itemize}
        Since $\log|d|,\log|c| \in \bE_0$ and
        \[\log|d|+ \log|z_1+c/d|=\log|z|=\log|c|+\log|1+z_1d/c|,\]
        the thesis follows.
    \end{proof}
\end{lemma}

\begin{lemma}\label{lem:logE<y>}
    Let $\bE_0\subseteq \bE$ be an $L(\log)$-substructure.
    Suppose that $(y_i)_{i<n}$ is a finite sequence in $\bE$ such that $(\val(y_i): i<n)$ is $\Lambda$-linearly independent over $\val(\bE_0)$ and $\bE_1 := \bE_0 \langle y_i: i< n\rangle_L$.
    Then for all $z \in \bE_1$, $\log |z| \in \bE_1 \langle \log |y_i|: i<n\rangle_L$. 
    \begin{proof}
        Notice that without loss of generality we can assume $y_i>0$ for all $i<n$. By Corollary~\ref{cor:rvpptcor2}, for every $z \in \bE_1$, there are $c \in \bE_0$, $\beta_0, \ldots, \beta_{n-1} \in \Lambda$, and $z_1 \in \bE_1 \cap \co$, such that
        $z = c y_0^{\beta_0} \cdots y_{n-1}^{\beta_{n-1}} (1+z_1)$. It follows that $\log |z| = \log |c| + \sum_{i<n} \beta_i \log|y_i| + \log (1+z_1)$.
    \end{proof}
\end{lemma}

\begin{lemma}\label{lem:log_scale_properties}
    Let $\bE' \subseteq \bE$ be an $L(\log, \exp)$-substructure and suppose that $n \in\omega$ and $(y_i)_{i<n} \in \bE^{n}$ is such that $\val(y_{i}) \notin \val (\bE')$ for all $i<n$ and $c_i:=\log |y_i| - y_{i+1} \in \bE'$ for all $i$ with $i<n-1$.
    Then for all $i<n$
    \begin{enumerate}
        \item for all $j<i$, $1\prec \log|y_i| \prec |y_i| \preceq \log|y_j| \prec y_j$;
        \item $\val(y_{i}) \notin \val (\bE') + \sum_{j<i} \val(y_j)\cdot \Lambda$.
    \end{enumerate}
    In particular
    \begin{enumerate}[resume]
        \item $\tp((y_i)_{i<n}/\bE')$ only depends on $\tp(y_0/\bE')$ and $(c_i)_{i<n-1}$.
    \end{enumerate}
    \begin{proof}
        Suppose $j<i<n$.
        Since $\val(y_{j+1})\notin \val (\bE')$, we have $y_{j+1} \preceq \log |y_j|$ and $c_j \preceq \log |y_j|$.
        Moreover since $\val(y_j)\notin \val(\bE')$, we have in particular $y_j\exp(-c_j) \in \bE^{\not\asymp 1}$, whence $y_{j+1}=\log|y_j|-c_j\succ 1$.
        It follows that $1 \prec \log |y_{j+1}| \prec y_{j+1} \preceq \log |y_j|$, which entails (1).

        As for (2), note that it suffices to show that for all $a\in \bE$ with $\val(a)\in \val(\bE')$
        \begin{equation}\label{eq:log_exp_main_lemma_1}
            |y_i| \ge |a| \prod_{j<i} |y_j|^{\beta_j} \Longleftrightarrow |y_i| \succ |a| \prod_{j<i} |y_j|^{\beta_j}.
        \end{equation}
        Set $a=a'u$ with $a' \in \bE'$ and $u \in \cO^{>\co}$, 
        \[ B:=\log|y_i|-\sum_{j<i} \beta_j y_{j+1}, \quad\text{and}\quad A\coloneqq\log |a'| + \sum_{j<i} \beta_j c_j \in \bE',
        \]
        we get that the left-hand side of (\ref{eq:log_exp_main_lemma_1}) is equivalent to $B+\log(u) \ge A$. But if any of the $\beta_j$ is non-zero, then by (1) $\val(B)<0$ and $\val(B) \notin \val(\bE')$, thus $(B+\cO) \cap \bE' =\emptyset$. Therefore $B +\log(u) \ge A$ if and only if $B\ge A$ if and only if $B+\cO \ge A$, which in turn is equivalent to the right-hand side of (\ref{eq:log_exp_main_lemma_1}).

        As for (3), first note that since $\bE'$ is $\log$-closed, $\log$ is increasing, and $(c_j)_{j<i} \in (\bE')^i$, we have that $\tp(y_i/\bE')$ is determined by $\tp(y_0/\bE')$ and $(c_j)_{j<i}$. Then note that by Corollary~\ref{cor:rvpptcor2} $t \in \bE' \langle y_j: j<i\rangle^{>0}$ can be written as $t\coloneqq a \prod_{j<i} |y_j|^{\beta_j}$ for some $a>0$ with $\val(a) \in \val(\bE)$. So, using the same notations as above, the validity of $|y_i| > t \coloneqq |a| \prod_{j<i} |y_j|^{\beta_j}$ is determined by $\tp(y_i/\bE')$ if $\beta_j=0$ for all $j<i$, and by the validity of $B - \log |y_i| \ge A$ otherwise, because in this latter case $B\succ \log|y_i|$. The validity of $B - \log |y_i| \ge A$, in turn, only depends on $\tp(y_j/\bE')$ for the least $j<i$ for which $\beta_j\neq 0$. Thus in both cases $\tp(|y_i|/ \bE'\langle y_j : j<i\rangle)$ is determined by $\tp(y_0/\bE')$ and $(c_j)_{j<i}$. Since $\tp(y_i/\bE')$ is also determined by $\tp(y_0/\bE')$ and $(c_j)_{j<i}$ and it determines the sign of $y_i$, we deduce (3).
    \end{proof}
\end{lemma}

\begin{remark}
    Sequences $(y_i)_{i<n}$ as in Lemma~\ref{lem:log_scale_properties} are, in one form or another a key tool in many proofs of quantifier elimination of $T(\exp, \log)$. They correspond to \emph{logarithmic scales} in the language of \cite{dries2002minimal}.
\end{remark}

\begin{lemma}\label{lem:maximal_log_scale}
    Let $\bE' \subseteq \bE$ be an $L(\exp, \log)$ substructure and let $x \in \bE\setminus \bE'$. Then one of the following holds
    \begin{enumerate}
        \item there is a sequence $(y_i)_{i\le n}$ such that $x\in y_0+\bE'$ and for all $i<n$, $\val(y_{i})\notin \val(\bE')$ and $\log|y_{i}| - y_{i+1} \in \bE'$ and $y_n$ is either wim or residual over $\bE'$; 
        \item there is a sequence $(y_i)_{i<\omega}$ such that for all $i<\omega$, $\val(y_{i})\notin \val(\bE')$ and $\log|y_{i}|-y_{i+1} \in \bE'$. Moreover $\bE' \langle y_i: i< \omega\rangle$ is $\log$-closed.
    \end{enumerate}
    \begin{proof}
        Notice that a sequence $(y_i)_{i\le n}$ is maximal among the sequences such that $y_0\in x+\bE'$ and for all $i<n$ $\val(y_{i})\notin \val(\bE')$ and $\log|y_{i}|-y_{i+1} \in \bE'$, if and only if $y_n$ is wim or residual.
        Thus if (1) does not hold, there is $(y_i)_{i<\omega}$ as in (2). The fact that $\bE'\langle y_i: i< \omega\rangle$ is $\log$-closed follows again from Lemma~\ref{lem:logE<y>} because the set $\{\val(y_i): i<n\}$ is linearly independent and again $\log|y_i| \in \bE \langle y_i: i< \omega\rangle$.
    \end{proof}
\end{lemma}

\begin{lemma}\label{lem:exp_step}
    Suppose $\bE_0\subseteq \bE$ is an $L(\log)$-substructure and let $\bE_1:= \bE_0 \langle \exp(\bE_0)\rangle$. Then:
    \begin{enumerate}
        \item $\log(\bE_0^{>0})$ is a $\Lambda$-linear subspace of $\bE_0$;
        \item $\bE_1$ is $\log$-closed and $\bE_1 = \bE_0 \langle \exp(B_0)\rangle$ for any $\Lambda$-linear basis $B_0$ of $\bE_0$ over $\log(\bE_0^{>0})$;
        \item If furthermore $\bE_*\supseteq \bE_0$ is a model of $T(\exp, \log)$ extending $(\bE_0, \log)$, and $\iota: \bE_1 \to \bE_*$ is an $L$-embedding such that $\iota(\exp(y))=\exp(\iota(y))$ for all $y \in \bE_0$, then $\iota$ is an $L(\log)$-embedding. 
    \end{enumerate}
    \begin{proof}
        (1) $\log(\bE_0^{>0})$ is a $\Lambda$-linear subspace because by compatibility of $\exp$, $\log$ is a $\Lambda$-linear isomorphism from the multiplicatively written ordered $\Lambda$-vector space $(\bE^{>0}, \cdot, <, \{(-)^\beta\}_{\beta \in \Lambda})$ to $(\bE, +, <, \{(-\cdot \beta)\}_{\beta \in \Lambda})$. 
        (2) $\bE_1=\bE_0 \langle \exp(B_0)\rangle$ because the multiplicative $\Lambda$-linear span of $\exp(B_0) \cup \bE_0^{>0} \subseteq \bE^{>0}$ is $\exp(\bE_0)$ by construction and $\Lambda$-linearity of $\exp$. To see that $\bE_1$ is $\log$-closed and that each $\iota$ as in (3) is $L(\log)$-preserving, let $z \in \bE_1$.
        By Corollary~\ref{cor:rvpptcor2} and the fact that $\res(\bE_0)=\res(\bE_1)$, there is $z_1 \in \Span_\Lambda (B_0)$, $c\in \bE_0$, and $z_2 \in \bE_1 \cap \co$ such that $z = \exp(z_1) c (1+z_2)$. Thus $\log|z| = z_1 + \log|c|+ \log (1+z_2) \in \bE_1$. Moreover
        \[\begin{aligned}
            \log |\iota(z)| = \log \iota(\exp(z_1)) + \log|c|+ \log (1+\iota(z_2)) =\\
            =\iota(z_1) + \log |c| + \log(1+\iota(z_2))= \iota(\log|z|),
        \end{aligned}\]
        because by the hypothesis on $\iota$, $\iota(\exp(z_1))=\exp(\iota(z_1))$.
    \end{proof}
\end{lemma}

We have said enough that the standard proof of quantifier elimination for $T(\exp, \log)$ is just a few arguments away. The proof of the proposition below is essentially the one given for $\bR_{an, \exp, \log}$ in \cite{dries1994elementary}.

\begin{proposition}\label{prop:qe_for_log_exp}
    $T(\exp,\log)$ is complete, has quantifier elimination and is universally axiomatized.
    \begin{proof}[Proof (Essentially the same argument of \cite{dries1994elementary}).]
        Since by Remark~\ref{rmk:T(exp,log)_extension_axiom_scheme}, the theory $T(\exp, \log)$ is universally axiomatized and has a prime model, it suffices to show that $T(\exp, \log)$ is model complete.
        
        We will do this by showing that given any extension $(\bE', \exp, \log) \subseteq (\bE, \exp, \log)$ of models of $T(\exp, \log)$, there is an elementary extension $\bE_*$ of $\bE'$ and an embedding $\iota: (\bE,\exp, \log) \to (\bE_*, \exp, \log)$ over $\bE'$.

        Let $\bE_*$ be any $|\bE|^+$-saturated extension and endow each $\bE'$, $\bE$ and $\bE_*$ with the $T$-convex valuation ring $\CH(\bP)$.

        It suffices to show that for every fixed $x \in \bE$, we can find an $L(\log, \exp)$-closed $\bE_\omega \subseteq \bE$ such that $\{x\} \cup \bE' \subseteq \bE_\omega$ and that there is an embedding $\iota_\omega: \bE_\omega \to \bE_*$ over $\bE'$.
        
        Observe that if we have an $L(\log)$-substructure $\bE_n\supseteq \bE'$ in $\bE$ and an $L(\log)$-embedding $\iota_n : \bE_n \to \bE_*$ over $\bE'$, then by Lemma~\ref{lem:exp_step} $\bE_{n+1}:= \bE_n \langle \exp(\bE_n)\rangle$ is an $L(\log)$-substructure. If $B_n$ is a $\Lambda$-linear basis of $\bE_n$ over $\log(\bE_n^{>0})$ and we set 
        $\iota_{n+1}(\exp(b)):=\exp(\iota_{n}(b))$ for all $b \in B_n$, we have
        \[(\iota_n)_*\tp_L(\exp(B_n)/\bE_n)=\tp_L(\iota_{n+1}\exp(B_n)/\iota_n\bE_n).\] 
        Thus $\iota_{n+1}\cup \iota_n$ extends to a $L$-embedding $\iota_{n+1}:\bE_{n+1}\to \bE_*$ such that $\iota_{n+1}(\exp(y))=\exp(\iota_{n+1}(y))$ for all $y \in \bE_{n}$. This is then an $L(\log)$-embedding by Lemma~\ref{lem:exp_step}. Thus given an $L(\log)$-substructure $\bE_0 \subseteq \bE$ and an $L(\log)$-embedding $\iota_0: \bE_0 \to \bE_*$ over $\bE'$, by considering the union of the inductively defined $\iota_n: \bE_{n} \to \bE_*$, we have that $\iota_\omega:=\bigcup_{n} \iota_n: \bigcup_n\bE_n \to \bE_*$ is an $L(\log)$-embedding and $\bE_\omega:= \bigcup_n\bE_n$ is $L(\log,\exp)$-closed.
        
        Let us now show that we can always find $\iota_0$ and $\bE_0$ as above and such that moreover $x \in \bE_\omega$. Recall Lemma~\ref{lem:maximal_log_scale} and consider its two cases.
        
        Suppose we are in case (1) of Lemma~\ref{lem:maximal_log_scale}. Thus there is $(y_i)_{i \le n}$ such that $\val(y_i)\notin \val (\bE')$ and $y_{i+1}-\log|y_i|\in \bE'$ for all $i<n$ and $y_n$ is wim or residual. By saturation we can find $z \in \bE_*$ such that $\tp(z/\bE')= \tp(y_n/\bE')$. Thus there is an $L(\bE)$-embedding $\iota: \bE'\langle y_n \rangle \to \bE_*$ such that $\iota(y_n)=z$. Since $\bE'\langle y_n \rangle$ is wim or residual over $\bE'$, by Lemmas~\ref{lem:log(wim)} and \ref{lem:log(res)}, it follows that $\iota : \bE'\langle y_n \rangle \to \bE_*$ is as required.
        Thus in this case we can set $\bE_0 = \bE' \langle y_n\rangle$ and $\iota_0:=\iota$ and we will certainly have $x \in \bE_\omega$.

        Suppose instead that there is no such maximal $(y_i)_{i \le n}$, so we can find $(y_i)_{i<\omega}$ such that $c_i:=\log |y_i|-y_{i+1}\in \bE'$ and $\val(y_i)\notin \val(\bE')$.
        Moreover $\bE_0:= \bE' \langle y_i : i<\omega \rangle_L$ is closed under logarithms.
        Let $z_0\in \bE_*$ be such that $\tp(z_0/\bE')=\tp(y_0/\bE')$ and define inductively $(z_i)_{i<\omega}$ by $z_{i+1}:=\log|z_i|-c_i$.
        
        By Lemma~\ref{lem:log_scale_properties}(3), $\varphi(y_i):=z_i$ defines an $L(\bE)$ elementary map $\varphi : (y_i:i<\omega) \to (z_i:i<\omega)$, thus we can extend $\varphi$ to a $L$-embedding over $\bE'$, $\iota_0: \bE_0 \to \bE \langle z_i:i<\omega\rangle$.
        
        To see that $\iota_0$ is $\log$-preserving observe that any $t \in \bE_0$ can be written as $t=t_0 \cdot (1+t_1) \cdot \prod_{j<n} y_j^{\beta_j}$ for some $n<\omega$, $t_0 \in \bE'$ and $t_1 \in \bE_0\cap \co$. Thus
        \[\begin{aligned}
            \log|\iota_0(t)| = \log |t_0| + \log (1+\iota_0(t_1)) + \sum_{j<n} \beta_j\log |z_j|=\\
            = \log |t_0| + \iota_0\log (1+t_1) + \sum_{j<n} \beta_j (z_{j+1}+c_j) = \iota_0 (\log|t|).
        \end{aligned}\]
        This completes the proof that $T(\exp, \log)$ is model complete.
    \end{proof}
\end{proposition}

\begin{corollary}\label{cor:equivalent_simply_exp}
    If $T'$ is a simply exponential theory, then there is a power-bounded $T$ as in Context~\ref{defn:T_with_restricted_exp} such that $T(\exp, \log)$ and $T'$ have a common expansion by definitions.
    \begin{proof}
        Let $L'$ be the language of $T'$ and $\exp \in L'$.
        Pick a sublanguage $L_0\subseteq L'$ such that $T_0:=T'|_{L_0}$ is power-bounded and defines all restrictions of $\exp$ to bounded $T$-definable intervals, and $T'$ is an expansion by definition of $T'|_{L_0\cup \{\exp\}}$.
        Since o-minimal theories have definable Skolem functions, $T_0$ has an expansion by definitions $T$ in some language $L\supseteq L_0$ that satisfies the assumptions in Context~\ref{defn:T_with_restricted_exp}. Now $T'$ can be expanded by the same definitions for the symbols in $L$ to a theory $T''$ in $L'':=L' \cup L \cup \{\log\}$. Notice that then $T''_{L\cup \{\exp, \log\}}\supseteq T(\exp, \log)$ because since $T'$ is o-minimal, $\exp$ is compatible with every power-bounded reduct of $T$ defining all the restrictions of $\exp$ to bounded $T$-definable intervals. Thus $T''$ is a common expansion by definitions of $T'$ and $T(\exp, \log)$. 
    \end{proof}
\end{corollary}

O-minimality of $T(\exp, \log)$ can also be proven with an analogous argument to the one given in \cite{dries1994elementary}. We won't give the full details because we won't use this Corollary.

\begin{corollary}
    $T(\exp, \log)$ is o-minimal.
    \begin{proof}[Proof Sketch]
        It suffices to show that the prime model $(\bP,\exp, \log)$ is o-minimal. For this in turn it suffices to show that for every $x$ in an elementary extension of $\bP$ and for every definable subset $S$ of $\bP$,
        either $\bP^{<x}\cap S$ and $\bP^{>x}\cap S$ are respectively a final and initial segments of $\bP^{<x}$ and $\bP^{>x}$, or $\bP^{<x}\cap S$ and $\bP^{>x}\cap S$ are respectively bounded above and bounded below in $\bP^{<x}$ and $\bP^{>x}$.
        
        By quantifier elimination, this in turn is equivalent to showing that every $L(\exp,\log)$-term-definable function at $x$ is either positive or negative or $0$ in some interval $(a,b)$ of $\bP$ such that $a<x<b$.
        This follows from compatibility of the exponential whenever $\tp(x/\bP)$ is not definable. So we can assume that $x>\bP$, but then exactly the same arguments of \cite[(5.9) to (5.12)]{dries1994elementary} work. The only extra care needed is related to the use of the intermediate value property for $\exp$ and $\log$ in \cite[(5.10) and (5.11)]{dries1994elementary}: for this notice that both $\exp$ and $\log$ have the intermediate value property on $\bP$, because they are $T$-definable and continuous on every bounded interval of $\bP$.
    \end{proof}
\end{corollary}

\subsection{Weakly immediate types in simply exponential theories}\label{ssec:simply_exp}
		
This section is dedicated to some basic observations on weakly immediate types in a simply exponential theory $T$.
This is meant as a step toward understanding $\lambda$-spherical completions for models of $T_\convex$. In particular it will allow us to show that wim-constructible extensions are 1-wim (Corollary~\ref{cor:wimc_1wim}) which is a step toward answering Question~\ref{quest:wimc-factors}~(1) for these theories.
    
\begin{definition}
	Let $T$ be an o-minimal exponential theory and $(\bE, \cO)\models T_\convex$. A \emph{generalized nested exponential over $\bE$} (\emph{gne} for short) $g: \bE \to \bE$, is an $\bE$-definable composition of translations, changes of signs and exponentials. More specifically:
	\begin{enumerate}
		\item $g(t):=t+c$ is a gne of height $0$ for every $c \in \bE$;
		\item if $h$ is a gne of height $m$, then $g(t):=c + \sigma \exp (h(t))$ for $c \in \bE$ and $\sigma \in \{\pm 1\}$ is a gne of height $m+1$. 
	\end{enumerate}
	Let $(\bU, \cO')\succeq (\bE,\cO)$ and $x \in \bU$, we say inductively that
	\begin{enumerate}
		\item $g$ is \emph{normal} at $x$ if $g(x)=x+c$ with $x\prec c \in \bE$ or $c=0$, or if $g(t)= c + \sigma \exp(h(t))$ with $h$ normal at $x$ and, again, $c=0$ or $c \succ \exp (h(x))$.
		\item $g$ is \emph{essential} at $x$ if $g(t)=t+c$ or if $g(t)= c+\sigma \exp(h(t))$ with $c\in \bE$, $h$ essential at $x$ and $(h(x)+\cO') \cap \bE =\emptyset$.
	\end{enumerate}
\end{definition}

\begin{lemma}\label{lem:gne_ppty}
    Let $T'=T(\exp, \log)$ for some power-bounded theory $T$ in a language $L$. Suppose $T'$ is o-minimal, $(\bE, \cO)\models T'_\convex$, and let $x$ be wim over $(\bE,\cO)$. Then for every $y \in\bE \langle x \rangle$ there is an $\bE$-definable gne $g$ and $z \in \bE \langle x \rangle_{T}$ such that $\tp(y/\bE)= \tp(g(z)/\bE)$ and $g$ is essential at $z$.
    \begin{proof}
        Let $\bE_0:= \bE \langle x \rangle_{T}$ and $\bE_{n+1}:=\bE_n\langle \exp(\bE_n)\rangle_{T}$. 
        Observe that by Lemmas~\ref{lem:log(wim)} and \ref{lem:exp_step}, all $\bE_n$ are closed under logarithms and that by Proposition~\ref{prop:qe_for_log_exp} $\bE\langle x \rangle_{T'} = \bigcup_{n<\omega} \bE_n$.
        
        By induction on $n$ we will show that for all $n$ and $y \in \bE_n$ there is an $\bE$-definable gne $g$ and a $y_0 \in \bE_0$ such that $y \equiv_{\bE} g(y_0)$ and $g$ is essential at $y_0$.

        If $y \in \bE_0$ this is obvious. Suppose the thesis holds for $n$ and let $y \in \bE_{n+1}\setminus \bE_n$.

        By Lemma~\ref{lem:exp_step}~(2) and Corollary~\ref{cor:rvpptcor2}, there is $c \in \bE_n$ such that $\val(y-c) \in \val(\exp(\bE_n))\setminus \val(\bE_n)$. It follows that
        $(\log|y-c|+\cO) \cap \bE_n\neq \emptyset$ and $(\log |y-c|+ \cO) \cap \log \bE_{n}^{>0} = \emptyset$. Thus we can find $\bar{y}\in \bE_n$ such that $(\log|y-c|-\bar{y}) \in \cO$. We will then have that $\bar{y} \notin \cO + \log \bE_n^{>0}$.
        It follows that $\tp_{L}(y/\bE_{n}) = \tp_{L}(c+ \sigma \exp(\bar{y})/\bE_{n})$.

        Since $\val(\exp(\bar{y})) \notin \val (\bE_n)\supseteq \val(\bE)$, it follows that either we can find $\bar{c} \in \bE$ such that $\tp(y/\bE)=\tp(\bar{c} + \sigma \exp(\bar{y})/\bE)$ or $\tp(y/\bE)=\tp(c/\bE)$.

        In both cases the thesis follows from the inductive hypothesis.
    \end{proof}
\end{lemma}

\begin{proposition}\label{prop:gne_ppty}
    Suppose $T$ is simply exponential, $(\bE, \cO)\models T_\convex$, and let $x$ be wim over $(\bE,\cO)$. Then for every $y \in\bE \langle x \rangle$ there is an $\bE$-definable gne $g$ such that $z=g^{-1}(y)$ is wim over $(\bE, \cO)$ and $g$ is essential at $z$.
    \begin{proof}
        By Corollary~\ref{cor:equivalent_simply_exp} we can assume that $T$ is $T_0(\exp, \log)$ for some power-bounded $T_0$ defining a restricted exponential. Then it suffices to apply Lemma~\ref{lem:gne_ppty} and Proposition~\ref{prop:rv-prop}. 
    \end{proof}
\end{proposition}

\begin{lemma}\label{lem:gnechange}
	Let $T$ be exponential and $(\bE, \cO)\prec (\bE_1, \cO_1) \prec (\bU, \cO_*)\models T_\convex$.
    Assume $z\in \bU\setminus \bE_1$ is wim over $(\bE,\cO)$, $g: \bE_1 \to \bE_1$ is an $\bE_1$-definable gne essential at $z$, and that $\tp(g(z)/\bE)$ is not realized in $\bE_1$. Then there is a $z_1\in \bU$ wim over $(\bE,\cO)$ and an $\bE$-definable gne $g_1$ essential at $z_1$ such that $\tp(g_1(z_1)/\bE)=\tp(g(z)/\bE)$.
	\begin{proof}
		By induction on the height of $g$. If $g(t)=t+c$, then since $\tp(z+c/\bE)$ is not realized in $\bE_1$ and $z+c$ is weakly $\cO$-immediate over $\bE_1$, $z_1:=z+c$ is weakly $\cO$-immediate over $\bE$.
			
		Let $g(t) = c + \sigma \exp(h(t))$. By inductive hypothesis, there is an essential $\bE$-definable gne $h_1$, and a $z_1$ weakly $\cO$-immediate over $\bE$ with $\tp(h(z)/\bE) = \tp(h_1(z_1)/\bE)$.
			
		If $\bE \cap (c + \exp(h(z))\co_1)=\emptyset$ then $\tp(g(z)/\bE)$ would be realized in $\bE_1$ by any element in $c+\bE_1^{\prec \exp(h(z))}$ contradicting that $\tp(g(z)/\bE)$ is not realized in $\bE_1$.
			
		Thus we can find $c_1 \in \bE$ such that $c-c_1 \prec \exp h(z)$. But then since $\val(\exp h(z)) \notin \val \bE_1$, $\tp(c+\exp(h(z))/\bE)=\tp(c_1+\exp(h(z))/\bE)=\tp(c_1+\exp(h_1(z_1))/\bE)$, so we can choose $g_1(t):=c_1 + \sigma \exp(h_1(t))$.
 	\end{proof}
\end{lemma}
	
\begin{corollary}
	Let $T$ be exponential and $(\bE, \cO)\prec (\bE_1, \cO_1) \prec (\bE_2, \cO_2)\models T_\convex$. Suppose $(\bE_1, \cO_1)$ is 1-wim over $(\bE,\cO)$ and $(\bE_2, \cO_2)$ is 1-wim over $(\bE_1, \cO_1)$. Then $(\bE_2, \cO_2)$ is 1-wim over $(\bE,\cO)$.
	\begin{proof}
		Let $z \in \bE_2$, so $z$ is weakly immediately generated over $(\bE_1, \cO_1)$.
            
        If $\tp(z/\bE)$ is realized in $\bE_1$, then $z$ is weakly immediately generated over $\bE$ and there is nothing to show. So assume $\tp(z/\bE)$ is not realized in $\bE_1$. 
        Thus by Lemma~\ref{lem:wimstat}, if $\tp(z/\bE_1)$ is weakly $\cO_1$-immediate, then $\tp(z/\bE)$ is weakly $\cO$-immediate.
            
        If instead $\tp(z/\bE_1)$ is not weakly $\cO_1$-immediate, then by Proposition~\ref{prop:gne_ppty}, $z=g(y)$ for some $\bE_1$-definable gne essential at $y$ and $y$ wim over $(\bE_1, \cO_1)$. By Lemma~\ref{lem:gnechange}, there is an $\bE$-definable gne $g_1$ and a $y_1$ wim over $\bE$ such that $g_1$ is essential at $y_1$ and $\tp(g(y)/\bE)=\tp(g_1(y_1)/\bE)$, so $z$ is weakly immediately generated over $\bE$ once again.
	\end{proof}
\end{corollary}
	
\begin{corollary}\label{cor:wimc_1wim}
	If $T$ is simply exponential then every wim-constructible extension of models of $T_\convex$ is 1-wim.
\end{corollary}

\bibliographystyle{abbrvnat}
\bibliography{Res}

@article{dries1994elementary,
 AUTHOR = {{\noopsort{dries}van den Dries}, Lou and Macintyre, Angus and Marker, David},
     TITLE = {The elementary theory of restricted analytic fields with
              exponentiation},
   JOURNAL = {Ann. of Math. (2)},
  FJOURNAL = {Annals of Mathematics. Second Series},
    VOLUME = {140},
      YEAR = {1994},
    NUMBER = {1},
     PAGES = {183--205},
      ISSN = {0003-486X,1939-8980},
   MRCLASS = {12L12 (03C10 03C62)},
  MRNUMBER = {1289495},
MRREVIEWER = {Thanases\ Pheidas},
       DOI = {10.2307/2118545},
       URL = {https://doi.org/10.2307/2118545},
}

@article{wilkie1996model,
 AUTHOR = {Wilkie, A. J.},
     TITLE = {Model completeness results for expansions of the ordered field
              of real numbers by restricted {P}faffian functions and the
              exponential function},
   JOURNAL = {J. Amer. Math. Soc.},
  FJOURNAL = {Journal of the American Mathematical Society},
    VOLUME = {9},
      YEAR = {1996},
    NUMBER = {4},
     PAGES = {1051--1094},
      ISSN = {0894-0347,1088-6834},
   MRCLASS = {03C62 (03C60 03C65 14P15)},
  MRNUMBER = {1398816},
MRREVIEWER = {Luc\ B\'{e}lair},
       DOI = {10.1090/S0894-0347-96-00216-0},
       URL = {https://doi.org/10.1090/S0894-0347-96-00216-0},
}

@article{dries1994real,
 AUTHOR = {{\noopsort{dries}van den Dries}, Lou and Miller, Chris},
     TITLE = {On the real exponential field with restricted analytic
              functions},
   JOURNAL = {Israel J. Math.},
  FJOURNAL = {Israel Journal of Mathematics},
    VOLUME = {85},
      YEAR = {1994},
    NUMBER = {1-3},
     PAGES = {19--56},
      ISSN = {0021-2172,1565-8511},
   MRCLASS = {03C10 (03C62 12J15)},
  MRNUMBER = {1264338},
MRREVIEWER = {A.\ J.\ Wilkie},
       DOI = {10.1007/BF02758635},
       URL = {https://doi.org/10.1007/BF02758635},
}

@article{kaplansky1942maximal,
    AUTHOR = {Kaplansky, Irving},
     TITLE = {Maximal fields with valuations},
   JOURNAL = {Duke Math. J.},
  FJOURNAL = {Duke Mathematical Journal},
    VOLUME = {9},
      YEAR = {1942},
     PAGES = {303--321},
      ISSN = {0012-7094,1547-7398},
   MRCLASS = {09.1X},
  MRNUMBER = {6161},
MRREVIEWER = {Saunders\ Mac Lane},
       URL = {http://projecteuclid.org/euclid.dmj/1077493226},
}

@article {kuhlmann1997exponentiation,
    AUTHOR = {Kuhlmann, Franz-Viktor and Kuhlmann, Salma and Shelah,
              Saharon},
     TITLE = {Exponentiation in power series fields},
   JOURNAL = {Proc. Amer. Math. Soc.},
  FJOURNAL = {Proceedings of the American Mathematical Society},
    VOLUME = {125},
      YEAR = {1997},
    NUMBER = {11},
     PAGES = {3177--3183},
      ISSN = {0002-9939,1088-6826},
   MRCLASS = {12J15 (06F20 12J25 13J05)},
  MRNUMBER = {1402868},
MRREVIEWER = {Joachim\ Gr\"ater},
       DOI = {10.1090/S0002-9939-97-03964-6},
       URL = {https://doi.org/10.1090/S0002-9939-97-03964-6},
}

@article{mourgues1993every,
    AUTHOR = {Mourgues, M.-H. and Ressayre, J. P.},
     TITLE = {Every real closed field has an integer part},
   JOURNAL = {J. Symbolic Logic},
  FJOURNAL = {The Journal of Symbolic Logic},
    VOLUME = {58},
      YEAR = {1993},
    NUMBER = {2},
     PAGES = {641--647},
      ISSN = {0022-4812,1943-5886},
   MRCLASS = {03C60 (12J15 16W80)},
  MRNUMBER = {1233929},
MRREVIEWER = {S.\ R.\ Kogalovski\u{\i}},
       DOI = {10.2307/2275224},
       URL = {https://doi.org/10.2307/2275224},
}

@article {daquino2012real,
    AUTHOR = {D'Aquino, Paola and Knight, Julia F. and Kuhlmann, Salma and
              Lange, Karen},
     TITLE = {Real closed exponential fields},
   JOURNAL = {Fund. Math.},
  FJOURNAL = {Fundamenta Mathematicae},
    VOLUME = {219},
      YEAR = {2012},
    NUMBER = {2},
     PAGES = {163--190},
      ISSN = {0016-2736,1730-6329},
   MRCLASS = {03C57 (03C60 03C70)},
  MRNUMBER = {2993471},
MRREVIEWER = {Giuseppina\ Terzo},
       DOI = {10.4064/fm219-2-6},
       URL = {https://doi.org/10.4064/fm219-2-6},
}

@incollection {ressayre1993integer,
    AUTHOR = {Ressayre, J.-P.},
     TITLE = {Integer parts of real closed exponential fields (extended
              abstract)},
 BOOKTITLE = {Arithmetic, proof theory, and computational complexity
              ({P}rague, 1991)},
    SERIES = {Oxford Logic Guides},
    VOLUME = {23},
     PAGES = {278--288},
 PUBLISHER = {Oxford Univ. Press, New York},
      YEAR = {1993},
      ISBN = {0-19-853690-9},
   MRCLASS = {03C60 (12J15)},
  MRNUMBER = {1236467},
}

@article {kuhlmann2005kappa,
    AUTHOR = {Kuhlmann, Salma and Shelah, Saharon},
     TITLE = {{$\kappa$}-bounded exponential-logarithmic power series
              fields},
   JOURNAL = {Ann. Pure Appl. Logic},
  FJOURNAL = {Annals of Pure and Applied Logic},
    VOLUME = {136},
      YEAR = {2005},
    NUMBER = {3},
     PAGES = {284--296},
      ISSN = {0168-0072,1873-2461},
   MRCLASS = {03C60 (06A05)},
  MRNUMBER = {2169687},
MRREVIEWER = {Martin\ Weese},
       DOI = {10.1016/j.apal.2005.04.001},
       URL = {https://doi.org/10.1016/j.apal.2005.04.001},
}

@book{kuhlmann2000ordered,
AUTHOR = {Kuhlmann, Salma},
     TITLE = {Ordered exponential fields},
    SERIES = {Fields Institute Monographs},
    VOLUME = {12},
 PUBLISHER = {American Mathematical Society, Providence, RI},
      YEAR = {2000},
     PAGES = {xviii+166},
      ISBN = {0-8218-0943-1},
   MRCLASS = {12J15 (03C60 06F25 12L12 26A12)},
  MRNUMBER = {1760173},
MRREVIEWER = {R.\ H.\ Redfield},
       DOI = {10.1090/fim/012},
       URL = {https://doi.org/10.1090/fim/012},
}

@article {macpherson2000weakly,
    AUTHOR = {Macpherson, Dugald and Marker, David and Steinhorn, Charles},
     TITLE = {Weakly o-minimal structures and real closed fields},
   JOURNAL = {Trans. Amer. Math. Soc.},
  FJOURNAL = {Transactions of the American Mathematical Society},
    VOLUME = {352},
      YEAR = {2000},
    NUMBER = {12},
     PAGES = {5435--5483},
      ISSN = {0002-9947,1088-6850},
   MRCLASS = {03C64 (03C60)},
  MRNUMBER = {1781273},
MRREVIEWER = {Carlo\ Toffalori},
       DOI = {10.1090/S0002-9947-00-02633-7},
       URL = {https://doi.org/10.1090/S0002-9947-00-02633-7},
}

@article{dries1997t,
AUTHOR = {{\noopsort{dries}van den Dries}, Lou},
     TITLE = {{$T$}-convexity and tame extensions. {II}},
   JOURNAL = {J. Symbolic Logic},
  FJOURNAL = {The Journal of Symbolic Logic},
    VOLUME = {62},
      YEAR = {1997},
    NUMBER = {1},
     PAGES = {14--34},
      ISSN = {0022-4812,1943-5886},
   MRCLASS = {03C60 (03C10 03C35 12L12)},
  MRNUMBER = {1450511},
MRREVIEWER = {O.\ V.\ Belegradek},
       DOI = {10.2307/2275729},
       URL = {https://doi.org/10.2307/2275729},
}

@article{dries1995t,
 AUTHOR = {{\noopsort{dries}van den Dries}, Lou and Lewenberg, Adam H.},
     TITLE = {{$T$}-convexity and tame extensions},
   JOURNAL = {J. Symbolic Logic},
  FJOURNAL = {The Journal of Symbolic Logic},
    VOLUME = {60},
      YEAR = {1995},
    NUMBER = {1},
     PAGES = {74--102},
      ISSN = {0022-4812,1943-5886},
   MRCLASS = {03C60 (03C10 03C35 12L12)},
  MRNUMBER = {1324502},
MRREVIEWER = {O.\ V.\ Belegradek},
       DOI = {10.2307/2275510},
       URL = {https://doi.org/10.2307/2275510},
}

@incollection {miller1993growth,
    AUTHOR = {Miller, Chris},
     TITLE = {A growth dichotomy for o-minimal expansions of ordered fields},
 BOOKTITLE = {Logic: from foundations to applications ({S}taffordshire,
              1993)},
    SERIES = {Oxford Sci. Publ.},
     PAGES = {385--399},
 PUBLISHER = {Oxford Univ. Press, New York},
      YEAR = {1996},
      ISBN = {0-19-853862-6},
   MRCLASS = {03C60 (03C50)},
  MRNUMBER = {1428013},
MRREVIEWER = {G.\ Cherlin},
}

@book{tyne2003t,
 AUTHOR = {Tyne, James Michael},
     TITLE = {T-levels and {T}-convexity},
      NOTE = {Thesis (Ph.D.)--University of Illinois at Urbana-Champaign},
 PUBLISHER = {ProQuest LLC, Ann Arbor, MI},
      YEAR = {2003},
     PAGES = {106},
      ISBN = {978-0496-34011-8},
   MRCLASS = {99-05},
  MRNUMBER = {2704495},
       IgnoreURL = {http://gateway.proquest.com/openurl?url_ver=Z39.88-2004&rft_val_fmt=info:ofi/fmt:kev:mtx:dissertation&res_dat=xri:pqdiss&rft_dat=xri:pqdiss:3086205},
}

@article{dries2000field,
 AUTHOR = {{\noopsort{dries}van den Dries}, Lou and Speissegger, Patrick},
     TITLE = {The field of reals with multisummable series and the
              exponential function},
   JOURNAL = {Proc. London Math. Soc. (3)},
  FJOURNAL = {Proceedings of the London Mathematical Society. Third Series},
    VOLUME = {81},
      YEAR = {2000},
    NUMBER = {3},
     PAGES = {513--565},
      ISSN = {0024-6115,1460-244X},
   MRCLASS = {03C64 (03C10 12L12 26E05)},
  MRNUMBER = {1781147},
MRREVIEWER = {Chris\ Miller},
       DOI = {10.1112/S0024611500012648},
       URL = {https://doi.org/10.1112/S0024611500012648},
}

@article{dries2002minimal,
  AUTHOR = {{\noopsort{dries}van den Dries}, L. and Speissegger, P.},
     TITLE = {O-minimal preparation theorems},
 BOOKTITLE = {Model theory and applications},
    SERIES = {Quad. Mat.},
    VOLUME = {11},
     PAGES = {87--116},
 PUBLISHER = {Aracne, Rome},
      YEAR = {2002},
      ISBN = {88-7999-411-5},
   MRCLASS = {03C64 (03C10 32B05)},
  MRNUMBER = {2159715},
MRREVIEWER = {Marcus\ Tressl},
}

@article{pillay1994definability,
AUTHOR = {Pillay, Anand},
     TITLE = {Definability of types, and pairs of {O}-minimal structures},
   JOURNAL = {J. Symbolic Logic},
  FJOURNAL = {The Journal of Symbolic Logic},
    VOLUME = {59},
      YEAR = {1994},
    NUMBER = {4},
     PAGES = {1400--1409},
      ISSN = {0022-4812,1943-5886},
   MRCLASS = {03C45 (03C40 03C50)},
  MRNUMBER = {1312317},
MRREVIEWER = {O.\ V.\ Belegradek},
       DOI = {10.2307/2275712},
       URL = {https://doi.org/10.2307/2275712},
}

@article{marker1994definable,
AUTHOR = {Marker, David and Steinhorn, Charles I.},
     TITLE = {Definable types in {O}-minimal theories},
   JOURNAL = {J. Symbolic Logic},
  FJOURNAL = {The Journal of Symbolic Logic},
    VOLUME = {59},
      YEAR = {1994},
    NUMBER = {1},
     PAGES = {185--198},
      ISSN = {0022-4812,1943-5886},
   MRCLASS = {03C45},
  MRNUMBER = {1264974},
MRREVIEWER = {Alexandre\ Ivanov},
       DOI = {10.2307/2275260},
       URL = {https://doi.org/10.2307/2275260},
}

@article{kaplan2023t,
    AUTHOR = {Kaplan, Elliot},
     TITLE = {{$T$}-convex {$T$}-differential fields and their immediate
              extensions},
   JOURNAL = {Pacific J. Math.},
  FJOURNAL = {Pacific Journal of Mathematics},
    VOLUME = {320},
      YEAR = {2022},
    NUMBER = {2},
     PAGES = {261--298},
      ISSN = {0030-8730,1945-5844},
   MRCLASS = {03C64 (12H05 12J10)},
  MRNUMBER = {4548312},
MRREVIEWER = {Athipat\ Thamrongthanyalak},
       DOI = {10.2140/pjm.2022.320.261},
       URL = {https://doi.org/10.2140/pjm.2022.320.261},
}

@misc{freni2024structure,
      title={On the structure of {$T$}-{$\lambda$}-spherical completions of o-minimal fields}, 
      author={Pietro Freni},
      year={2024},
      eprint={2407.07442},
      archivePrefix={arXiv},
      primaryClass={math.LO},
      url={https://arxiv.org/abs/2407.07442}, 
}

@article {mennuni2020product,
    AUTHOR = {Mennuni, Rosario},
     TITLE = {Product of invariant types modulo domination-equivalence},
   JOURNAL = {Arch. Math. Logic},
  FJOURNAL = {Archive for Mathematical Logic},
    VOLUME = {59},
      YEAR = {2020},
    NUMBER = {1-2},
     PAGES = {1--29},
      ISSN = {0933-5846,1432-0665},
   MRCLASS = {03C45},
  MRNUMBER = {4050032},
MRREVIEWER = {Amador\ Martin-Pizarro},
       DOI = {10.1007/s00153-019-00676-9},
       URL = {https://doi.org/10.1007/s00153-019-00676-9},
}

@article {tressl2006pseudo,
    AUTHOR = {Tressl, Marcus},
     TITLE = {Pseudo completions and completions in stages of o-minimal
              structures},
   JOURNAL = {Arch. Math. Logic},
  FJOURNAL = {Archive for Mathematical Logic},
    VOLUME = {45},
      YEAR = {2006},
    NUMBER = {8},
     PAGES = {983--1009},
      ISSN = {0933-5846,1432-0665},
   MRCLASS = {03C64 (12J10 12J15)},
  MRNUMBER = {2271334},
MRREVIEWER = {Chris\ Miller},
       DOI = {10.1007/s00153-006-0022-2},
       URL = {https://doi.org/10.1007/s00153-006-0022-2},
}

@article {dickmann1987elimination,
    AUTHOR = {Dickmann, M. A.},
     TITLE = {Elimination of quantifiers for ordered valuation rings},
   JOURNAL = {J. Symbolic Logic},
  FJOURNAL = {The Journal of Symbolic Logic},
    VOLUME = {52},
      YEAR = {1987},
    NUMBER = {1},
     PAGES = {116--128},
      ISSN = {0022-4812,1943-5886},
   MRCLASS = {03C10 (03C60 12J20 12L12 13L05)},
  MRNUMBER = {877859},
MRREVIEWER = {J.\ M.\ Plotkin},
       DOI = {10.2307/2273866},
       URL = {https://doi.org/10.2307/2273866},
}

@unpublished{wreo36168,
           month = {September},
            note = {Unpublished},
          school = {University of Leeds},
           title = {Structural Investigations in some Classes of o-minimal Fields},
            year = {2024},
             url = {https://etheses.whiterose.ac.uk/id/eprint/36168/},
          author = {Freni, Pietro}
}

@article {dries1998correction,
    AUTHOR = {{\noopsort{dries}van den Dries}, Lou},
     TITLE = {Correction to: ``{$T$}-convexity and tame extensions. {II}''
              [{J}. {S}ymbolic {L}ogic {\bf 62} (1997), no. 1, 14--34;
              {MR}1450511 (98h:03048)]},
   JOURNAL = {J. Symbolic Logic},
  FJOURNAL = {The Journal of Symbolic Logic},
    VOLUME = {63},
      YEAR = {1998},
    NUMBER = {4},
     PAGES = {1597},
      ISSN = {0022-4812,1943-5886},
   MRCLASS = {03C60 (03C10 03C35 12L12)},
  MRNUMBER = {1665787},
       DOI = {10.2307/2586669},
       URL = {https://doi.org/10.2307/2586669},
}

@article {aschenbrenner2018maximal,
    AUTHOR = {Aschenbrenner, Matthias and van den Dries, Lou and van der
              Hoeven, Joris},
     TITLE = {Maximal immediate extensions of valued differential fields},
   JOURNAL = {Proc. Lond. Math. Soc. (3)},
  FJOURNAL = {Proceedings of the London Mathematical Society. Third Series},
    VOLUME = {117},
      YEAR = {2018},
    NUMBER = {2},
     PAGES = {376--406},
      ISSN = {0024-6115,1460-244X},
   MRCLASS = {12H05 (16W60)},
  MRNUMBER = {3851327},
MRREVIEWER = {Akira\ Masuoka},
       DOI = {10.1112/plms.12128},
       URL = {https://doi.org/10.1112/plms.12128},
}

@article {berarducci2023exponential,
    AUTHOR = {Berarducci, Alessandro and Kuhlmann, Salma and Mantova,
              Vincenzo and Matusinski, Micka\"el},
     TITLE = {Exponential fields and {C}onway's omega-map},
   JOURNAL = {Proc. Amer. Math. Soc.},
  FJOURNAL = {Proceedings of the American Mathematical Society},
    VOLUME = {151},
      YEAR = {2023},
    NUMBER = {6},
     PAGES = {2655--2669},
      ISSN = {0002-9939,1088-6826},
   MRCLASS = {03C64 (16W60)},
  MRNUMBER = {4576327},
MRREVIEWER = {Ricardo\ Bianconi},
       DOI = {10.1090/proc/14577},
       URL = {https://doi.org/10.1090/proc/14577},
}

@article {dries2001fields,
    AUTHOR = {{\noopsort{dries}van den Dries}, Lou and Ehrlich, Philip},
     TITLE = {Fields of surreal numbers and exponentiation},
   JOURNAL = {Fund. Math.},
  FJOURNAL = {Fundamenta Mathematicae},
    VOLUME = {167},
      YEAR = {2001},
    NUMBER = {2},
     PAGES = {173--188},
      ISSN = {0016-2736,1730-6329},
   MRCLASS = {03C64 (12J99)},
  MRNUMBER = {1816044},
       DOI = {10.4064/fm167-2-3},
       URL = {https://doi.org/10.4064/fm167-2-3},
}

@article {dries2001erratum,
    AUTHOR = {{\noopsort{dries}van den Dries}, Lou and Ehrlich, Philip},
     TITLE = {Erratum to: ``{F}ields of surreal numbers and
              exponentiation''},
   JOURNAL = {Fund. Math.},
  FJOURNAL = {Fundamenta Mathematicae},
    VOLUME = {168},
      YEAR = {2001},
    NUMBER = {3},
     PAGES = {295--297},
      ISSN = {0016-2736,1730-6329},
   MRCLASS = {03C64 (12J99)},
  MRNUMBER = {1853411},
       DOI = {10.4064/fm168-3-5},
       URL = {https://doi.org/10.4064/fm168-3-5},
}

@article {bradley-williams2023spherically,
    AUTHOR = {Bradley-Williams, David B. and Halupczok, Immanuel},
     TITLE = {Spherically complete models of {H}ensel minimal valued fields},
   JOURNAL = {MLQ Math. Log. Q.},
  FJOURNAL = {MLQ. Mathematical Logic Quarterly},
    VOLUME = {69},
      YEAR = {2023},
    NUMBER = {2},
     PAGES = {138--146},
      ISSN = {0942-5616,1521-3870},
   MRCLASS = {03C64 (03C60)},
  MRNUMBER = {4635177},
MRREVIEWER = {Mihai\ D.\ Prunescu},
}

@phdthesis{foster2010power,
  edition = {},
  number = {},
  journal = {},
  pages = {},
  publisher = {Oxford University, UK},
  school = {Oxford University, UK},
  title = {Power functions and exponentials in o-minimal expansions of fields},
  volume = {},
  author = {Foster, T},
  editor = {},
  year = {2010},
  series = {}
}

@inproceedings{parusinski1994lipschitz,
  title={Lipschitz stratification of subanalytic sets},
  author={Parusi{\'n}ski, Adam},
  booktitle={Annales scientifiques de l'Ecole normale sup{\'e}rieure},
  volume={27},
  number={6},
  pages={661--696},
  year={1994}
}

@article{nowak2007proof,
author = {Krzysztof Jan Nowak},
journal = {Annales Polonici Mathematici},
language = {eng},
number = {1},
pages = {75-85},
title = {A proof of the valuation property and preparation theorem},
url = {http://eudml.org/doc/281119},
volume = {92},
year = {2007},
}

@inproceedings{lion1997theoreme,
  title={Th{\'e}oreme de pr{\'e}paration pour les fonctions logarithmico-exponentielles},
  author={Lion, Jean-Marie and Rolin, Jean-Philippe},
  booktitle={Annales de l'institut Fourier},
  volume={47},
  number={3},
  pages={859--884},
  year={1997}
}

@misc{freni2024t2,
      title={$T$-convexly valued o-minimal fields are definably spherically complete}, 
      author={Pietro Freni},
      year={2024},
      eprint={2411.16706},
      archivePrefix={arXiv},
      primaryClass={math.LO},
      url={https://arxiv.org/abs/2411.16706}, 
}

@article {tressl2005model,
    AUTHOR = {Tressl, Marcus},
     TITLE = {Model completeness of o-minimal structures expanded by
              {D}edekind cuts},
   JOURNAL = {J. Symbolic Logic},
  FJOURNAL = {The Journal of Symbolic Logic},
    VOLUME = {70},
      YEAR = {2005},
    NUMBER = {1},
     PAGES = {29--60},
      ISSN = {0022-4812,1943-5886},
   MRCLASS = {03C64 (03C50)},
  MRNUMBER = {2119122},
MRREVIEWER = {Patrick\ Speissegger},
       DOI = {10.2178/jsl/1107298509},
       URL = {https://doi.org/10.2178/jsl/1107298509},
}

@book {tent2012course,
    AUTHOR = {Tent, Katrin and Ziegler, Martin},
     TITLE = {A course in model theory},
    SERIES = {Lecture Notes in Logic},
    VOLUME = {40},
 PUBLISHER = {Association for Symbolic Logic, La Jolla, CA; Cambridge
              University Press, Cambridge},
      YEAR = {2012},
     PAGES = {x+248},
      ISBN = {978-0-521-76324-0},
   MRCLASS = {03C35 (03-01 03C45)},
  MRNUMBER = {2908005},
MRREVIEWER = {David\ Evans},
       DOI = {10.1017/CBO9781139015417},
       URL = {https://doi.org/10.1017/CBO9781139015417},
}

@misc{berarducci2026elementary,
      title={On the elementary theory of the real exponential field}, 
      author={Alessandro Berarducci and Francesco Gallinaro},
      year={2026},
      eprint={2603.08365},
      archivePrefix={arXiv},
      primaryClass={math.LO},
      url={https://arxiv.org/abs/2603.08365}, 
}

@misc{daquino2015note,
      title={A note on $\aleph_{\alpha}$-saturated o-minimal expansions of real closed fields}, 
      author={Paola D'Aquino and Salma Kuhlmann},
      year={2015},
      eprint={1112.4078},
      archivePrefix={arXiv},
      primaryClass={math.LO},
      url={https://arxiv.org/abs/1112.4078}, 
}

\end{document}